\documentclass[english,11pt,aps,letter,superscriptaddress,nofootinbib,pre]{revtex4}
\pdfoutput=1
\usepackage[T1]{fontenc}
\usepackage[latin9]{inputenc}
\usepackage{xcolor}
\usepackage{pdfcolmk}
\usepackage{amsmath}
\usepackage{amssymb}
\usepackage{graphicx}
\usepackage{esint}
\PassOptionsToPackage{normalem}{ulem}
\usepackage{ulem}

\makeatletter

\providecolor{lyxadded}{rgb}{1,0,0}
\providecolor{lyxdeleted}{rgb}{0,0,1}

\@ifundefined{textcolor}{}
{%
 \definecolor{BLACK}{gray}{0}
 \definecolor{WHITE}{gray}{1}
 \definecolor{RED}{rgb}{1,0,0}
 \definecolor{GREEN}{rgb}{0,1,0}
 \definecolor{BLUE}{rgb}{0,0,1}
 \definecolor{CYAN}{cmyk}{1,0,0,0}
 \definecolor{MAGENTA}{cmyk}{0,1,0,0}
 \definecolor{YELLOW}{cmyk}{0,0,1,0}
}

\usepackage{ae,aecompl}

\newcommand{\sfrac}[2]{\mathchoice
  {\kern0em\raise.5ex\hbox{\the\scriptfont0 #1}\kern-.15em/
   \kern-.15em\lower.25ex\hbox{\the\scriptfont0 #2}}
  {\kern0em\raise.5ex\hbox{\the\scriptfont0 #1}\kern-.15em/
   \kern-.15em\lower.25ex\hbox{\the\scriptfont0 #2}}
  {\kern0em\raise.5ex\hbox{\the\scriptscriptfont0 #1}\kern-.2em/
   \kern-.15em\lower.25ex\hbox{\the\scriptscriptfont0 #2}}
  {#1\!/#2}}

\makeatletter
\DeclareMathSizes{\@xipt}{10}{6}{5}
\makeatother

\makeatother

\usepackage{babel}
\begin{document}
\global\long\def\V#1{\boldsymbol{#1}}
\global\long\def\M#1{\boldsymbol{#1}}
\global\long\def\Set#1{\mathbb{#1}}

\global\long\def\D#1{\Delta#1}
\global\long\def\d#1{\delta#1}

\global\long\def\norm#1{\left\Vert #1\right\Vert }
\global\long\def\abs#1{\left|#1\right|}

\global\long\def\grad{\M{\nabla}}
\global\long\def\avv#1{\langle#1\rangle}
\global\long\def\av#1{\left\langle #1\right\rangle }

\global\long\def\myhalf{\sfrac{1}{2}}
\global\long\def\mythreehalves{\sfrac{3}{2}}

\title{Low Mach Number Fluctuating Hydrodynamics of Diffusively Mixing Fluids}

\author{Aleksandar Donev}

\email{donev@courant.nyu.edu}

\selectlanguage{english}%

\affiliation{Courant Institute of Mathematical Sciences, New York University,
New York, NY 10012}

\author{Andy Nonaka}

\affiliation{Center for Computational Science and Engineering, Lawrence Berkeley
National Laboratory, Berkeley, CA, 94720}

\author{Yifei Sun}

\affiliation{Courant Institute of Mathematical Sciences, New York University,
New York, NY 10012}

\affiliation{Leon H. Charney Division of Cardiology, Department of Medicine, New
York University School of Medicine, New York, NY 10016}

\author{Thomas G. Fai}

\affiliation{Courant Institute of Mathematical Sciences, New York University,
New York, NY 10012}

\author{Alejandro L. Garcia }

\affiliation{Department of Physics and Astronomy, San Jose State University, San
Jose, California, 95192}

\author{John B. Bell}

\affiliation{Center for Computational Science and Engineering, Lawrence Berkeley
National Laboratory, Berkeley, CA, 94720}
\begin{abstract}
We formulate low Mach number fluctuating hydrodynamic equations appropriate
for modeling diffusive mixing in isothermal mixtures of fluids with
different density and transport coefficients. These equations represent
a coarse-graining of the microscopic dynamics of the fluid molecules
in both space and time, and eliminate the fluctuations in pressure
associated with the propagation of sound waves by replacing the equation
of state with a local thermodynamic constraint. We demonstrate that
the low Mach number model preserves the spatio-temporal spectrum of
the slower diffusive fluctuations. We develop a strictly conservative
finite-volume spatial discretization of the low Mach number fluctuating
equations in both two and three dimensions and construct several explicit
Runge-Kutta temporal integrators that strictly maintain the equation
of state constraint. The resulting spatio-temporal discretization
is second-order accurate deterministically and maintains fluctuation-dissipation
balance in the linearized stochastic equations. We apply our algorithms
to model the development of giant concentration fluctuations in the
presence of concentration gradients, and investigate the validity
of common simplifications such as neglecting the spatial non-homogeneity
of density and transport properties. We perform simulations of diffusive
mixing of two fluids of different densities in two dimensions and
compare the results of low Mach number continuum simulations to hard-disk
molecular dynamics simulations. Excellent agreement is observed between
the particle and continuum simulations of giant fluctuations during
time-dependent diffusive mixing.
\end{abstract}
\maketitle
\setlength{\abovedisplayskip}{0.35ex}\setlength{\belowdisplayskip}{0.35ex}

\section{Introduction}

Stochastic fluctuations are intrinsic to fluid dynamics because fluids
are composed of molecules whose positions and velocities are random
at thermodynamic scales. Because they span the whole range of scales
from the microscopic to the macroscopic \cite{DiffusionRenormalization_PRL,FractalDiffusion_Microgravity},
fluctuations need to be consistently included in all levels of description.
Stochastic effects are important for flows in new microfluidic, nanofluidic
and microelectromechanical devices \cite{Nanofluidics_Review}; novel
materials such as nanofluids \cite{Nanofluids_Review}; biological
systems such as lipid membranes \cite{SIBM_Biomembrane}, Brownian
molecular motors \cite{BrownainMotor_Peskin}, nanopores \cite{Nanopore_Fluctuations};
as well as processes where the effect of fluctuations is amplified
by strong non-equilibrium effects, such as ultra clean combustion,
capillary dynamics \cite{CapillaryNanowaves,StagerredFluct_Inhomogeneous},
and hydrodynamic instabilities \cite{BreakupNanojets,DropletSpreading,DropFormationFluctuations}.

One can capture thermal fluctuations using direct particle level calculations.
But even coarse-grained particle methods \cite{ParticleMesoscaleHydrodynamics,SHSD_PRL,DiffusionRenormalization_PRL}
are computationally expensive because the dynamics of individual particles
has time scales significantly shorter than hydrodynamic time scales.
Alternatively, thermal fluctuations can be included in the Navier-Stokes
equations through stochastic forcing terms, as proposed by Landau
and Lifshitz \cite{Landau:Fluid} and later extended to fluid mixtures
\cite{FluctHydroNonEq_Book}. The basic idea of \emph{fluctuating
hydrodynamics} is to add a \emph{stochastic flux} corresponding to
each dissipative (irreversible, diffusive) flux \cite{OttingerBook}.
This ensures that the microscopic conservation laws and thermodynamic
principles are obeyed while also maintaining fluctuation-dissipation
balance. Specifically, the equilibrium thermal fluctuations have the
Gibbs-Boltzmann distribution dictated by statistical mechanics. Fluctuating
hydrodynamics is a useful tool in understanding complex fluid flows
far from equilibrium \cite{FluctHydroNonEq_Book} but theoretical
calculations are often only feasible after ignoring nonlinearities,
inhomogeneities in density, temperature, and transport properties,
surface dynamics, gravity, unsteady flow patterns, and other important
effects. In the past decade fluctuating hydrodynamics has been applied
to study a number of nontrivial practical problems \cite{LB_SoftMatter_Review,StagerredFluct_Inhomogeneous,LB_OrderParameters,SELM};
however, the numerical methods used are far from the comparable state-of-the-art
for deterministic solvers.

Previous computational studies of the effect of thermal fluctuations
in fluid mixtures \cite{StagerredFluct_Inhomogeneous,Bell:09,LB_OrderParameters}
have been based on the compressible fluid equations and thus require
small time steps to resolve fast sound waves (pressure fluctuations).
Recently, some of us developed finite-volume methods for the incompressible
equations of fluctuating hydrodynamics \cite{LLNS_Staggered}, which
eliminate the stiffness arising from the separation of scales between
the acoustic and vortical modes \cite{IncompressibleLimit_Majda,ZeroMachCombustion}.
For inhomogeneous fluids with non-constant density, diffusive mass
and heat fluxes create local expansion and contraction of the fluid
and the incompressibility constraint should be replaced by a ``quasi-incompressibility''
constraint \cite{ZeroMachCombustion,Cahn-Hilliard_QuasiIncomp}. The
resulting \emph{low-Mach number} equations have been used for some
time to model deterministic flows with thermo-chemical effects \cite{ZeroMach_Buoyancy,ZeroMachCombustion},
and several conservative finite-volume techniques have been developed
for solving equations of this type \cite{LowMachAdaptive,ZeroMach_Klein,LaminarFlowChemistry,LowMach_FiniteDifference,LowMachAcoustics}.
To our knowledge, thermal fluctuations have not yet been incorporated
in low Mach number models.

In this work we extend the staggered-grid, finite-volume approach
developed in Ref. \cite{LLNS_Staggered} to isothermal mixtures of
fluids with unequal densities. The imposition of the quasi-incompressibility
constraint poses several nontrivial mathematical and computational
challenges. At the mathematical level, the traditional low Mach number
asymptotic expansions \cite{IncompressibleLimit_Majda,ZeroMachCombustion}
assume spatio-temporal smoothness of the flow and thus do not directly
apply in the stochastic context. At the computational level, enforcing
the quasi-incompressibility or equation of state (EOS) constraint
in a conservative and stable manner requires specialized spatio-temporal
discretizations. By careful selection of the analytical form of the
EOS constraint and the spatial discretization of the advective fluxes
we are able to maintain strict local conservation and enforce the
EOS to within numerical tolerances. In the present work, we employ
an explicit projection-based temporal discretizations because of the
substantial complexity of designing and implementing semi-implicit
discretizations of the momentum equation for spatially-inhomogeneous
fluids \cite{StokesKrylov}.

Thermal fluctuations exhibit unusual features in systems out of thermodynamic
equilibrium. Notably, external gradients can lead to \emph{enhancement}
of thermal fluctuations and to \emph{long-range} correlations between
fluctuations \cite{DSMC_Fluctuations_Shear,Mareschal:92,LongRangeCorrelations_MD,Zarate:04,FluctHydroNonEq_Book}.
Sharp concentration gradients present during diffusive mixing lead
to the development of macroscopic or \emph{giant fluctuations }\cite{GiantFluctuations_Theory,TemperatureGradient_Cannell,GiantFluctuations_ThinFilms}
in concentration, which have been observed using light scattering
and shadowgraphy techniques \cite{GiantFluctuations_Nature,GiantFluctuations_Cannell,FractalDiffusion_Microgravity}.
These experimental studies have found good but imperfect agreement
between the predictions of a simplified fluctuating hydrodynamic theory
and experiments. Computer simulations are, in principle, an ideal
tool for studying such complex time-dependent processes in the presence
of nontrivial boundary conditions without making the sort of approximations
necessary for analytical calculations, such as assuming spatially-constant
density and transport coefficients and spatially-uniform gradients.
On the other hand, the multiscale (more precisely, \emph{many-scale})
nature of the equations of fluctuating hydrodynamics poses many mathematical
and computational challenges that are yet to be addressed. Notably,
it is necessary to develop temporal integrators that can accurately
and robustly handle the large separation of time scales between different
physical processes, such as mass and momentum diffusion. The computational
techniques we develop here form the foundation for incorporating additional
physics, such as heat transfer and internal energy fluctuations, phase
separation and interfacial dynamics, and chemical reactions.

We begin Section \ref{sec:Equations} by formulating the fluctuating
low Mach number equations for an isothermal binary fluid mixture.
We present both a traditional pressure (constrained) formulation and
a gauge (unconstrained) formulation. We analyze the spatio-temporal
spectrum of the thermal fluctuations in the linearized equations and
demonstrate that the low Mach equations eliminate the fast (sonic)
pressure fluctuations but maintain the correct spectrum of the slow
(diffusive) fluctuations. In Section \ref{sec:TemporalIntegration}
we develop projected Runge-Kutta schemes for solving the spatially-discretized
equations, including a midpoint and a trapezoidal second-order predictor-corrector
scheme, and a third-order three-stage scheme. In Section \ref{sec:SpatialDiscretization}
we describe a spatial discretization of the equations that strictly
maintains the equation of state constraint and also obeys a fluctuation-dissipation
balance principle \cite{LLNS_S_k}. In Section \ref{sec:GiantFluct}
we study the steady-state spectrum of giant concentration fluctuations
in the presence of an applied concentration gradient in a mixture
of two dissimilar fluids, and test the applicability of common approximations
that neglect spatial inhomogeneities. In Section \ref{sec:MixingMD}
we study the dynamical evolution of giant interface fluctuations during
diffusive mixing of two dissimilar fluids, using both hard-disk molecular
dynamics and low Mach number fluctuating hydrodynamics. We find excellent
agreement between the two, providing a strong support for the usefulness
of the fluctuating low Mach number equations as a coarse-grained model
of complex fluid mixtures. In Section \ref{sec:Conclusions} we offer
some concluding remarks and point out several outstanding challenges
for the future. Several technical calculations and procedures are
detailed in Appendices.

\section{\label{sec:Equations}Low Mach Number Equations}

The compressible equations of fluctuating hydrodynamics were proposed
some time ago \cite{Landau:Fluid} and have since been studied and
applied successfully to a variety of situations \cite{FluctHydroNonEq_Book}.
The presence of rapid pressure fluctuations due to the propagation
of sound waves leads to stiffness that makes it computationally expensive
to solve the fully compressible equations numerically, especially
for typical liquids. It is therefore important to develop fluctuating
hydrodynamics equations that capture the essential physics in cases
where acoustics can be neglected. 

It is important to note that the equations of fluctuating hydrodynamics
are to be interpreted as a mesoscopic coarse-grained representation
of the mass, momentum and energy transport which occurs at microscopic
scales through molecular interactions (collisions). As such, these
equations implicitly contain a mesoscopic coarse-graining length and
time scale that is larger than molecular scales \cite{DiscreteLLNS_Espanol}.
While a coarse-graining scale does not appear explicitly in the formal
stochastic partial differential equations (SPDEs) written in this
section (but note that it can be if desired \cite{DiffusionJSTAT}),
it does explicitly enter in the spatio-temporal discretization described
in Section \ref{sec:SpatialDiscretization} through the grid spacing
(equivalently, the volume of the grid, or more precisely, the number
of molecules per grid cell) and time step size. This changes the appropriate
interpretation of convergence of numerical methods to a continuum
limit in the presence of fluctuations and nonlinearities \cite{DFDB}.
Only for the linearized equations of fluctuating hydrodynamics \cite{FluctHydroNonEq_Book}
can the formal SPDEs be given a precise continuum meaning \cite{LLNS_S_k}.

Developing coarse-grained models that only resolve the relevant spatio-temporal
scales is a well-studied but still \emph{ad hoc} procedure that requires
substantial \emph{a priori} physical insight \cite{OttingerBook}.
More precise mathematical mode-elimination procedures \cite{AdiabaticElimination_1,MoriZwanzig_ConstrainedMD}
are technically involved and often purely formal, especially in the
context of SPDEs \cite{DiffusionJSTAT}. Here we follow a heuristic
approach to constructing fluctuating low Mach number equations, starting
from the well-known deterministic low Mach equations (which can be
obtained via asymptotic analysis \cite{IncompressibleLimit_Majda,ZeroMachCombustion})
and then adding fluctuations in a manner consistent with fluctuation-dissipation
balance. Alternatively, our low Mach number equations can be seen
as a formal asymptotic limit in which the noise terms are formally
treated as smooth forcing terms; a more rigorous derivation is nontrivial
and is deferred for future work.

\subsection{Compressible Equations}

The starting point of our investigations is the system of isothermal
compressible equations of fluctuating hydrodynamics for the density
$\rho(\V r,t)$, velocity $\V v(\V r,t)$, and mass concentration
$c(\V r,t)$ for a mixture of two fluids in $d$ dimensions. In terms
of mass and momentum densities the equations can be written as conservation
laws \cite{OttingerBook,FluctHydroNonEq_Book,Bell:09}, 
\begin{align}
\partial_{t}\rho+\grad\cdot\left(\rho\V v\right)= & \;0\nonumber \\
\partial_{t}\left(\rho\V v\right)+\grad\cdot\left(\rho\V v\V v^{T}\right)= & -\grad P+\grad\cdot\left[\eta\left(\grad\V v+\grad^{T}\V v\right)+\left(\kappa-\frac{2}{d}\eta\right)\left(\grad\cdot\V v\right)\M I+\M{\Sigma}\right]+\rho\V g\nonumber \\
\partial_{t}\left(\rho_{1}\right)+\grad\cdot\left(\rho_{1}\V v\right)= & \grad\cdot\left[\rho\chi\left(\grad c+K_{P}\grad P\right)+\M{\Psi}\right],\label{LLNS_primitive}
\end{align}
where $\rho_{1}=\rho c$ is the density of the first component, $\rho_{2}=(1-c)\rho$
is the density of the second component, $P(\rho,c;T)$ is the equation
of state for the pressure at the reference temperature $T=T_{0}=\mbox{const.}$,
and $\V g$ is the gravitational acceleration. Temperature fluctuations
are neglected in this study but can be accounted for using a similar
approach. The shear viscosity $\eta$, bulk viscosity $\kappa$, mass
diffusion coefficient $\chi$, and baro-diffusion coefficient $K_{P}$,
in general, depend on the state. The baro-diffusion coefficient $K_{P}$
above {[}denoted with $k_{P}/P$ in Ref. \cite{Bell:09}, see Eq.
(A.17) in that paper{]} is not a transport coefficient but rather
determined from thermodynamics \cite{Landau:StatPhys1},
\begin{equation}
K_{P}=\frac{\left(\partial\mu/\partial P\right)_{c}}{\left(\partial\mu/\partial c\right)_{P}}=-\rho^{-2}\frac{\left(\partial\rho/\partial c\right)_{P}}{\left(\partial\mu/\partial c\right)_{P}}=\frac{\left(\partial P/\partial c\right)_{\rho}}{\rho^{2}c_{T}^{2}\mu_{c}},\label{eq:K_P}
\end{equation}
where $\mu$ is the chemical potential of the mixture at the reference
temperature, $\mu_{c}=\left(\partial\mu/\partial c\right)_{P}$, and
$c_{T}^{2}=\left(\partial P/\partial\rho\right)_{c}$ is the isothermal
speed of sound. The capital Greek letters denote stochastic momentum
and mass fluxes that are formally modeled as \cite{LLNS_Staggered}
\begin{align}
\M{\Sigma}=\sqrt{\eta k_{B}T}\left(\M{\mathcal{W}}+\M{\mathcal{W}}^{T}-\frac{2}{d}\mathrm{Tr}\,\M{\mathcal{W}}\right)+\sqrt{\frac{2\kappa k_{B}T}{d}}\,\mathrm{Tr}\,\M{\mathcal{W}}\mbox{ and } & \M{\Psi}=\sqrt{2\chi\rho\mu_{c}^{-1}k_{B}T}\;\widetilde{\M{\mathcal{W}}},\label{stoch_flux_covariance}
\end{align}
where $k_{B}$ is Boltzmann's constant, and $\M{\mathcal{W}}(\V r,t)$
and $\widetilde{\M{\mathcal{W}}}(\V r,t)$ are standard zero mean,
unit variance random Gaussian tensor and vector fields with uncorrelated
components,
\[
\avv{\mathcal{W}_{ij}(\V r,t)\mathcal{W}_{kl}(\V r^{\prime},t')}=\delta_{ik}\delta_{jl}\;\delta(t-t^{\prime})\delta(\V r-\V r^{\prime}),
\]
and similarly for $\widetilde{\M{\mathcal{W}}}$.

\subsection{Low Mach Equations}

At mesoscopic scales, in typical liquids, sound waves are much faster
than momentum diffusion and can usually be eliminated from the fluid
dynamics description. Formally, this corresponds to taking the zero
Mach number singular limit $c_{T}\rightarrow\infty$ of the system
(\ref{LLNS_primitive}) by performing an asymptotic analysis as the
Mach number $\text{Ma}=U/c_{T}\rightarrow0$, where $U$ is a reference
flow velocity. The limiting dynamics can be obtained by performing
an asymptotic expansion in the Mach number \cite{IncompressibleLimit_Majda}.
In a deterministic setting this analysis shows that the pressure can
be written in the form
\[
P(\V r,t)=P_{0}(t)+\pi(\V r,t)
\]
where $\pi=O\left(\text{Ma}^{2}\right)$. The low Mach number equations
can then be obtained by making the anzatz that the thermodynamic behavior
of the system is captured by the reference pressure, $P_{0}$, and
$\pi$ captures the mechanical behavior while not affecting the thermodynamics.
We note that when the system is sufficiently large or the gravitational
forcing is sufficiently strong, assuming a spatial constant reference
pressure is not valid. In those cases, the reference pressure represents
a global hydrostatic balance, $\nabla P_{0}=\rho_{0}\V g$ (see \cite{AnelasticApproximation}
for details of the construction of these types of models). Here, however,
we will restrict consideration to cases where gravity causes negligible
changes in the thermodynamic state across the domain.

In this case, the reference pressure constrains the system so that
the evolution of $\rho$ and $c$ remains consistent with the thermodynamic
equation of state
\begin{equation}
P\left(\rho\left(\V r,t\right),c\left(\V r,t\right);T\right)=P_{0}\left(t\right).\label{eq:general_EOS}
\end{equation}
This constraint means that any change in concentration (equivalently,
$\rho_{1}$) must be accompanied by a corresponding change in density,
as would be observed in a system at thermodynamic equilibrium held
at the fixed reference pressure and temperature. This implies that
variations in density are coupled to variations in composition. Note
that we do not account for temperature variations in our isothermal
model.

The equation for $\rho_{1}$ can be written in primitive (non-conservation)
form as the concentration equation
\begin{equation}
\rho\frac{Dc}{Dt}=\rho D_{t}c=\rho\left(\partial_{t}c+\V v\cdot\grad c\right)=\grad\cdot\V F,\label{eq:div_v_derivation}
\end{equation}
where the non-advective (diffusive and stochastic) fluxes are denoted
with
\[
\V F=\rho\chi\grad c+\M{\Psi}.
\]
Note that there is no barodiffusion flux because barodiffusion is
of thermodynamic origin (as seen from (\ref{eq:K_P}) \cite{FluctHydroNonEq_Book})
and involves the gradient of the \emph{thermodynamic} pressure $\nabla P_{0}=0$.
By differentiating the EOS constraint along a Lagrangian trajectory
we obtain 
\begin{equation}
\frac{D\rho}{Dt}=\beta\rho\frac{Dc}{Dt}=\beta\grad\cdot\V F=\partial_{t}\rho+\V v\cdot\grad\rho=-\rho\nabla\cdot\V v,\label{eq:Drho_Dt}
\end{equation}
where the solutal expansion coefficient 
\[
\beta\left(c\right)=\frac{1}{\rho}\left(\frac{\partial\rho}{\partial c}\right)_{P_{0}}
\]
is determined by the specific form of the EOS. 

Equation (\ref{eq:Drho_Dt}) shows that the EOS constraint can be
re-written as a constraint on the divergence of velocity,
\begin{equation}
\rho\grad\cdot\V v=-\beta\,\grad\cdot\V F.\label{eq:div_v}
\end{equation}
Note that the usual incompressibility constraint is obtained when
the density is not affected by changes in concentration, $\beta=0$.
When $\beta\neq0$ changes in composition (concentration) due to diffusion
cause local expansion and contraction of the fluid and thus a nonzero
$\grad\cdot\V v$. It is important at this point to consider the boundary
conditions. For a closed system, such as a periodic domain or a system
with rigid boundaries, we must ensure that the integral of $\grad\cdot\V v$
over the domain is zero. This is consistent with (\ref{eq:div_v})
if $\beta/\rho$ is constant, so that we can rewrite (\ref{eq:div_v})
in the form $\grad\cdot\V v=-\grad\cdot\left(\left(\beta/\rho\right)\V F\right)$.
In this case $P_{0}$ does not vary in time. If $\beta/\rho$ is not
constant, then for a closed system the reference pressure $P_{0}$
must vary in time to enforce that the total fluid volume remains constant.
Here we will assume that $\beta/\rho=\text{const.}$, and we will
give a specific example of an EOS that obeys this condition.

The asymptotic low Mach analysis of (\ref{LLNS_primitive}) is standard
and follows the procedure outlined in Ref. \cite{IncompressibleLimit_Majda},
formally treating the stochastic forcing as smooth. This analysis
leads to the \emph{isothermal low Mach number} equations for a binary
mixture of fluids in conservation form,
\begin{align}
\partial_{t}\left(\rho\V v\right)+\nabla\pi=-\grad\cdot\left(\rho\V v\V v^{T}\right) & +\grad\cdot\left[\eta\left(\grad\V v+\grad^{T}\V v\right)+\M{\Sigma}\right]+\rho\V g\equiv\V f(\rho,\V v,c,t)\label{eq:momentum_eq}\\
\partial_{t}\left(\rho_{1}\right)=-\grad\cdot\left(\rho_{1}\V v\right)+ & \grad\cdot\V F\equiv h(\rho,\V v,c,t)\label{eq:rho1_eq}\\
\partial_{t}\left(\rho_{2}\right)=-\grad\cdot\left(\rho_{2}\V v\right) & -\grad\cdot\V F\label{eq:rho2_eq}\\
\mbox{such that }\grad\cdot\V v= & -\left(\rho^{-1}\beta\right)\,\grad\cdot\V F\equiv S(\rho,c,t).\label{eq:div_v_constraint}
\end{align}
The gradient of the non-thermodynamic component of the pressure $\pi$
(Lagrange multiplier) appears in the momentum equation as a driving
force that ensures the EOS constraint (\ref{eq:div_v_constraint})
is obeyed. We note that the bulk viscosity term gives a gradient term
that can be absorbed in $\pi$ and therefore does not explicitly need
to appear in the equations. By adding the two density equations (\ref{eq:rho1_eq},\ref{eq:rho2_eq})
we get the usual continuity equation for the total density,
\begin{equation}
\partial_{t}\rho=-\grad\cdot\left(\rho\V v\right)\label{eq:rho_eq}
\end{equation}
Our conservative numerical scheme is based on Eqs. (\ref{eq:momentum_eq},\ref{eq:rho1_eq},\ref{eq:div_v_constraint},\ref{eq:rho_eq}).

In Appendix \ref{sec:LinearizedAnalysis}, we apply the standard linearized
fluctuating hydrodynamics analysis to the low Mach number equations.
This gives expressions for the equilibrium and nonequilibrium static
and dynamic covariances (spectra) of the fluctuations in density and
concentration as a function of wavenumber and wavefrequency. Specifically,
the dynamic structure factor in the low Mach number approximation
has the form
\[
S_{\rho,\rho}\left(\V k,\omega\right)=\av{\left(\widehat{\d{\rho}}\right)\left(\widehat{\d{\rho}}\right)^{\star}}=\beta^{2}\left(\rho\mu_{c}^{-1}k_{B}T\right)\frac{2\chi k^{2}}{\omega^{2}+\chi^{2}k^{4}}.
\]
The linearized analysis shows that the low Mach number equations reproduce
the slow fluctuations (small $\omega$) in density and concentration
(central Rayleigh peak in the dynamic structure factor \cite{FluctHydroNonEq_Book,LLNS_S_k})
as in the full compressible equations (see Section \ref{sub:CompressibleSpectra}),
while eliminating the fast isentropic pressure fluctuations (side
Brillouin peaks) from the dynamics.

The fluctuations in velocity, however, are different between the compressible
and low Mach number equations. In the compressible equations, the
dynamic structure factor for the longitudinal component of velocity
decays to zero as $\omega\rightarrow\infty$ because it has two sound
(Brillouin) peaks centered around $\omega\approx\pm c_{T}k$, in addition
to the central diffusive (Rayleigh) peak. The low Mach number equations
reproduce the central peak (slow fluctuations) correctly, replacing
the side peaks with a flat spectrum for large $\omega$, which is
unphysical as it formally makes the velocity white in time. The low
Mach equations should therefore be used only for time scales larger
than the sound propagation time.

The fact that the velocity fluctuations are white in space and in
time poses a further challenge in interpreting the nonlinear low Mach
number equations, and in particular, numerical schemes may not converge
to a sensible limit as the time step goes to zero. In practice, just
as the spatial discretization of the equations imposes a spatial smoothing
or regularization of the fluctuations, the temporal discretization
of the equations imposes a temporal smoothing and filters the problematic
large frequencies. In the types of problems we study in this work
the problem concentration fluctuations can be neglected, $\hat{\M{\Psi}}\approx\V 0$,
because the concentration fluctuations are dominated by nonequilibrium
effects. If $\hat{\M{\Psi}}=\V 0$ the problematic white-in-time longitudinal
component of velocity disappears.

\subsubsection{Model Equation of State}

In general, the EOS constraint (\ref{eq:general_EOS}) is a non-linear
constraint. In this work we consider a specific linear EOS, 
\begin{equation}
\frac{\rho_{1}}{\bar{\rho}_{1}}+\frac{\rho_{2}}{\bar{\rho}_{2}}=\frac{c\rho}{\bar{\rho}_{1}}+\frac{(1-c)\rho}{\bar{\rho}_{2}}=1,\label{eq:EOS_quasi_incomp}
\end{equation}
where $\bar{\rho}_{1}$ and $\bar{\rho}_{2}$ are the densities of
the pure component fluids ($c=1$ and $c=0$, respectively), giving
\begin{equation}
\beta=\rho\left(\frac{1}{\bar{\rho}_{2}}-\frac{1}{\bar{\rho}_{1}}\right)=\frac{\bar{\rho}_{1}-\bar{\rho}_{2}}{c\bar{\rho}_{2}+(1-c)\bar{\rho}_{1}}.\label{eq:beta_simple}
\end{equation}
It is important that for this specific form of the EOS $\beta/\rho$
is a material constant independent of the concentration. The density
dependence (\ref{eq:beta_simple}) on concentration arises if one
assumes that the two fluids do not change volume upon mixing. This
is a reasonable assumption for liquids that are not too dissimilar
at the molecular level. Surprisingly the EOS (\ref{eq:EOS_quasi_incomp})
is also valid for a mixture of ideal gases, since
\[
P=P_{1}+P_{2}=P_{0}=nk_{B}T=\left(n_{1}+n_{2}\right)k_{B}T=\left(\frac{\rho_{1}}{m_{1}}+\frac{\rho_{2}}{m_{2}}\right)k_{B}T,
\]
where $m$ is molecular mass and $n=\rho/m$ is the number density.
This is exactly of the form (\ref{eq:EOS_quasi_incomp}) with $\bar{\rho}_{1}=m_{1}P_{0}/\left(k_{B}T\right)=nm_{1}$
and $\bar{\rho}_{2}=nm_{2}$.

Even if the specific EOS (\ref{eq:EOS_quasi_incomp}) is not a very
good approximation over the entire range of concentration $0\leq c\leq1$,
(\ref{eq:EOS_quasi_incomp}) may be a very good approximation over
the range of concentrations of interest if $\bar{\rho}_{1}$ and $\bar{\rho}_{2}$
are adjusted accordingly. In this case $\bar{\rho}_{1}$ and $\bar{\rho}_{2}$
are not the densities of the pure component fluids but rather fitting
parameters that approximate the true EOS in the range of concentrations
of interest. For small variations in concentration around some reference
concentration $\bar{c}$ and density $\bar{\rho}$ one can approximate
$\beta\approx\bar{\rho}^{-1}\left(\partial\rho/\partial c\right)_{\bar{c}}$
by a constant and determine appropriate values of $\bar{\rho}_{1}$
and $\bar{\rho}_{2}$ from (\ref{eq:beta_simple}) and the EOS (\ref{eq:EOS_quasi_incomp})
evaluated at the reference state. Our specific form choice of the
EOS will aid significantly in the construction of simple conservative
spatial discretizations that strictly maintain the EOS without requiring
complicated nonlinear iterative corrections.

\subsubsection{Boundary Conditions}

Several different types of boundary conditions can be imposed for
the low Mach number equations, just as for the more familiar incompressible
equations. The simplest case is when periodic boundary conditions
are used for all of the variables. We briefly describe the different
types of conditions that can be imposed at a physical boundary with
normal direction $n$.

For the concentration (equivalently, $\rho_{1}$), either Neumann
(zero mass flux) or Dirichlet (fixed concentration) boundary conditions
can be imposed. Physically, a Neumann condition corresponds to a physical
boundary that is impermeable to mass, while Dirichlet conditions correspond
to a permeable membrane that connects the system to a large reservoir
held at a specified concentration. In the case of Neumann conditions
for concentration, both the normal component of the diffusive flux
$F_{n}=0$ and the advective flux $\rho_{1}v_{n}=0$ vanish at the
boundary, implying that the normal component of velocity must vanish,
$v_{n}=0$. For Dirichlet conditions on the concentration, however,
there will, in general, be a nonzero normal diffusive flux $F_{n}$
through the boundary. This diffusive flux for concentration will induce
a corresponding mass flux, as required to maintain the equation of
state near the boundary. From the condition (\ref{eq:div_v_constraint}),
we infer the proper boundary condition for the normal component of
velocity to be
\begin{equation}
v_{n}=-\left(\rho^{-1}\beta\right)F_{n}.\label{eq:v_n_BC}
\end{equation}
This condition expresses the notion that there is no net volume change
for the fluid in the domain. Note that no additional boundary conditions
can be specified for $\rho$ since its boundary conditions follow
from those on $c$ via the EOS constraint.

For the tangential component of velocity $\V v_{\tau}$, we either
impose a no-slip condition $\V v_{\tau}=0$, or a free slip boundary
condition in which the tangential component of the normal viscous
stress vanishes,
\[
\eta\left(\frac{\partial v_{n}}{\partial\V{\tau}}+\frac{\partial\V v_{\tau}}{\partial n}\right)=\V 0.
\]
In the case of zero normal mass flux, $v_{n}=0$, the free slip condition
simplifies to a Neumann condition for the tangential velocity, $\partial\V v_{\tau}/\partial n=0$.

\subsection{Gauge Formalism}

The low Mach number system of equations (\ref{eq:momentum_eq},\ref{eq:rho1_eq},\ref{eq:div_v_constraint},\ref{eq:rho_eq})
is a \emph{constrained} problem. For the purposes of analysis and
in particular for constructing higher-order temporal integrators,
it is useful to rewrite the constrained low Mach number equations
as an \emph{unconstrained} initial value problem. In the incompressible
case, $\grad\cdot\V v=0$, we can write the constrained Navier-Stokes
equations as an unconstrained system by eliminating the pressure using
a projection operator formalism. The constraint $\grad\cdot\V v=0$
is a constant linear constraint and independent of the state and of
time. However, in the low Mach number equations the velocity-divergence
constraint $\grad\cdot\V v=-\beta D_{t}c$ depends on concentration,
and also on time when there are additional (stochastic or deterministic)
forcing terms in the concentration equation. Treating this type of
system requires a more general vector field decomposition. This more
general vector field decomposition provides the basis for a projection-based
discretization of the constrained system. We also introduce a gauge
formulation of the system \cite{GaugeIncompressible_E} that casts
the evolution as a nonlocal unconstrained system that is analytically
equivalent to the orignal constrained evolution. The gauge formulation
allows us to develop higher-order method-of-lines temporal integration
algorithms.

\subsubsection{Vector Field Decomposition}

The velocity in the low Mach number equations can be split into two
components,
\[
\V v=\V u+\grad\zeta,
\]
where $\grad\cdot\V u=0$ is a divergence-free (solenoidal or vortical)
component, and therefore
\[
\grad\cdot\V v=\grad^{2}\zeta=S(\rho,c,t).
\]
This is a Poisson problem for $\zeta$ that is well-posed for appropriate
boundary conditions on $\V v$. Specifically, periodic boundary conditions
on $\V v$ imply periodic boundary conditions for $\V u$ and $\zeta$.
At physical boundaries where a Dirichlet condition (\ref{eq:v_n_BC})
is specified for the normal component of the velocity, we set $u_{n}=0$
and use Neumann conditions for the Poisson solve, $\partial\zeta/\partial n=v_{n}$.

We can now define a more general vector field decomposition that plays
the role of the Hodge decomposition in incompressible flow. Given
a vector field $\tilde{\V v}$ and a density $\rho$ we can decompose
$\tilde{\V v}$ into three components 
\[
\tilde{\V v}=\V u+\grad\zeta+\rho^{-1}\grad\psi.
\]
This decomposition can be obtained by using the condition $\grad\cdot\V u=0$
and $\grad^{2}\zeta=S$, which allows us to define a density-weighted
Poisson equation for $\psi$, 
\[
\grad\cdot\left(\rho^{-1}\grad\psi\right)=-\grad\cdot\left(\tilde{\V v}-\grad\zeta\right)=-\grad\cdot\tilde{\V v}+S(\rho,c,t).
\]
Let $\M L_{\rho}^{-1}$ denote the solution operator to the density-dependent
Poisson problem, formally,
\[
\M L_{\rho}^{-1}=\left[\grad\cdot\left(\rho^{-1}\grad\right)\right]^{-1},
\]
and also define a density-dependent projection operator $\M{\mathcal{P}}_{\rho}$
defined through its action on a vector field $\V w$,
\[
\M{\mathcal{P}}_{\rho}\V w=\V w-\rho^{-1}\grad\left[\M L_{\rho}^{-1}\left(\grad\cdot\V w\right)\right].
\]
This is a well-known variable density generalization \cite{almgren-iamr}
of the constant-density projection operator $\M{\mathcal{P}}\V w=\V w-\grad\left[\grad^{-2}\left(\grad\cdot\V w\right)\right]$.
We can now write
\[
\V u=\M{\mathcal{P}}_{\rho}\left(\tilde{\V v}-\grad\zeta\right)=\M{\mathcal{P}}_{\rho}\tilde{\V v}+\rho^{-1}\grad\left[\M L_{\rho}^{-1}S(\rho,c,t)\right]-\grad\zeta.
\]
This gives
\[
\V v=\V u+\grad\zeta=\M{\mathcal{R}}_{S}\left(\tilde{\V v}\right),
\]
where we have introduced an affine transformation $\M{\mathcal{R}}_{S}(\rho,c,t)$
that depends on $\rho$, $c$ and $t$ through $S(\rho,c,t)$, and
is defined via its action on a vector field $\V w$, 
\begin{equation}
\M{\mathcal{R}}_{S}\left(\V w\right)=\V w-\rho^{-1}\grad\left[\M L_{\rho}^{-1}\left(\grad\cdot\V w-S\right)\right].\label{P_tilde_v}
\end{equation}
Note that application of $\M{\mathcal{R}}_{S}$ requires only one
Poisson solve and does not actually require computing $\zeta$.

\subsubsection{Gauge Formulation}

The low Mach number system (\ref{eq:momentum_eq},\ref{eq:rho1_eq},\ref{eq:rho_eq},\ref{eq:div_v_constraint})
has the form
\begin{eqnarray}
\partial_{t}\rho & = & -\grad\cdot\left(\rho\V v\right)\nonumber \\
\partial_{t}\V m+\grad\pi & = & \V f\left(c,\V v,t\right)\nonumber \\
\partial_{t}\rho_{1} & = & h(c,\V v,t)\nonumber \\
\grad\cdot\V v & = & S(\rho,c,t),\label{eq:div_v_S}
\end{eqnarray}
where $\V m=\rho\V v$ is the momentum density, and $\V f$, $h$
and $S$ are as defined in (\ref{eq:momentum_eq},\ref{eq:rho1_eq},\ref{eq:div_v_constraint}).
At present, we will assume that these functions are smooth functions
of time, which is only justified in the presence of stochastic forcing
terms in a linearized setting. We note that, for the constrained system,
$\rho$ is not an independent variable because of the EOS constraint
(\ref{eq:EOS_quasi_incomp}); however, we will retain the evolution
of $\rho$ with the implicit understanding that the evolution must
be constrained so that $\rho$ and $c$ remain consistent with (\ref{eq:EOS_quasi_incomp}).

To define the gauge formulation, we introduce a new variable
\[
\tilde{\V m}=\rho\tilde{\V v}=\V m+\grad\psi,
\]
where $\psi$ is a \emph{gauge} variable. We note that $\psi$ is
not uniquely determined; however, the specific choice does not matter.
If we choose the gauge so that $\partial_{t}\psi=\pi$ then the momentum
equation in (\ref{eq:div_v_S}) is equivalent to
\[
\partial_{t}\tilde{\V m}=\V f(\rho,\V v,c,t).
\]
The appropriate boundary conditions for $\psi$ are linked to the
boundary conditions on $\V v$; we set $\psi$ to be periodic if $\V v$
is periodic, and employ a homogeneous Neumann (natural) boundary condition
$\partial\psi/\partial n=0$ if a Dirichlet condition (\ref{eq:v_n_BC})
is specified for the normal component of the velocity $v_{n}$. Note
that in the spatially-discrete staggered formulation that we employ,
the homogeneous Neumann condition follows automatically from the boundary
conditions on velocity used to define the appropriate divergence and
gradient operators in the interior of the domain.

If we know $\tilde{\V m}$ and $\rho$, we can then define $\tilde{\V v}=\tilde{\V m}/\rho$
and compute $\V v=\M{\mathcal{R}}_{S}\left(\tilde{\V v}\right)$,
where $\M{\mathcal{R}}_{S}$ is defined in (\ref{P_tilde_v}). Thus
by using the gauge formulation we can formally write the low Mach
number equations in the form of an unconstrained initial value problem
\begin{eqnarray}
\partial_{t}\tilde{\V m} & = & \V f\left(\rho(c),\M{\mathcal{R}}_{S}\left(\tilde{\V v}\right),c,t\right)\label{eq:m_gauge}\\
\partial_{t}\rho_{1} & = & h\left(\rho(c),\M{\mathcal{R}}_{S}\left(\tilde{\V v}\right),c,t\right).\label{eq:rho1_gauge}
\end{eqnarray}
The utility of the gauge formulation is that in fact, we do not need
to know $\psi$ in order to determine $\V v$. Therefore, the time
evolution equation for $\psi$ does not actually need to be solved,
and in particular, $\pi$ does not need to be computed. Futhermore,
by adopting the gauge formulation, we can directly use a method of
lines approach for spatially-discretizing the system (\ref{eq:m_gauge},\ref{eq:rho1_gauge}),
and then apply standard Runge-Kutta temporal integrators to the resulting
system of ordinary (stochastic) differential equations.

It is important to emphasize that the actual independent physical
variables in the low Mach formulation (\ref{eq:m_gauge},\ref{eq:rho1_gauge})
are the vortical (solenoidal) component of velocity $\V u$ and the
concentration $c$. The density $\rho=\rho(c)$ and the velocity $\V v=\V u+\grad\left[\grad^{-2}S(\rho,c,t)\right]$
are determined from $\V u$ and $c$ and the constraints; hence they
can formally be eliminated from the system, as can be seen in the
linearized analysis in Appendix \ref{sec:LinearizedAnalysis}, which
shows that fluctuations in the vortical velocity modes are decoupled
from the longitudinal fluctuations.

\section{\label{sec:TemporalIntegration}Temporal Integration}

Our spatio-temporal discretization follows a ``method of lines''
approach in which we first discretize the equations (\ref{eq:momentum_eq},\ref{eq:rho1_eq},\ref{eq:div_v_constraint},\ref{eq:rho_eq})
in space and then integrate the resulting semi-continuum equations
in time. Our uniform staggered-grid spatial discretization of the
low Mach number equations is relatively standard and is described
in Section \ref{sec:SpatialDiscretization}. The main difficulty is
the temporal integration of the resulting equations in the presence
of the EOS constraint. Our temporal integrators are based on the gauge
formulation (\ref{eq:m_gauge},\ref{eq:rho1_gauge}) of the low Mach
equations. The gauge formulation is unconstrained and enables us to
use standard temporal integrators for initial-value problems. In the
majority of this section, we assume that all of the fields and differential
operators have already been spatially discretized and focus on the
temporal integration of the resulting initial-value problem.

Because in the present schemes we handle both diffusive and advective
fluxes explicitly, the time step size $\D t$ is restricted by well-known
CFL conditions. For fluctuating hydrodynamics applications the time
step is typically limited by momentum diffusion,
\[
\alpha_{\nu}=\frac{\nu\D t}{\D x^{2}}<\frac{1}{2d},
\]
where $d$ is the number of spatial dimensions and $\D x$ is the
grid spacing. The design and implementation of numerical methods that
handle momentum diffusion semi-implicitly, as done in Ref. \cite{LLNS_Staggered}
for incompressible flow, is substantially more difficult for the low
Mach number equations because it requires a variable coefficient implicit
fluid solver. We have recently developed an efficient Stokes solver
for solving variable-density and variable-viscosity time-dependent
and steady Stokes problems \cite{StokesKrylov}, and in future work
we will employ this solver to construct a semi-implicit temporal integrator
for the low Mach number equations.

Our temporal discretization will make use of the special form of the
EOS and the discretization of mass advection described in Section
\ref{sub:Advection} in order to strictly maintain the EOS relation
(\ref{eq:EOS_quasi_incomp}) between density and concentration in
each cell at \emph{all} intermediate values. Therefore, no additional
action is needed to enforce the EOS constraint after an update of
$\rho_{1}$ and $\rho$. This is, however, only true to within the
accuracy of the Poisson solver and also roundoff, and it is possible
for a slow drifting off the EOS to occur over many time steps. In
Section \ref{sub:DriftCorrection}, we describe a correction that
prevents such drifting and ensures that the EOS is obeyed at all times
to essentially roundoff tolerance. For simplicity, we will often omit
the explicit update for the density $\rho$ and instead focus on updating
$\rho_{1}$ and the momentum density $\V m=\rho\V v$, with the understanding
that $\rho$ is updated whenever $\rho_{1}$ is.

\subsection{Euler Scheme}

The foundation for our higher-order explicit temporal integrators
is the first-order Euler method applied to the gauge formulation (\ref{eq:m_gauge},\ref{eq:rho1_gauge}).

\subsubsection{Gauge-Free Euler Update}

We use a superscript to denote the time step and the point in time
where a given term is evaluated, e.g., $f^{n}\equiv f_{D}\left(\rho^{n},\V v^{n},c^{n},t^{n}\right)$
where $f_{D}$ denotes the spatial discretization of $f$ with analogous
definitions for $h^{n}$ and $S^{n}$. We also denote the time step
size with $\D t=t^{n+1}-t^{n}$. Assume that at the beginning of timestep
$n$ we know $\tilde{\V m}^{n}$ and we can then compute 
\[
\V v^{n}=\M{\mathcal{R}}_{S}^{n}\left(\tilde{\V v}^{n}\right)
\]
 by enforcing the constraint (\ref{eq:div_v_S}). Here $\M{\mathcal{R}}_{S}^{n}$
denotes the affine transformation (\ref{P_tilde_v}) with all terms
evaluated at the beginning of the time step, so that $\grad\cdot\V v^{n}=S^{n}$.
An Euler step for the low Mach equations then consists of the update
\begin{eqnarray}
\rho_{1}^{n+1} & = & \rho_{1}^{n}+\D t\, h^{n}\nonumber \\
\tilde{\V m}^{n+1} & = & \tilde{\V m}^{n}+\D t\,\V f^{n},\label{eq:Euler_lagged}
\end{eqnarray}
together with an update of the density $\rho^{n+1}$ consistent with
$\rho_{1}^{n+1}$.

At the beginning of the next time step, $\V v^{n+1}$ will be calculated
from $\tilde{\V m}^{n+1}$ by applying $\M{\mathcal{R}}_{S}^{n+1}$,
and it is only $\V v^{n+1}$ that will actually be used during time
step $n+1$. We therefore do not need to explicitly store $\tilde{\V m}^{n+1}$
and can instead replace it with $\M m^{n+1}=\rho^{n+1}\V v^{n+1}$
without changing any of the observable results. This is related to
the fact that the gauge is \emph{de facto} arbitrary and, in the present
setting, the gauge formulation is simply a formalism to put the equations
in an unconstrained form suitable for method of lines discretization.
The difference between $\tilde{\V m}$ and $\M m$ is a (discrete)
gradient of a scalar. Since our temporal integrators only use linear
combinations of the intermediate values, the difference between the
final result for $\tilde{\V m}^{n+1}$ and $\M m^{n}$ is also a gradient
of a scalar and replacing $\tilde{\V m}^{n+1}$ with $\M m^{n+1}$
simply amounts to redefining the (arbitrary) gauge variable. For these
reasons, the Euler advance,
\begin{eqnarray}
\rho_{1}^{n+1} & = & \rho_{1}^{n}+\D t\, h^{n}\nonumber \\
\M m^{n+1} & = & \rho^{n+1}\,\M{\mathcal{R}}_{S}^{n+1}\left[\left(\rho^{n+1}\right)^{-1}\left(\V m^{n}+\D t\,\V f^{n}\right)\right],\label{eq:Euler_step_complex}
\end{eqnarray}
is analytically equivalent to (\ref{eq:Euler_lagged}). We will use
this form as the foundation for our temporal integrators. The equivalence
to the gauge form implies that the update specified by (\ref{eq:Euler_step_complex})
can be viewed as an explicit update in spite of the formal dependence
of the update on the solution at both old and new time levels.

\subsubsection{Stochastic Forcing}

Thermal fluctuations cannot be straightforwardly incorporated in (\ref{eq:Euler_step_complex})
because it is not clear how to define $\M{\mathcal{R}}_{S}^{n+1}$.
In the deterministic setting, $S$ is a function of concentration
and density and can be evaluated pointwise at time level $n+1$. When
the white-in-time stochastic concentration flux $\M{\Psi}$ is included,
however, $S$ cannot be evaluated at a particular point of time. Instead,
one must think of $\M{\Psi}$ as representing the \emph{average} stochastic
flux over a given time interval $\d t$, which can be expressed in
terms of the increments $\sqrt{\d t}\,\widetilde{\V W}$ of the underlying
Wiener processes,
\[
\M{\Psi}\left(\d t,\,\widetilde{\V W}\right)=\sqrt{\frac{2\chi\rho\mu_{c}^{-1}k_{B}T}{\d t\,\D V}}\,\widetilde{\V W},
\]
where $\widetilde{\V W}$ is a collection of normal variates generated
using a pseudo-random number generator, and $\D V$ is the volume
of the hydrodynamic cells. Similarly, the average stochastic momentum
flux over a time step is modeled as
\[
\M{\Sigma}\left(\d t,\,\M W\right)=\sqrt{\frac{\eta k_{B}T}{\d t\,\D V}}\,\left(\M W+\M W^{T}\right),
\]
where $\M W$ are normal random variates. As described in more detail
in Ref. \cite{LLNS_Staggered}, stochastic fluxes are spatially discretized
by generating normal variates on the faces of the grid on which the
corresponding variable is discretized, independently at each time
step. As mentioned earlier, the volume of the grid cell appears here
because it expresses the spatial coarse graining length scale (i.e.,
the degree of coarse-graining for which a fluid element with discrete
molecules can be modeled by continuous density fields) implicit in
the equations of fluctuating hydrodynamics. Similarly, the time interval
$\d t\sim\D t$ expresses the typical time scale at which the mass
and momentum transfer can be modeled with low Mach number hydrodynamics.

With this in mind, we first evaluate the velocity divergence associated
with the constraint using the particular sample of $\M{\Psi}$, 
\[
S=-\left(\rho^{-1}\beta\right)\,\grad\cdot\left[\rho\chi\grad c+\M{\Psi}\left(\d t,\,\widetilde{\V W}\right)\right].
\]
We then define a discrete affine operator $\M{\mathcal{R}}_{F}\left(\d t,\,\widetilde{\V W}\right)$
in terms of its action on the momentum $\V m$ 
\[
\left[\M{\mathcal{R}}_{F}\left(\d t,\,\widetilde{\V W}\right)\right]\left(\V m\right)=\rho\M{\mathcal{R}}_{S}\left(\rho^{-1}\V m\right).
\]
Using this shorthand notation, the momentum update in (\ref{eq:Euler_step_complex})
in the presence of thermal fluctuations can be written as
\[
\M m^{n+1}=\left[\M{\mathcal{R}}_{F}^{n+1}\left(\D t,\,\widetilde{\V W}^{n+1}\right)\right]\left(\V m^{n}+\D t\,\V f^{n}\right).
\]
Observe that this is a conservative momentum update since the application
of $\M{\mathcal{R}}_{F}$ subtracts the (discrete) gradient of a scalar
from the momentum. In actual implementation, it is preferable to apply
$\M{\mathcal{R}}_{F}^{n+1}$ at the beginning of the time step $n+1$
instead of at the end of time step $n$, once the value $S^{n+1}$
is computed from the diffusive and stochastic fluxes for the concentration.

\subsubsection{Euler-Maruyama Update}

Following the above discussion, we can write an Euler-Maruyama temporal
integrator for the low Mach number equations in the shorthand notation,
\begin{eqnarray}
\M m^{n} & = & \left[\M{\mathcal{R}}_{F}^{n}\left(\D t,\,\widetilde{\V W}^{n}\right)\right]\left(\tilde{\V m}^{n}\right)\nonumber \\
\rho_{1}^{n+1} & = & \rho_{1}^{n}+\D t\,\bar{h}^{n}+\check{h}^{n}\left(\D t,\,\widetilde{\V W}^{n}\right)\nonumber \\
\tilde{\V m}^{n+1} & = & \V m^{n}+\D t\,\bar{\V f}^{n}+\check{\V f}^{n}\left(\D t,\,\M W^{n}\right),\label{eq:Euler_step}
\end{eqnarray}
where $\M W^{n}$ and $\widetilde{\V W}^{n}$ are collections of standard
normal variates generated using a pseudo-random number generator independently
at each time step. Here the deterministic increments are written using
the shorthand notation,
\begin{eqnarray*}
\bar{\V f} & = & \grad\cdot\left[-\rho\V v\V v^{T}+\eta\left(\grad\V v+\grad^{T}\V v\right)\right]+\rho\V g\\
\bar{h} & = & \grad\cdot\left(-\rho_{1}\V v+\rho\chi\grad c\right).
\end{eqnarray*}
The stochastic increments are written in terms of
\begin{eqnarray*}
\check{\V f}\left(\d t,\,\M W\right) & = & \left[\grad\cdot\M{\Sigma}\left(\d t,\,\M W\right)\right]\delta t=\grad\cdot\left[\sqrt{\frac{\eta\left(k_{B}T\right)\d t}{\D V}}\,\left(\M W+\M W^{T}\right)\right]\\
\check{h}\left(\d t,\,\widetilde{\V W}\right) & = & \left[\grad\cdot\M{\Psi}\left(\d t,\,\widetilde{\V W}\right)\right]\delta t=\grad\cdot\left[\sqrt{\frac{2\chi\rho\mu_{c}^{-1}\left(k_{B}T\right)\d t}{\D V}}\,\widetilde{\V W}\right],
\end{eqnarray*}
where $\widetilde{\V W}$ and $\M W$ are vectors of standard Gaussian
variables \cite{DFDB}.

\subsection{\label{sub:Higher-Order-Temporal}Higher-Order Temporal Integrators}

A good strategy for composing higher-order temporal integrators for
the low Mach number equations is to use a linear combination of several
projected Euler steps of the form (\ref{eq:Euler_step}). In this
way, the higher-order integrators inherit the properties of the Euler
step. In our case, this will be very useful in constructing conservative
discretizations that strictly maintain the EOS constraint and only
evaluate fluxes at states that strictly obey the EOS constraint.

The incorporation of stochastic forcing in the Runge-Kutta temporal
integrators that we use is described in Refs. \cite{LLNS_S_k,DFDB};
here we only summarize the resulting schemes. We note that the stochastic
terms should be considered additive noise, even though we evaluate
them using an instantaneous state like multiplicative noise \cite{LLNS_Staggered}.

\subsubsection{Explicit Trapezoidal Rule}

A weakly second-order temporal integrator for (\ref{eq:m_gauge},\ref{eq:rho1_gauge})
is provided by the \emph{explicit trapezoidal rule}, in which we first
take a predictor Euler step
\begin{eqnarray}
\M m^{n} & = & \left[\M{\mathcal{R}}_{F}^{n}\left(\D t,\,\widetilde{\V W}^{n}\right)\right]\left(\tilde{\V m}^{n}\right)\nonumber \\
\rho_{1}^{\star,n+1} & = & \rho_{1}^{n}+\D t\,\bar{h}^{n}+\check{h}^{n}\left(\D t,\,\widetilde{\V W}^{n}\right)\\
\tilde{\V m}^{\star,n+1} & = & \V m^{n}+\D t\,\bar{\V f}^{n}+\check{\V f}^{n}\left(\D t,\,\M W^{n}\right).\label{eq:trapezoidal_predictor}
\end{eqnarray}
The corrector step is a linear combination of the predictor and another
Euler update,
\begin{eqnarray}
\M m^{\star,n+1} & = & \left[\M{\mathcal{R}}_{F}^{\star,n+1}\left(\D t,\,\widetilde{\V W}^{n}\right)\right]\left(\tilde{\V m}^{\star,n+1}\right)\nonumber \\
\rho_{1}^{n+1} & = & \frac{1}{2}\rho_{1}^{n}+\frac{1}{2}\left[\rho_{1}^{\star,n+1}+\D t\,\bar{h}^{\star,n+1}+\check{h}^{\star,n+1}\left(\D t,\,\widetilde{\V W}^{n}\right)\right]\\
\tilde{\V m}^{n+1} & = & \frac{1}{2}\V m^{n}+\frac{1}{2}\left[\M m^{\star,n+1}+\D t\,\bar{\V f}^{\star,n+1}+\check{\V f}^{\star,n+1}\left(\D t,\,\M W^{n}\right)\right],\label{eq:trapezoidal_corrector}
\end{eqnarray}
and reuses the same random numbers $\M W^{n}$ and $\widetilde{\V W}^{n}$
as the predictor step.

Note that both the predicted and the corrected values for density
and concentration obey the EOS. We numerically observe that the trapezoidal
rule does exhibit a slow but systematic numerical drift in the EOS,
and therefore it is necessary to use the correction procedure described
in Section \ref{sub:DriftCorrection} at the end of each time step.
The analysis in Ref. \cite{DFDB} indicates that for the incompressible
case the trapezoidal scheme exhibits second-order weak accuracy in
the nonlinear and linearized settings.

\subsubsection{Explicit Midpoint Rule}

An alternative second-order scheme is the\emph{ explicit midpoint
rule}, which can be summarized as follows. First we take a projected
Euler step to estimate midpoint values (denoted here with superscript
\textbf{$\star,n+\myhalf$}),
\begin{eqnarray}
\M m^{n} & = & \left[\M{\mathcal{R}}_{F}^{n}\left(\frac{\D t}{2},\,\widetilde{\V W}_{1}^{n}\right)\right]\left(\tilde{\V m}^{n}\right)\nonumber \\
\rho_{1}^{\star,n+\myhalf} & = & \rho_{1}^{n}+\frac{\D t}{2}\,\bar{h}^{n}+\check{h}^{n}\left(\frac{\D t}{2},\,\widetilde{\V W}_{1}^{n}\right)\nonumber \\
\tilde{\M m}^{\star,n+\myhalf} & = & \V m^{n}+\frac{\D t}{2}\,\bar{\V f}^{n}+\check{\V f}^{n}\left(\frac{\D t}{2},\,\M W_{1}^{n}\right).\label{eq:midpoint_predictor}
\end{eqnarray}
and then we complete the time step with another Euler-like update
\begin{eqnarray}
\M m^{\star,n+\myhalf} & = & \left[\M{\mathcal{R}}_{F}^{n+\myhalf}\left(\D t,\,\widetilde{\V W}^{n}\right)\right]\left(\tilde{\V m}^{\star,n+\myhalf}\right)\nonumber \\
\rho_{1}^{n+1} & = & \rho_{1}^{n}+\D t\,\bar{h}^{\star,n+\myhalf}+\check{h}^{\star,n+\myhalf}\left(\D t,\,\widetilde{\V W}^{n}\right)\nonumber \\
\tilde{\V m}^{n+1} & = & \V m^{n}+\D t\,\bar{\V f}^{\star,n+\myhalf}+\check{\V f}^{\star,n+\myhalf}\left(\D t,\,\M W^{n}\right),\label{eq:midpoint_corrector}
\end{eqnarray}
where the standard Gaussian variates
\[
\widetilde{\V W}^{n}=\frac{\widetilde{\V W}_{1}^{n}+\widetilde{\V W}_{2}^{n}}{\sqrt{2}},
\]
and the vectors of standard normal variates $\widetilde{\V W}_{1}^{n}$
and $\widetilde{\V W}_{2}^{n}$ are independent, and similarly for
$\M W_{1}^{n}$ and $\M W_{2}^{n}$. Note that $\widetilde{\V W}_{1}^{n}$
and $\M W_{1}^{n}$ are used in \emph{both} the predictor and the
corrector stages, while $\widetilde{\V W}_{2}^{n}$ and $\M W_{2}^{n}$
are used in the corrector only. Physically, the random numbers $\M W_{1}^{n}/\sqrt{2}$
(and similarly for $\widetilde{\V W}_{1}^{n}$) correspond to the
increments of the underlying Wiener processes $\D{\V{\mathcal{B}}}_{1}=\sqrt{\D t/2}\,\M W_{1}^{n}$
over the first half of the time step, and the random numbers $\M W_{2}^{n}/\sqrt{2}$
correspond to the Wiener increments for the second half of the timestep
\cite{DFDB}.

Note that both the midpoint and the endpoint values for density and
concentration obey the EOS. We numerically observe that the midpoint
rule does not exhibit a systematic numerical drift in the EOS, and
can therefore be used without the correction procedure described in
Section \ref{sub:DriftCorrection}. The analysis in Ref. \cite{DFDB}
indicates that for the incompressible case the midpoint scheme exhibits
second-order weak accuracy in the nonlinear setting. Furthermore,
in the linearized setting it reproduces the steady-state covariances
of the fluctuating fields to third order in the time step size.

\subsubsection{Three-Stage Runge-Kutta (RK3) Rule}

We have also tested and implemented the three-stage Runge Kutta scheme
that was used in Refs. \cite{LLNS_S_k,LLNS_Staggered}. This scheme
can be expressed as a linear combination of three Euler steps. The
first stage is a predictor Euler step,
\begin{eqnarray}
\M m^{n} & = & \left[\M{\mathcal{R}}_{F}^{n}\left(\D t,\,\widetilde{\V W}^{n}\right)\right]\left(\tilde{\V m}^{n}\right)\nonumber \\
\rho_{1}^{\star} & = & \rho_{1}^{n}+\D t\,\bar{h}^{n}+\check{h}^{n}\left(\D t,\,\widetilde{\V W}^{n}\right)\\
\tilde{\V m}^{\star} & = & \V m^{n}+\D t\,\bar{\V f}^{n}+\check{\V f}^{n}\left(\D t,\,\M W^{n}\right).\label{eq:trapezoidal_predictor-1}
\end{eqnarray}
The second stage is a midpoint predictor
\begin{eqnarray}
\M m^{\star} & = & \left[\M{\mathcal{R}}_{F}^{\star}\left(\D t,\,\widetilde{\V W}^{\star,n}\right)\right]\left(\tilde{\V m}^{\star}\right)\nonumber \\
\rho_{1}^{\star\star} & = & \frac{3}{4}\rho_{1}^{n}+\frac{1}{4}\left[\rho_{1}^{\star}+\D t\,\bar{h}^{\star}+\check{h}^{\star}\left(\D t,\,\widetilde{\V W}^{\star,n}\right)\right]\\
\tilde{\V m}^{\star\star} & = & \frac{3}{4}\V m^{n}+\frac{1}{4}\left[\M m^{\star}+\D t\,\bar{\V f}^{\star}+\check{\V f}^{\star}\left(\D t,\,\M W^{\star,n}\right)\right],\label{eq:RK3_predictor}
\end{eqnarray}
and a final corrector stage completes the time step
\begin{eqnarray}
\M m^{\star\star} & = & \left[\M{\mathcal{R}}_{F}^{\star\star}\left(\D t,\,\widetilde{\V W}^{\star\star,n}\right)\right]\left(\tilde{\V m}^{\star\star}\right)\nonumber \\
\rho_{1}^{n+1} & = & \frac{1}{3}\rho_{1}^{n}+\frac{2}{3}\left[\rho_{1}^{\star\star}+\D t\,\bar{h}^{\star\star}+\check{h}^{\star\star}\left(\D t,\,\widetilde{\V W}^{\star\star,n}\right)\right]\\
\tilde{\M m}^{n+1} & = & \frac{1}{3}\V m^{n}+\frac{2}{3}\left[\M m^{\star\star}+\D t\,\bar{\V f}^{\star\star}+\check{\V f}^{\star\star}\left(\D t,\,\M W^{\star\star,n}\right)\right].\label{eq:RK3_corrector}
\end{eqnarray}
Here the stochastic fluxes between different stages are related to
each other via 
\begin{align}
\V W^{n}= & \V W_{1}^{n}+\frac{\left(2\,\sqrt{2}+\sqrt{3}\right)}{5}\V W_{2}^{n}\nonumber \\
\V W^{\star,n}= & \V W_{1}^{n}+\frac{\left(-4\,\sqrt{2}+3\,\sqrt{3}\right)}{5}\V W_{2}^{n}\nonumber \\
\V W^{\star\star,n}= & \V W_{1}^{n}+\frac{\left(\sqrt{2}-2\,\sqrt{3}\right)}{10}\V W_{2}^{n},\label{RK3_optimal}
\end{align}
where $\V W_{1}^{n}$ and $\V W_{2}^{n}$ are independent and generated
independently at each RK3 step, and similarly for $\widetilde{\V W}$.
The weights of $\V W_{2}^{n}$ are chosen to maximize the weak order
of accuracy of the scheme while still using only two random samples
of the stochastic fluxes per time step \cite{DFDB}.

The RK3 method is third-order accurate deterministically, and stable
even in the absence of diffusion/viscosity (i.e., for advection-dominated
flows). Note that the predicted, the midpoint and the endpoint values
for density and concentration all obey the EOS. We numerically observe
that the RK3 scheme does exhibit a systematic numerical drift in the
EOS, and therefore it is necessary to use the correction procedure
described in Section \ref{sub:DriftCorrection} at the end of each
time step. The analysis in Ref. \cite{DFDB} indicates that for the
incompressible case the RK3 scheme exhibits second-order weak accuracy
in the nonlinear setting. In the linearized setting it reproduces
the steady-state covariances of the fluctuating fields to third order
in the time step size.

\subsection{\label{sub:DriftCorrection}EOS drift}

While in principle our temporal integrators should strictly maintain
the EOS, roundoff errors and the finite tolerance employed in the
iterative Poisson solver lead to a small drift in the constraint that
can, depending on the specific scheme, lead to an exponentially increasing
violation of the EOS over many time steps. In order to maintain the
EOS at all times to within roundoff tolerance, we periodically apply
a globally-conservative $L_{2}$ projection of $\rho$ and $\rho_{1}$
onto the linear EOS constraint.

This projection step consists of correcting $\rho_{1}$ in cell $k$
using 
\[
\left(\rho_{1}\right)_{k}\leftarrow A\left(\rho_{1}\right)_{k}-B\left(\rho_{2}\right)_{k}-\frac{1}{N}\sum_{k^{\prime}}\left[A\left(\rho_{1}\right)_{k^{\prime}}-B\left(\rho_{2}\right)_{k^{\prime}}\right]+\frac{1}{N}\sum_{k^{\prime}}\left(\rho_{1}\right)_{k^{\prime}},
\]
where $N$ is the number of hydrodynamic cells in the system and 
\[
A=\frac{\bar{\rho}_{1}^{2}}{\bar{\rho}_{1}^{2}+\bar{\rho}_{2}^{2}},\quad B=\frac{\bar{\rho}_{1}\bar{\rho}_{2}}{\bar{\rho}_{1}^{2}+\bar{\rho}_{2}^{2}}.
\]
Note that the above update, while nonlocal in nature, conserves the
total mass $\sum_{k^{\prime}}\left(\rho_{1}\right)_{k^{\prime}}$.
A similar update applies to $\rho_{2}$, or equivalently, $\rho=\rho_{1}+\rho_{2}$.

\section{\label{sec:SpatialDiscretization}Spatial Discretization}

The spatial discretization we employ follows closely the spatial discretization
of the constant-coefficient incompressible equations described in
Ref. \cite{LLNS_Staggered}. Therefore, we focus here on the differences,
specifically, the use of conserved variables, the handling of the
variable-density projection and variable-coefficient diffusion, and
the imposition of the low Mach number constraint. Note that the handling
of the stochastic momentum and mass fluxes is identical to that described
in Ref. \cite{LLNS_Staggered}.

\begin{figure}[tb]
\centering{}\includegraphics[width=6in]{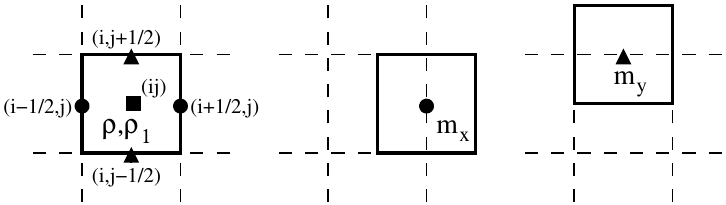}
\caption{\label{fig:grid}Staggered (MAC) finite-volume discretization on a
uniform Cartesian two-dimensional grid. (\emph{Left}) Control volume
and flux discretization for cell-centered scalar fields, such as densities
$\rho$ and $\rho_{1}$. (\emph{Middle}) Control volume for the $x$-component
of face-centered vector fields, such as $m_{x}$ (\emph{Right}) Control
volume for the $y$-component of face-centered vector fields, such
as $m_{y}$. }
\end{figure}

For simplicity of notation, we focus on two dimensional problems,
with straightforward generalization to three spatial dimensions. Our
spatial discretization follows the commonly-used MAC approach \cite{HarWel65},
in which the scalar conserved quantities $\rho$ and $\rho_{1}$ are
defined on a regular Cartesian grid. The vector conserved variables
$\V m=\rho\V v$ are defined on a staggered grid, such that the $k^{{\rm th}}$
component of momentum is defined on the faces of the scalar variable
Cartesian grid in the $k^{{\rm th}}$ direction, see Fig. \ref{fig:grid}.
For simplicity of notation, we often denote the different components
of velocity as $\V v=\left(u,v\right)$ in two dimensions and $\V v=\left(u,v,w\right)$
in three dimensions. The terms ``cell-centered'', ``edge-centered'',
and ``face-centered'' refer to spatial locations relative to the
underlying scalar grid. Our discretization is based on calculating
fluxes on the faces of a finite-volume grid and is thus locally conservative.
It is important to note, however, that for the MAC grid different
control volumes are used for the scalars and the components of the
momentum, see Fig. \ref{fig:grid}.

From the cell-centered $\rho$ and $\rho_{1}$ we can define other
cell-centered scalar quantities, notably, the concentration $c_{i,j}=\left(\rho_{1}\right)_{i,j}/\rho_{i,j}$
and the transport quantities $\chi_{i,j}$ and $\eta_{i,j}$, which
typically depend on the local density $\rho_{i,j}$ and concentration
$c_{i,j}$ (and temperature for non-isothermal models), and can, in
general, also depend on the spatial position of the cell $\left(x,y\right)=\left(i\D x,\, j\D y\right)$.
In order to define velocities we need to interpret the continuum relationship
$\V m=\rho\V v$ on the staggered grid. This is done by defining face-centered
scalar quantities obtained as an arithmetic average of the corresponding
cell-centered quantities in the two neighboring cells. Specifically,
we define 
\begin{equation}
\rho_{i+\myhalf,j}=\frac{\rho_{i,j}+\rho_{i+1,j}}{2},\quad u_{i+\myhalf,j}=\frac{\left(m_{x}\right)_{i+\myhalf,j}}{\rho_{i+\myhalf,j}},\label{eq:rho_face}
\end{equation}
except at physical boundaries, where the value is obtained from the
imposed boundary conditions (see Section \ref{sub:BoundaryConditions}).
Arithmetic averaging is only one possible interpolation from cells
to faces \cite{almgren-iamr}. In general, other forms of averaging
such as a harmonic or geometric average or higher-order, wider stencils
\cite{Projection4thOrder_SDC,LLNS_S_k} can be used. Most components
of the spatial discretization can easily be generalized to other choices
of interpolation. As we explain later, the use of linear averaging
simplifies the construction of conservative advection.

\subsection{\label{sub:Diffusion}Diffusion}

In this section we describe the spatial discretization of the diffusive
mass flux term $\nabla\cdot\rho\chi\nabla c$ in (\ref{eq:rho1_eq}).
The discretization is based on conservative centered differencing
\cite{LLNS_S_k,DFDB},

\begin{equation}
(\nabla\cdot\rho\chi\nabla c)_{i,j}=\Delta x^{-1}\left[\left(\rho\chi\frac{\partial c}{\partial x}\right)_{i+\myhalf,j}-\left(\rho\chi\frac{\partial c}{\partial x}\right)_{i-\myhalf,j}\right]+\Delta y^{-1}\left[\left(\rho\chi\frac{\partial c}{\partial y}\right)_{i,j+\myhalf}-\left(\rho\chi\frac{\partial c}{\partial y}\right)_{i,j-\myhalf}\right],\label{eq:discrete_diffusion}
\end{equation}
 where, for example,
\begin{equation}
\left(\rho\chi\frac{\partial c}{\partial x}\right)_{i+\myhalf,j}=\left(\rho_{i+\myhalf,j}\right)\left(\chi_{i+\myhalf,j}\right)\left(\frac{c_{i+1,j}-c_{i,j}}{\Delta x}\right),\label{eq:discrete_diff_flux}
\end{equation}
and $\chi_{i+\myhalf,j}$ is an interpolated face-centered diffusion
coefficient, for example, as done for $\rho$ in Eq. (\ref{eq:rho_face}),
\[
\chi_{i+\myhalf,j}=\frac{\chi_{i,j}+\chi_{i+1,j}}{2},
\]
except at physical boundaries, where the value is obtained from the
imposed boundary conditions.

Regardless of the specific form of the interpolation operator, the
same face-centered diffusion coefficient $\chi_{i+\myhalf,j}$ must
be used when calculating the magnitude of the stochastic mass flux
on face $\left(i+\myhalf,j\right)$,
\[
\left(\Psi_{x}\right)_{i+\myhalf,j}=\sqrt{2\chi_{i+\myhalf,j}\left(\rho\mu_{c}^{-1}\right)_{i+\myhalf,j}\, k_{B}T}\;\widetilde{\mathcal{W}}_{i+\myhalf,j}.
\]
This matches the covariance of the discrete stochastic mass increments
$\nabla\cdot\M{\Psi}$ with the discretization of the diffusive dissipation
operator $\nabla\cdot\rho\chi\nabla$ given in (\ref{eq:discrete_diffusion},\ref{eq:discrete_diff_flux}).
This matching ensures discrete fluctuation-dissipation balance in
the linearized setting \cite{LLNS_S_k}. Specifically, at thermodynamic
equilibrium the static covariance of the concentration is determined
from the equilibrium value of $\left(\rho\mu_{c}^{-1}\right)$ (thermodynamics)
independently of the particular values of the transport coefficients
(dynamics), as seen in (\ref{eq:S_equilibrium}) and dictated by statistical
mechanics principles.

\subsection{Viscous Terms}

In Ref. \cite{LLNS_Staggered} a Laplacian form of the viscous term
$\eta\grad^{2}\V v$ is assumed, which is not applicable when viscosity
is spatially varying and $\grad\cdot\V v=S\neq0$. In two dimensions,
the divergence of the viscous stress tensor in the momentum equation
(\ref{eq:momentum_eq}), neglecting bulk viscosity effects, is 
\begin{eqnarray}
\grad\cdot\left[\eta\left(\grad\V v+\grad^{T}\V v\right)\right] & = & \left[\begin{array}{c}
2\frac{\partial}{\partial x}\left(\eta\frac{\partial u}{\partial x}\right)+\frac{\partial}{\partial y}\left(\eta\frac{\partial u}{\partial y}+\eta\frac{\partial v}{\partial x}\right)\\
2\frac{\partial}{\partial y}\left(\eta\frac{\partial v}{\partial y}\right)+\frac{\partial}{\partial x}\left(\eta\frac{\partial v}{\partial x}+\eta\frac{\partial u}{\partial y}\right)
\end{array}\right].\label{eq:div_viscous_stress}
\end{eqnarray}
The discretization of the viscous terms requires $\eta$ at cell-centers
and edges (note that in two dimensions the edges are the same as the
nodes $\left(i+\myhalf,\, j+\myhalf\right)$ of the grid). The value
of $\eta$ at a node is interpolated as the arithmetic average of
the four neighboring cell-centers,
\[
\eta_{i+\myhalf,j+\myhalf}=\frac{1}{4}\left(\eta_{i,j}+\eta_{i+1,j+1}+\eta_{i+1,j}+\eta_{i,j+1}\right),
\]
except at physical boundaries, where the values are obtained from
the prescribed boundary conditions. The different viscous friction
terms are discretized by straightforward centered differences. Explicitly,
for the $x$-component of momentum 
\[
\left[\frac{\partial}{\partial x}\left(\eta\frac{\partial u}{\partial x}\right)\right]_{i+\myhalf,j}=\D x^{-1}\left[\left(\eta\frac{\partial u}{\partial x}\right)_{i+1,j}-\left(\eta\frac{\partial u}{\partial x}\right)_{i,j}\right]
\]
 with 
\[
\left(\eta\frac{\partial u}{\partial x}\right)_{i,j}=\eta_{i,j}\left(\frac{u_{i+\myhalf,j}-u_{i-\myhalf,j}}{\Delta x}\right).
\]
Similarly, for the term involving a second derivative in $y$, 
\[
\left[\frac{\partial}{\partial y}\left(\eta\frac{\partial u}{\partial y}\right)\right]_{i+\myhalf,j}=\D y^{-1}\left[\left(\eta\frac{\partial u}{\partial y}\right)_{i+\myhalf,j+\myhalf}-\left(\eta\frac{\partial u}{\partial y}\right)_{i+\myhalf,j-\myhalf}\right],
\]
 with 
\[
\left(\eta\frac{\partial u}{\partial y}\right)_{i+\myhalf,j+\myhalf}=\eta_{i+\myhalf,j+\myhalf}\left(\frac{u_{i+\myhalf,j+1}-u_{i+\myhalf,j}}{\Delta y}\right).
\]
A similar construction is used for the mixed-derivative term,
\[
\left[\frac{\partial}{\partial y}\left(\eta\frac{\partial v}{\partial x}\right)\right]_{i+\myhalf,j}=\D y^{-1}\left[\left(\eta\frac{\partial v}{\partial x}\right)_{i+\myhalf,j+\myhalf}-\left(\eta\frac{\partial v}{\partial x}\right)_{i+\myhalf,j-\myhalf}\right],
\]
 with 
\[
\left(\eta\frac{\partial v}{\partial x}\right)_{i+\myhalf,j+\myhalf}=\eta_{i+\myhalf,j+\myhalf}\left(\frac{v_{i+1,j+\myhalf}-v_{i,j+\myhalf}}{\Delta x}\right).
\]

The stochastic stress tensor discretization is described in more detail
in Ref. \cite{LLNS_Staggered} and applies in the present context
as well. For the low Mach number equations, just as for the compressible
equations, the symmetric form of the stochastic stress tensor must
be used in order to ensure discrete fluctuation-dissipation balance
between the viscous dissipation and stochastic forcing. Additionally,
when $\eta$ is not spatially uniform the same interpolated viscosity
$\eta_{i+\myhalf,j+\myhalf}$ as used in the viscous terms must be
used when calculating the amplitude in the stochastic forcing $\sqrt{\eta k_{B}T}$
at the edges (nodes) of the grid.

\subsection{\label{sub:Advection}Advection}

It is challenging to construct spatio-temporal discretizations that
conserve the total mass while remaining consistent with the equation
of state \cite{ZeroMach_Klein,LowMach_DiscreteCompatibility,LowMach_FiniteDifference},
as ensured in the continuum context by the constraint (\ref{eq:div_v_constraint}).
We demonstrate here how the special linear form of the constraint
(\ref{eq:EOS_quasi_incomp}) can be exploited in the discrete context.
Following Ref. \cite{LLNS_Staggered}, we spatially discretize the
advective terms in (\ref{eq:rho1_eq}) using a centered (skew-adjoint
\cite{ConservativeDifferences_Incompressible}) discretization,
\begin{equation}
\left[\nabla\cdot\left(\rho_{1}\V v\right)\right]_{i,j}=\Delta x^{-1}\left[\left(\rho_{1}\right)_{i+\myhalf,j}u_{i+\myhalf,j}-\left(\rho_{1}\right)_{i-\myhalf,j}u_{i-\myhalf,j}\right]+\Delta y^{-1}\left[\left(\rho_{1}\right)_{i,j+\myhalf}v_{i,j+\myhalf}-\left(\rho_{1}\right)_{i,j-\myhalf}v_{i,j-\myhalf}\right],\label{eq:rho_adv}
\end{equation}
and similarly for (\ref{eq:rho_eq}). We would like this discrete
advection to maintain the equation of state (\ref{eq:EOS_quasi_incomp})
at the discrete level, that is, maintain the constraint relating $\left(\rho_{1}\right)_{i,j}$
and $\left(\rho_{2}\right)_{i,j}$ in every cell $\left(i,j\right)$.

Because the different dimensions are decoupled and the divergence
is simply the sum of the one-dimensional difference operators, it
is sufficient to consider (\ref{eq:rho1_eq}) in one spatial dimension.
The method of lines discretization is given by the system of ODEs,
one differential equation per cell $i$,
\[
\left(\partial_{t}\rho_{1}\right)_{i}=\D x^{-1}\left(F_{i+\myhalf}-F_{i-\myhalf}\right)-\D x^{-1}\left[\left(\rho_{1}\right)_{i+\myhalf}u_{i+\myhalf}-\left(\rho_{1}\right)_{i-\myhalf}u_{i-\myhalf}\right],
\]
and similarly for $\left(\partial_{t}\rho_{2}\right)_{i}$. As a shorthand,
denote the quantity that appears in (\ref{eq:EOS_quasi_incomp}) with
\[
\delta=\frac{\rho_{1}}{\bar{\rho}_{1}}+\frac{\rho_{2}}{\bar{\rho}_{2}}=1.
\]
If we use the linear interpolation (\ref{eq:rho_face}) to calculate
face-centered densities, then because of the linearity of the EOS
the face-centered densities obey the EOS if the cell-centered ones
do, since $\delta_{i+\myhalf}=(\delta_{i}+\delta_{i+1})/2=1$. The
rate of change of $\delta$ in cell $i$ is
\begin{eqnarray*}
\D x\left(\partial_{t}\delta\right)_{i} & = & \left(\rho^{-1}\beta\right)\left(F_{i+\myhalf}-F_{i-\myhalf}\right)-\left[\delta_{i+\myhalf}u_{i+\myhalf}-\delta_{i-\myhalf}u_{i-\myhalf}\right]\\
 & = & \left(\rho^{-1}\beta\right)\left(F_{i+\myhalf}-F_{i-\myhalf}\right)-\left(u_{i+\myhalf}-u_{i-\myhalf}\right)=0.
\end{eqnarray*}

This simple calculation shows that the EOS constraint $\delta=1$
is obeyed discretely in each cell at all times if it is initially
satisfied and the velocities used to advect mass obey the discrete
version of the constraint (\ref{eq:div_v_constraint}),
\begin{eqnarray}
S_{i,j} & = & \Delta x^{-1}\left(u_{i+\myhalf,j}-u_{i-\myhalf,j}\right)+\Delta y^{-1}\left(v_{i,j+\myhalf}-v_{i,j-\myhalf}\right)\label{eq:div_v_MAC}\\
 & = & \left(\frac{1}{\bar{\rho}_{1}}-\frac{1}{\bar{\rho}_{2}}\right)\left[\Delta x^{-1}\left(F_{i+\myhalf,j}-F_{i-\myhalf,j}\right)+\Delta y^{-1}\left(F_{i,j+\myhalf}-F_{i,j-\myhalf}\right)\right],\nonumber 
\end{eqnarray}
in two dimensions. Our algorithm ensures that advective terms are
always evaluated using a discrete velocity field that obeys this constraint.
This is accomplished by using a discrete projection operator, as we
describe in the next section.

The spatial discretization of the advection terms in the momentum
equation (\ref{eq:momentum_eq}) is constructed using centered differences
on the corresponding shifted (staggered) grid, as described in Ref.
\cite{LLNS_Staggered}. For example, for the $x$-component of momentum
$m_{x}=\rho u$, 
\begin{equation}
\left[\nabla\cdot\left(m_{x}\V v\right)\right]_{i+\myhalf,j}=\D x^{-1}\left[(m_{x}u)_{i+1,j}-(m_{x}u)_{i,j}\right]+\D y^{-1}\left[(m_{x}v)_{i+\myhalf,j+\myhalf}-(m_{x}v)_{i+\myhalf,j-\myhalf}\right],\label{eq:mom_adv}
\end{equation}
 where simple averaging is used to interpolate momenta to the cell
centers and edges (nodes) of the grid, for example, 
\begin{equation}
(m_{x}u)_{i,j}=(m_{x})_{i,j}u_{i,j}=\left(\frac{\left(m_{x}\right)_{i-\myhalf,j}+\left(m_{x}\right)_{i+\myhalf,j}}{2}\right)\left(\frac{u_{i-\myhalf,j}+u_{i+\myhalf,j}}{2}\right).
\end{equation}
Because of the linearity of the interpolation procedure, the interpolated
discrete velocity used to advect $m_{x}$ obeys the constraint (\ref{eq:div_v_MAC})
on the shifted grid, with a right-hand side $S_{i+\myhalf,j}$ interpolated
using the same arithmetic average used to interpolate the velocities.
In particular, in the incompressible case all variables, including
momentum, are advected using a discretely divergence-free velocity,
ensuring discrete fluctuation-dissipation balance \cite{LLNS_Staggered,DFDB}.

It is well-known that the centered discretization of advection we
employ here is not robust for advection-dominated flows, and higher-order
limiters and upwinding schemes are generally preferred in the deterministic
setting \cite{bellColellaGlaz:1989}. However, these more robust advection
schemes add artificial dissipation, which leads to a violation of
discrete fluctuation-dissipation balance \cite{DFDB}. In Appendix
\ref{AppendixFiltering} we describe an alternative filtering procedure
that can be used to handle strong advection while continuing to use
centered differencing.

\subsection{Discrete Projection}

We now briefly discuss the spatial discretization of the affine operator
$\M{\mathcal{R}}_{S}$ defined by (\ref{P_tilde_v}), as used in our
explicit temporal integrators. The discrete projection takes a face-centered
(staggered) discrete velocity field $\tilde{\V v}=\left(\tilde{u},\,\tilde{v}\right)$
and a velocity divergence $S$ and projects $\V v=\M{\mathcal{R}}_{S}\left(\tilde{\V v}\right)$
onto the constraint (\ref{eq:div_v_MAC}) in a conservative manner.
Specifically, the projection consists of finding a cell-centered discrete
scalar field $\phi$ such that
\[
\rho\V v=\rho\tilde{\V v}-\grad\phi,\mbox{ and }\grad\cdot\V v=S,
\]
where the gradient is discretized using centered differences, e.g.,
\begin{equation}
v_{i+\myhalf,j}=\tilde{v}_{i+\myhalf,j}-\left(\frac{1}{\rho_{i+\myhalf,j}}\right)\left(\frac{\phi_{i+1,j}-\phi_{i,j}}{\Delta x}\right).\label{eq:u_minus_gradp}
\end{equation}
The pressure correction $\phi$ is the solution to the variable-coefficient
discrete Poisson equation,
\begin{eqnarray}
\frac{1}{\Delta x}\left[\left(\frac{1}{\rho_{i+\myhalf,j}}\right)\left(\frac{\phi_{i+1,j}-\phi_{i,j}}{\Delta x}\right)-\left(\frac{1}{\rho_{i-\myhalf,j}}\right)\left(\frac{\phi_{i,j}-\phi_{i,j-1}}{\Delta x}\right)\right]\nonumber \\
+\frac{1}{\Delta y}\left[\left(\frac{1}{\rho_{i,j+\myhalf}}\right)\left(\frac{\phi_{i,j+1}-\phi_{i,j}}{\Delta y}\right)-\left(\frac{1}{\rho_{i,j-\myhalf}}\right)\left(\frac{\phi_{i,j}-\phi_{i,j-1}}{\Delta y}\right)\right]\nonumber \\
=S_{i,j}-\left[\left(\frac{\tilde{u}_{i+\myhalf,j}-\tilde{u}_{i-\myhalf,j}}{\Delta x}\right)+\left(\frac{\tilde{v}_{i,j+\myhalf}-\tilde{v}_{i,j-\myhalf}}{\Delta y}\right)\right],\label{eq:discrete_Poisson}
\end{eqnarray}
which can be solved efficiently using a standard multigrid approach
\cite{almgren-iamr}.

\subsection{\label{sub:BoundaryConditions}Boundary Conditions}

The handling of different types of boundary conditions is relatively
straightforward when a staggered grid is used and the physical boundaries
are aligned with the cell boundaries for the scalar grid. Interpolation
is not used to obtain values for faces, nodes or edges of the grid
that lie on a physical boundary, since this would require ``ghost''
values at cell centers lying outside of the physical domain. Instead,
whenever a value of a physical variable is required at a face, node,
or edge lying on a physical boundary, the boundary condition is used
to obtain that value. Similarly, centered differences for the diffusive
and viscous fluxes that require values outside of the physical domain
are replaced by one-sided differences that only use values from the
interior cell bordering the boundary and boundary values.

For example, if the concentration is specified at the face $\left(i+\myhalf,\, j\right)$,
the diffusive flux discretization (\ref{eq:discrete_diff_flux}) is
replaced with
\[
\left(\rho\chi\frac{\partial c}{\partial x}\right)_{i+\myhalf,j}=\left(\rho_{i+\myhalf,j}\right)\left(\chi_{i+\myhalf,j}\right)\left(\frac{c_{i+\myhalf,j}-c_{i,j}}{\Delta x/2}\right),
\]
where $c_{i+\myhalf,j}$ is the specified boundary value, the density
$\rho_{i+\myhalf,j}$ is obtained from $c_{i+\myhalf,j}$ using the
EOS constraint, and the diffusion coefficient $\chi_{i+\myhalf,j}$
is calculated at the specified values of concentration and density.
Similar straightforward one-sided differencing is used for the viscous
fluxes. As discussed in Ref. \cite{LLNS_Staggered}, the use of second-order
one-sided differencing is not required to achieve global second-order
accuracy, and would make the handling of the stochastic fluxes more
complicated because it leads to a non-symmetric discrete Laplacian.
Note that for the nonlinear low Mach number equations our approach
is subtly different from linearly extrapolating the value in the ghost
cell $c_{i+1,j}=2c_{i+\myhalf,j}-c_{i}$. Namely, the extrapolated
value might be unphysical, and it might not be possible to evaluate
the EOS or transport coefficients at the extrapolated concentration.
For Neumann-type or zero-flux boundary conditions, the corresponding
diffusive flux is set to zero for any faces of the corresponding control
volume that lie on physical boundaries, and values in cells outside
of the physical domain are never required. The corresponding handling
of the stochastic fluxes is discussed in detail in Ref. \cite{LLNS_Staggered}.

The evaluation of advective fluxes for the scalars requires normal
components of the velocity at the boundary. For faces of the grid
that lie on a physical boundary, the normal component of the velocity
is determined from the value of the diffusive mass flux at that face
using (\ref{eq:v_n_BC}). Therefore, these velocities are not independent
variables and are not solved for or modified by the projection $\M{\mathcal{R}}_{S}$.
Specifically, the discrete pressure $\phi$ is only defined at the
cell centers in the interior of the grid, and the discrete Poisson
equation (\ref{eq:discrete_Poisson}) is only imposed on the interior
faces of the grid. Therefore, no explicit boundary conditions for
$\phi$ are required when the staggered grid is used, and the natural
homogeneous Neumann conditions are implied. Advective momentum fluxes
are only evaluated on the interior faces and thus do not use any values
outside of the physical domain.

\subsection{Summary of Euler-Maruyama Method}

By combining the spatial discretization described above with one of
the temporal integators described in Section \ref{sec:TemporalIntegration},
we can obtain a finite-volume solver for the fluctuating low Mach
equations. For the benefit of the reader, here we summarize our implementation
of a single Euler step (\ref{eq:Euler_step}). This forms the core
procedure that the higher-order Runge-Kutta schemes employ several
times during one time step.
\begin{enumerate}
\item Generate the vectors of standard Gaussian variates $\M W^{n}$ and
$\widetilde{\V W}^{n}$.
\item Calculate diffusive and stochastic fluxes for $\rho_{1}$ using (\ref{eq:discrete_diff_flux}),
\[
\V F^{n}=\left(\rho\chi\grad c\right)^{n}+\M{\Psi}^{n}\left(\D t,\,\widetilde{\V W}^{n}\right).
\]

\item Solve the Poisson problem (\ref{eq:discrete_Poisson}) with 
\[
S^{n}=-\left(\frac{1}{\bar{\rho}_{1}}-\frac{1}{\bar{\rho}_{2}}\right)\,\grad\cdot\V F^{n}
\]
to obtain the velocity $\V v^{n}$ from $\tilde{\V v}^{n}=\tilde{\V m}^{n}/\rho^{n}$
using (\ref{eq:u_minus_gradp}), enforcing $\grad\cdot\V v^{n}=S^{n}$.
\item Calculate viscous and stochastic momentum fluxes using (\ref{eq:div_viscous_stress}),
\[
\grad\cdot\left[\eta\left(\grad\V v+\grad^{T}\V v\right)\right]^{n}+\grad\cdot\left[\M{\Sigma}^{n}\left(\D t,\,\M W^{n}\right)\right].
\]

\item Calculate external forcing terms for the momentum equation, such as
the contribution $-\rho^{n}\V g$ due to gravity.
\item Calculate advective fluxes for mass and momentum using (\ref{eq:rho_adv})
and (\ref{eq:mom_adv}).
\item Update mass and momentum densities, including advective, diffusive,
stochastic and external forcing terms, to obtain $\rho^{n+1}$, $\rho_{1}^{n+1}$
and $\tilde{\V m}^{n+1}$. Note that this update preserves the EOS
constraint as explained in Section \ref{sub:Advection}.
\end{enumerate}
We have tested and validated the accuracy of our methods and numerical
implementation using a series of standard deterministic tests, as
well as by examining the equilibrium spectrum of the concentration
and velocity fluctuations \cite{LLNS_S_k,LLNS_Staggered,DFDB}. The
next two sections present further verification and validation in the
context of nonequilibrium systems.

\section{\label{sec:GiantFluct}Giant Concentration Fluctuations}

Advection of concentration by thermal velocity fluctuations in the
presence of large concentration gradients leads to the appearance
of \emph{giant fluctuations} of concentration, as has been studied
theoretically and experimentally for more than a decade \cite{GiantFluctuations_Theory,GiantFluctuations_Cannell,FractalDiffusion_Microgravity,GiantFluctuations_Summary}.
These giant fluctuations were previously simulated in the absence
of gravity in three dimensions by some of us in Ref. \cite{LLNS_Staggered},
and good agreement was found with experimental results \cite{FractalDiffusion_Microgravity}.
In those previous studies the incompressible equations were used,
that is, it was assumed that concentration was a passively-advected
scalar. However, it is more physically realistic to account for the
fact that the properties of the fluid, notably the density and the
transport coefficients, depend on the concentration. In Ref. \cite{GiantFluctuations_Cannell}
a series of experiments were performed to study the temporal evolution
of giant concentration fluctuations during the diffusive mixing of
water and glycerol, starting with a glycerol mass fraction of $c=0.39$
in the bottom half of the experimental domain, and $c=0$ in the top
half. Because it is essentially impossible to analytically solve the
full system of fluctuating equations in the presence of spatial inhomogeneity
and nontrivial boundary conditions, the existing theoretical analysis
of the diffusive mixing process \cite{GiantFluctuations_Theory} makes
a quasi-periodic constant-coefficient incompressible approximation.

For simplicity, in this section we focus on a time-independent problem
and study the spectrum of steady-state concentration fluctuations
in a mixture under gravity in the presence of a constant concentration
gradient. This extends the study reported in Ref. \cite{LLNS_Staggered}
to account for the fact that the density, viscosity, and diffusion
coefficient depend on the concentration. For simplicity, we do two-dimensional
simulations, since there are no qualitative differences between the
spectrum of concentration fluctuations in two and three dimensions
\cite{LLNS_Staggered} (note, however, that in real space, unlike
in Fourier space, the effect of the fluctuations on the transport
is very different in two and three dimensions). Furthermore, in these
simulations we do not include a stochastic flux in the concentration
equation, i.e., we set $\M{\Psi}=0$, so that all fluctuations in
the concentration arise from being out of thermodynamic equilibrium.
With this approximation we do not need to model the chemical potential
of the mixture and obtain $\mu_{c}$. This formulation is justified
by the fact that it is known experimentally that the nonequilibrium
fluctuations are much larger than the equilibrium ones for the conditions
we consider \cite{GiantFluctuations_Cannell}.

In the simple linearized theory presented in Section \ref{sub:GiantTheory}
several approximations are made. The first one is that a quasi-periodic
approximation is used even though the actual system is not periodic
in the $y$ direction. This source of error has already been studied
numerically in Ref. \cite{LLNS_Staggered}. We also use a Boussinesq
approximation where it is assumed that $\bar{\rho}_{1}=\rho_{0}+\D{\rho}/2$
and $\bar{\rho}_{2}=\rho_{0}-\D{\rho}/2$, where $\D{\rho}$ is a
small density difference between the two fluids, $\D{\rho}/\rho_{0}\ll1$,
so that density is approximately constant and $\beta\ll1$. More precisely,
in the Boussinesq model the gravity term in the velocity equation
only enters through the product $\beta g$ so the approximation consists
of taking the limit $\beta\rightarrow0$ and $g\rightarrow\infty$
while keeping the product $\beta g$ fixed. The final approximation
made in the simple theory is that the transport coefficients, i.e.,
the viscosity and diffusion coefficients, are assumed to be constant.
Here we evaluate the validity of the constant-coefficient constant-density
approximation ($\rho,\,\eta$ and $\chi$ constant, $\beta\rightarrow0$),
as well as the constant-density (Boussinesq) approximation alone ($\rho$
constant, $\beta\rightarrow0$, but variable $\eta,\,\chi$), by comparing
with the solution to the complete low Mach number equations ($\rho,\,\eta,\,\chi$
and $\beta$ variable).

\subsection{Simulation Parameters}

We base our parameters on the experimental studies of diffusive mixing
in a water-glycerol mixture, as reported in Ref. \cite{GiantFluctuations_Cannell}.
The physical domain is $1\,\textrm{cm}\times0.25\,\textrm{cm}$ discretized
on a uniform $128\times32$ two dimensional grid, with a thickness
of $1\,\textrm{cm}$ along the $z$ direction. Gravity is applied
in the negative $y$ (vertical) direction. Reservoir boundary conditions
(\ref{eq:v_n_BC}) are applied in the $y$-direction and periodic
boundary conditions in the $x$-direction. We set the concentration
to $c=0.39$ on the bottom boundary and $c=0$ on the top boundary,
and apply no-slip boundary conditions for the velocity at both boundaries.
The initial condition is $c(t=0)=0.39\left(y/0.25-1\right)$, which
is close to the deterministic steady-state profile. A very good fit
to the experimental equation of state (dependence of density on concentration
at standard temperature and pressure) over the whole range of concentrations
of interest is provided by the EOS (\ref{eq:EOS_quasi_incomp}) with
the density of water set to $\bar{\rho}_{2}=1\,\textrm{g}/\textrm{cm}^{3}$
and the density of glycerol set to $\bar{\rho}_{1}=1.29\,\textrm{g}/\textrm{cm}^{3}$.
In these simulations the magnitude of the velocity fluctuations is
very small and we did not use filtering (see Appendix \ref{AppendixFiltering}).

Experimentally, the dependence of viscosity on glycerol mass fraction
has been fit to an exponential function \cite{GiantFluctuations_Cannell},
which we approximate with a quadratic function over the range of concentrations
of interest,
\begin{equation}
\eta(c)=\rho(c)\nu(c)=\rho_{0}\nu_{0}\exp(2.06c+2.32c^{2})\approx\rho_{0}\nu_{0}(1.0+0.66c+12c^{2}),\label{eq:nu_c}
\end{equation}
where $\rho_{0}=1\,\textrm{g}/\textrm{cm}^{3}$ and experimental measurements
estimate $\nu_{0}\approx10^{-2}\,\textrm{cm}^{2}/\textrm{s}$. The
diffusion coefficient dependence on the concentration has been studied
experimentally \cite{WaterGlycerolDiffusion}, but is in fact strongly
affected by thermal fluctuations and spatial confinement \cite{DiffusionRenormalization,Nanopore_FluctuationsPRE,StokesLaw}.
We approximate the dependence assuming a Stokes-Einstein relation
\cite{StokesLaw}, which is in reasonable agreement with the experimental
results in Ref. \cite{WaterGlycerolDiffusion} over the range of concentrations
of interest here,
\begin{equation}
\chi(c)=\frac{\chi_{0}\eta_{0}}{\eta(c)}\approx\chi_{0}(1.0-2.2c+1.2c^{2}),\label{eq:chi_c}
\end{equation}
where experimental estimates \cite{WaterGlycerolDiffusion} for water-glycerol
mixtures give $\chi_{0}\approx10^{-5}\,\textrm{cm}^{2}/\textrm{s}$,
with a Schmidt number $\mbox{Sc}=\nu/\chi\approx10^{3}$. This very
large separation of scales between mass and momentum diffusion is
not feasible to simulate with our explicit temporal integration methods.
Referring back to the simplified theory (\ref{eq:S_c_c}), which in
this case can be simplified further to
\begin{equation}
S_{c,c}\left(k_{x},k_{y}=0\right)=\av{\left(\widehat{\d c}\right)\left(\widehat{\d c}\right)^{\star}}\approx\left(\frac{\nu}{\nu+\chi}\right)\frac{k_{B}T}{\left(\chi\eta k_{x}^{4}+h_{\parallel}\rho g\beta\right)}\, h_{\parallel}^{2},\label{eq:S_c_c_simp}
\end{equation}
we see that for $\nu\gg\chi$ the shape of the spectrum of the steady-state
concentration fluctuations, and in particular, the cutoff wavenumber
due to gravity, is determined from the product $\chi\eta$ and not
$\chi$ and $\eta$ individually. Therefore, as also done in Ref.
\cite{LLNS_Staggered}, we choose $\chi_{0}$ and $\nu_{0}$ so that
$\chi\left(\bar{c}\right)\eta\left(\bar{c}\right)$ is kept at the
physical value of $10^{-7}\,\textrm{g}\cdot\textrm{cm}/\mbox{s}^{2}$
but the Schmidt number is reduced by two orders of magnitude, $S_{c}=\rho_{0}^{-1}\eta\left(\bar{c}\right)/\chi\left(\bar{c}\right)=10$,
where $\bar{c}=0.39/2$ is an estimate of the average concentration.
The condition $\eta\left(\bar{c}\right)\approx10^{-3}\,\textrm{g}/\left(\textrm{cm}\cdot\textrm{s}\right)$
and $\chi\left(\bar{c}\right)\approx10^{-4}\,\textrm{cm}^{2}/\textrm{s}$
gives our simulation parameters $\nu_{0}\approx6.1\times10^{-4}\,\textrm{cm}^{2}/\textrm{s}$
and $\chi_{0}\approx1.6\times10^{-4}\,\textrm{cm}^{2}/\textrm{s}$.

The physical value for gravity is $g\approx10^{3}\,\textrm{cm}/\textrm{s}^{2}$
and the solutal expansion coefficient $\beta\left(\bar{c}\right)\approx0.234$
follows from $\bar{\rho}_{1}$ and $\bar{\rho}_{2}$. When employing
the Boussinesq approximation, in which gravity only enters through
the product $\beta g,$ we set $\rho_{1}=1.054$ and $\rho_{2}=1.044$
so that $\beta=0.01$ and increase gravity by the corresponding factor
to $g=2.34\cdot10^{4}\,\textrm{cm}/\textrm{s}^{2}$ in order to keep
$\beta g$ fixed at the physical value. We also performed simulations
with a weaker gravity, $g\approx10^{2}\,\textrm{cm}/\textrm{s}^{2}$,
which enhances the nonequilibrium fluctuations, as well as no gravity,
which makes the fluctuations truly giant \cite{FractalDiffusion_Microgravity}.

\subsection{Results}

We employ the explicit midpoint temporal integrator (which we recall
is third-order accurate for static covariances) and set $\Delta t=0.005\, s$,
which results in a diffusive Courant number $\nu\Delta t/\D x^{2}\approx0.1$.
We skip the first 50,000 time steps (about 5 diffusion crossing times)
and then collect samples from the subsequent 50,000 time steps. We
repeat this eight times to increase the statistical accuracy and estimate
error bars. To compare to the theory (\ref{eq:S_c_c}), we set the
concentration gradient to $h_{\parallel}=0.39/0.25\mbox{ cm}^{-1}$
and evaluate $\rho\approx1.05\frac{\textrm{g}}{\textrm{cm}^{3}}$
at $c=0.39/2$ from the equation of state. When computing the theory,
we account for errors in the discrete approximation to the continuum
Laplacian by using the effective wavenumber
\begin{equation}
k_{\perp}=k_{x}\frac{\sin\left(k_{x}\Delta x/2\right)}{\left(k_{x}\Delta x/2\right)}\label{eq:modified_kx}
\end{equation}
instead of the actual discrete wavenumber $k_{x}$ \cite{LLNS_Staggered}.

\begin{figure}
\begin{centering}
\includegraphics[width=0.75\textwidth]{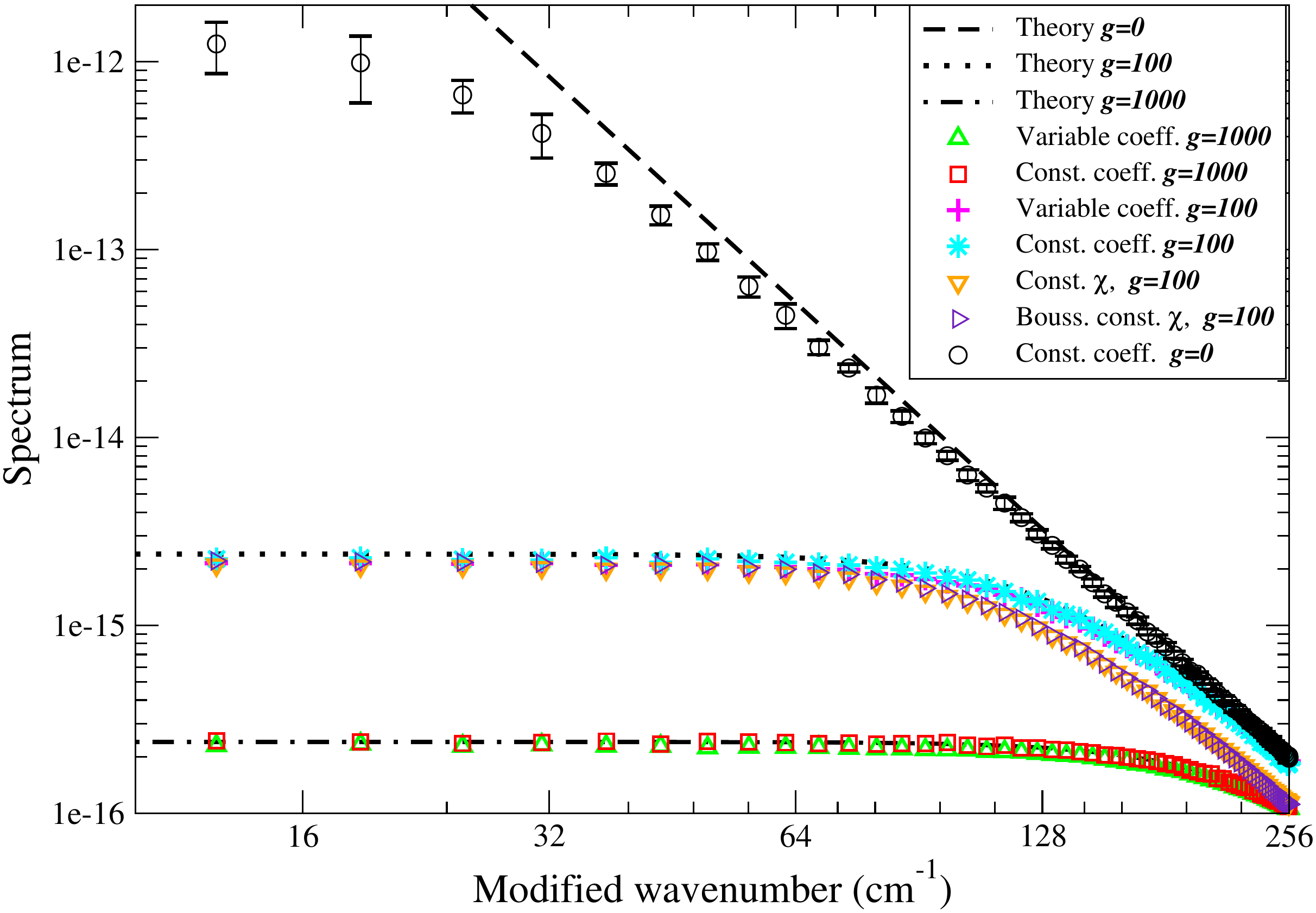}
\par\end{centering}

\caption{\label{fig:ConstantGradient}Comparison between the simple theory
(\ref{eq:S_c_c_simp}) (lines) and numerical results (symbols). Results
are shown for standard gravity $g\approx10^{3}\,\frac{\textrm{cm}}{\textrm{s}^{2}}$
(the cutoff wavenumber $k_{g}\approx246\mbox{ cm}^{-1}$), for the
complete variable-coefficient variable-density low Mach model (green
upward triangles) and the constant-coefficient constant-density approximation
(red squares). Also shown are results for a weaker gravity, $g\approx10^{2}\,\frac{\textrm{cm}}{\textrm{s}^{2}}$
(the cutoff wavenumber $k_{g}\approx138\mbox{ cm}^{-1}$), for the
complete low Mach model (magenta pluses) and the constant-coefficient
constant-density approximation (cyan stars). For comparison, results
for $g\approx10^{2}\,\frac{\textrm{cm}}{\textrm{s}^{2}}$ with variable
viscosity $\eta(c)$ but constant diffusion coefficient $\chi(c)=\chi_{0}$
are also shown, for variable density (orange downward triangles) and
the constant-density (Boussinesq) approximation (indigo right-facing
triangles). Finally, results for no gravity are shown in the constant-coefficient
approximation (black circles).}
\end{figure}

The results for the static spectrum of concentration fluctuations
$S_{c,c}\left(k_{x},k_{y}=0\right)$ as a function of the modified
wavenumber $k_{x}$ (\ref{eq:modified_kx}) are shown in Fig. \ref{fig:ConstantGradient}.
When there is no gravity, we see the characteristic giant fluctuation
power-law spectrum of the fluctuations, modulated at small wavenumbers
due to the presence of the physical boundaries \cite{LLNS_Staggered}.
When gravity is present, fluctuations at wavenumber below the cutoff
$k_{g}=\left[h_{\parallel}\rho g\beta/\left(\eta\chi\right)\right]^{1/4}$
are suppressed. If we use a constant-coefficient approximation, in
which we reduce $\beta=0.01$ so that $\rho\approx\rho\left(\bar{c}\right)$
and also fix the transport coefficients at $\eta(c)=\eta\left(\bar{c}\right)$
and $\chi(c)=\chi(\bar{c})$, we observe good agreement with the quasi-periodic
theory (\ref{eq:S_c_c_simp}). When we make the transport coefficients
dependent on the concentration as in (\ref{eq:nu_c},\ref{eq:chi_c}),
we observe a rather small change in the spectrum. This is perhaps
not unexpected because the simplified theory (\ref{eq:S_c_c_simp})
shows that only the product $\chi\eta$, and not $\chi$ and $\eta$
individually, matters. Since we used the Stokes-Einstein relation
$\chi(c)\eta(c)=\rho_{0}\chi_{0}\nu_{0}=\text{const}.$ to select
the concentration dependence of the diffusion coefficient, the value
of $\chi\eta$ is constant throughout the physical domain. For comparison,
in Fig. \ref{fig:ConstantGradient} we show results from a simulation
where we keep the concentration dependence of the viscosity (\ref{eq:nu_c})
but set the diffusion coefficient to a constant value, $\chi(c)=\chi_{0}$,
and we observe a more significant change in the spectrum. Further
employing the Boussinesq approximation makes little difference showing
that the primary effect here comes from the dependence of the transport
coefficients on concentration.

This shows that under the sort of parameters present in the experiments
on diffusive mixing in water-glycerol mixture, it is reasonable to
make the Boussinesq incompressible approximation; however, the spatial
dependence of the viscosity and diffusion coefficient cannot in general
be ignored if quantitative agreement is desired. In particular, time-dependent
quantities such as dynamic spectra \cite{GiantFluctuations_Theory,SoretDiffusion_Croccolo}
depend on the individual values of $\chi$ and $\eta$ and not just
their product, and are thus expected to be more sensitive to the details
of their concentration dependence. Even though the constant-coefficient
approximation gives qualitatively the correct shape and a better choice
of the constant transport coefficients may improve its accuracy, there
is no obvious or simple procedure to \emph{a priori} estimate what
parameters should be used (but see \cite{GiantFluctuations_Theory}
for a proposal to average the constant-coefficient theory over the
domain). A direct comparison with experimental results is not possible
until multiscale temporal integrators capable of handling the extreme
separation of time scales between mass and momentum diffusion are
developed. At present this has only been accomplished in the constant-coefficient
incompressible limit ($\beta=0$) \cite{StokesLaw}, and it remains
a significant challenge to accomplish the same for the complete low
Mach number system.

\section{\label{sec:MixingMD}Diffusive Mixing in Hard-Disk and Hard-Sphere
Fluids}

In this section we study the appearance of giant fluctuations during
\emph{time-dependent} diffusive mixing. As a validation of the low
Mach number fluctuating equations and our algorithm, we perform simulations
of diffusive mixing of two fluids of different densities in two dimensions.
We find excellent agreement between the results of low Mach number
(continuum) simulations and hard-disk molecular dynamics (particle)
simulations. This nontrivial test clearly demonstrates the usefulness
of low Mach number models as a coarse-grained mesoscopic model for
problems where sound waves can be neglected.

\begin{figure}
\begin{centering}
\includegraphics[width=0.4\textwidth]{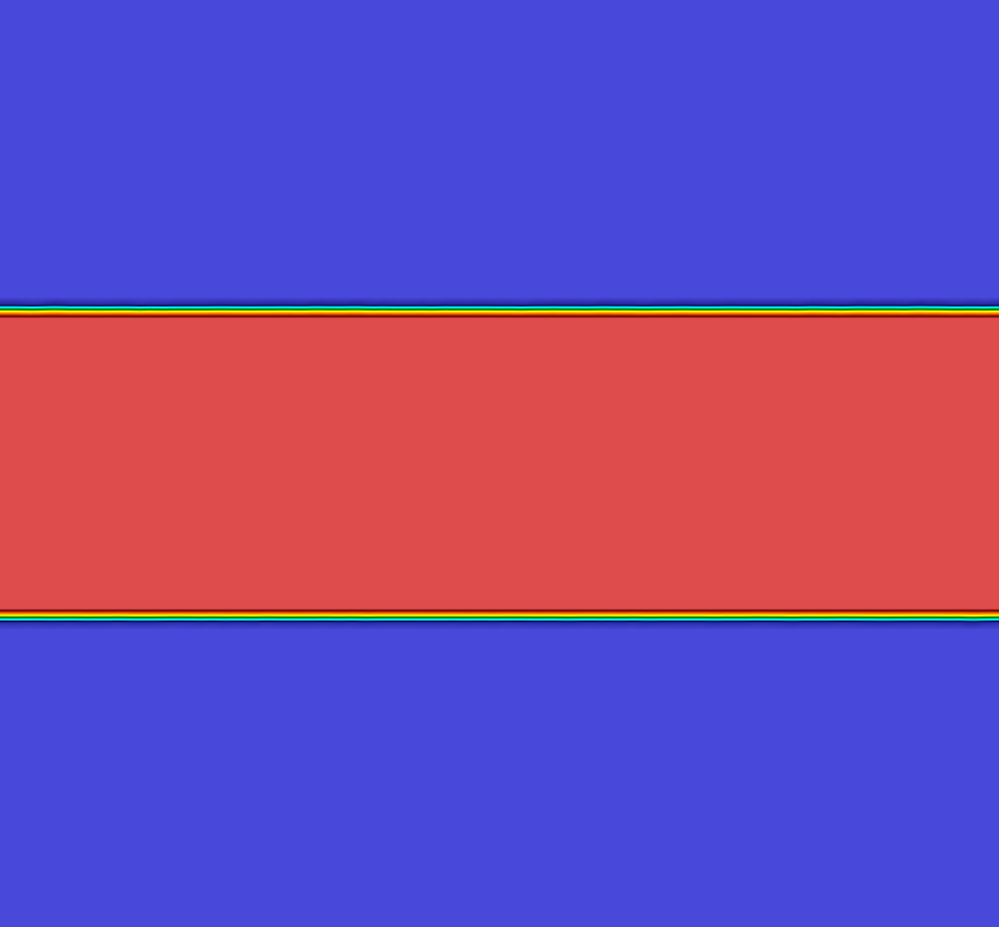}\hspace{0.2cm}\includegraphics[width=0.4\textwidth]{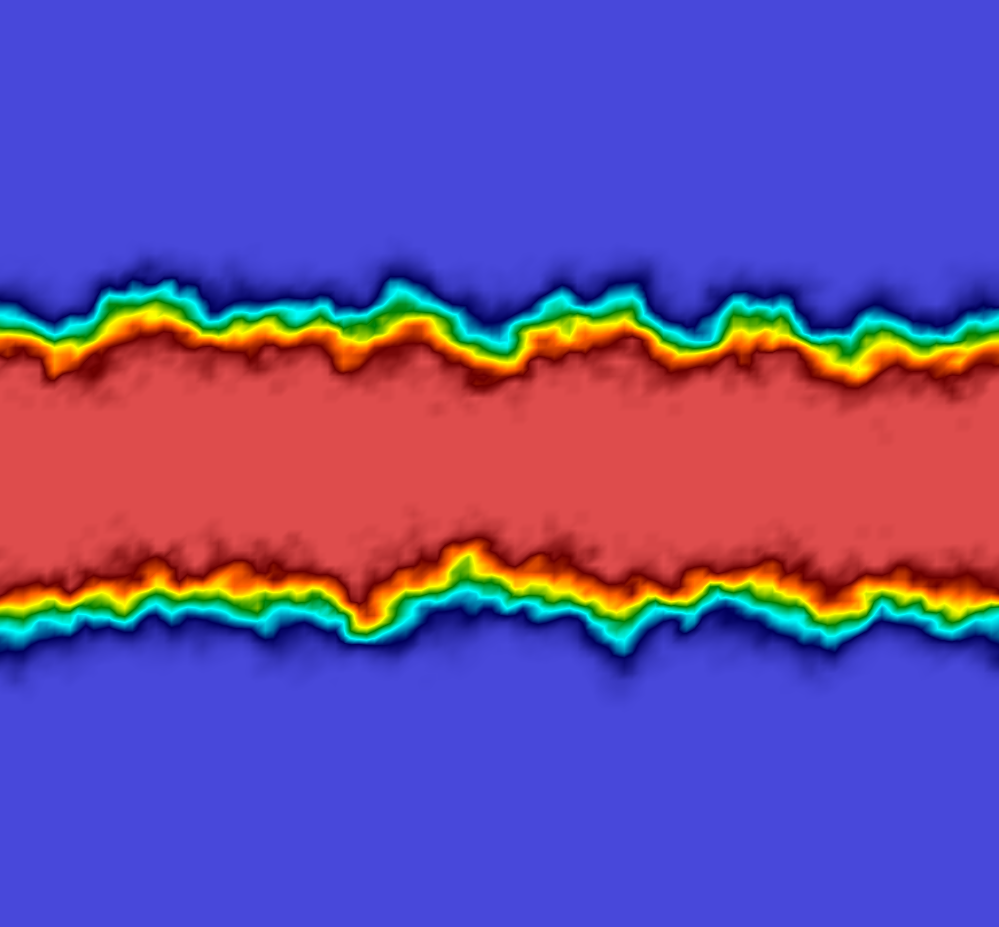}\vspace{0.2cm}
\par\end{centering}

\begin{centering}
\includegraphics[width=0.4\textwidth]{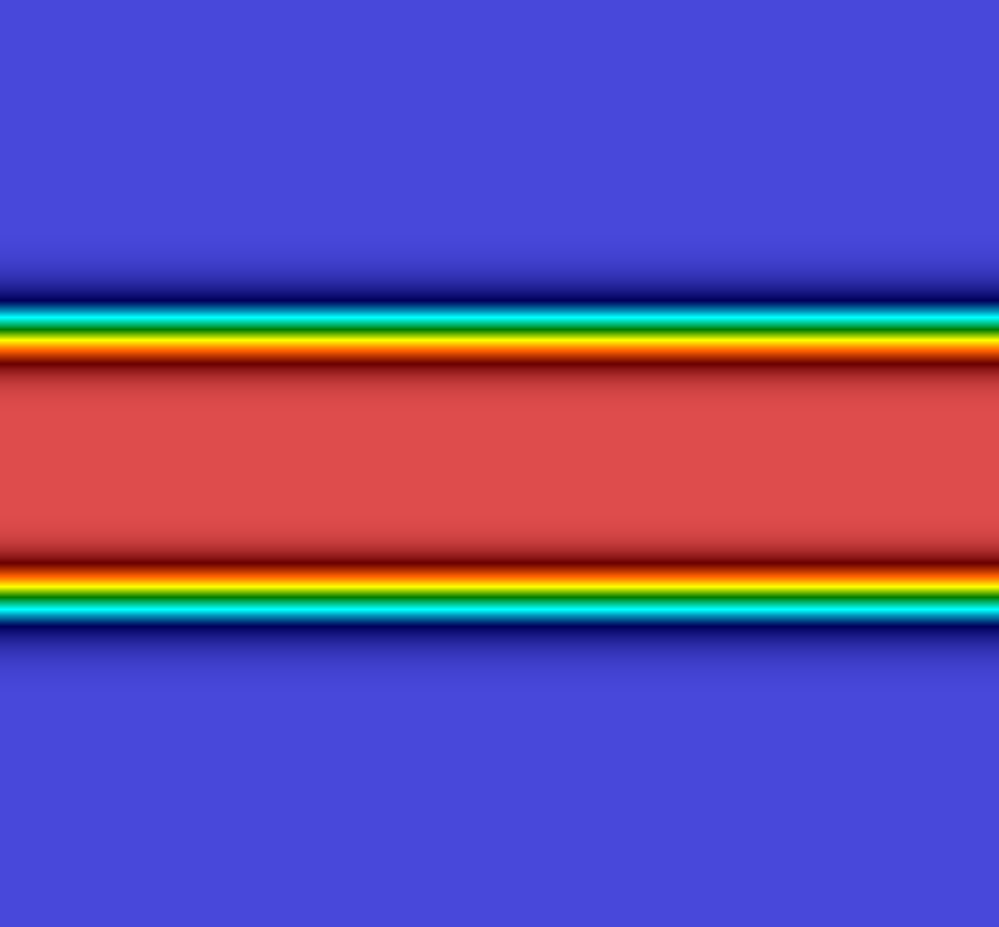}\hspace{0.2cm}\includegraphics[width=0.4\textwidth]{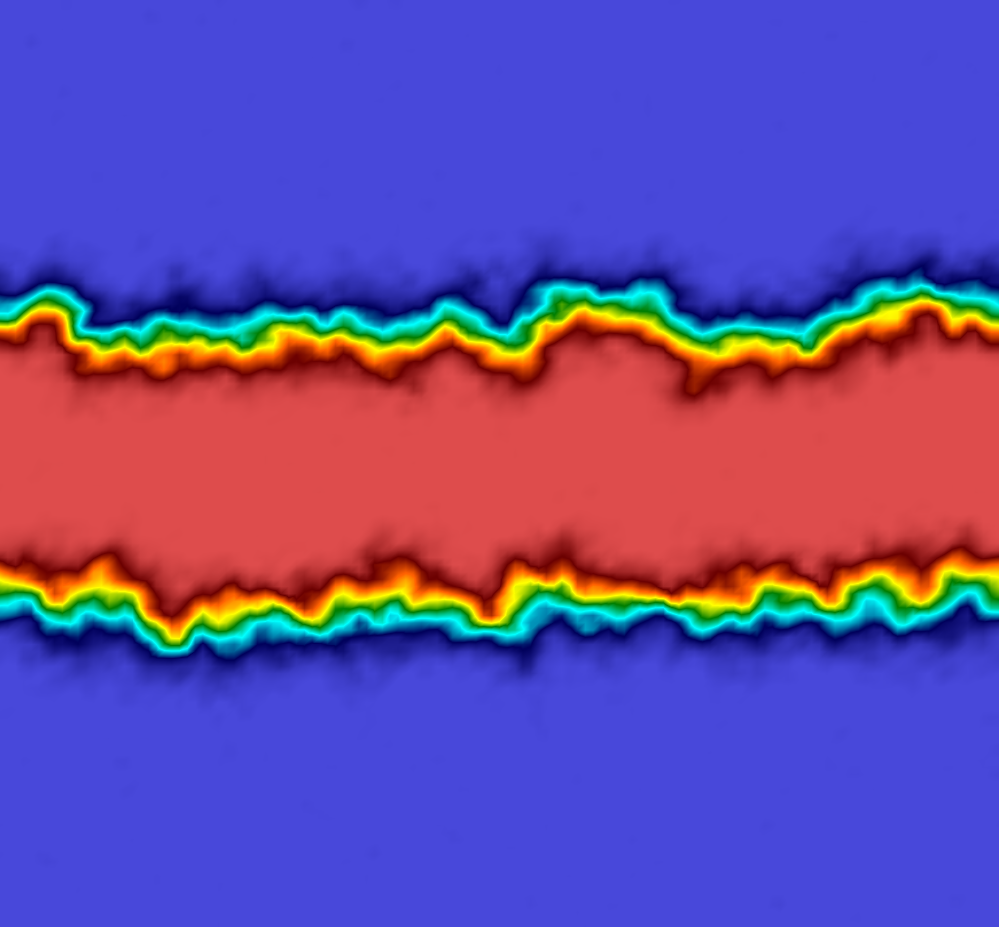}
\par\end{centering}

\caption{\label{fig:MixingIllustrationMass4}Diffusive mixing between two fluids
of unequal densities, $R=\rho_{2}/\rho_{1}$=4, with coloring based
on concentration, red for the pure first component, $c=1$, and blue
for the pure second component, $c=0$. A smoothed shading is used
for the coloring to eliminate visual discretization artifacts. The
simulation domain is periodic and contains $128^{2}$ hydrodynamic
(finite volume) cells. The top left panel shows the initial configuration,
which is the same for all simulations reported here. The top right
panel shows the final configuration at time $t=5,800$ as obtained
using molecular dynamics. The bottom left panel shows the final configuration
obtained using deterministic hydrodynamics, while the right panel
shows the final configuration obtained using fluctuating hydrodynamics.}
\end{figure}

Our simulation setup is illustrated in Fig. \ref{fig:MixingIllustrationMass4}.
We consider a periodic square box of length $L$ along both the $x$
(horizontal) and $y$ (vertical) directions, and initially place all
of the fluid of species one (colored red) in the middle third of the
domain, i.e., we set $c=1$ for $L/3\leq y\leq2L/3$, and $c=0$ otherwise,
as shown in the top left panel of the figure. The two fluids mix diffusively
and at the end of the simulation the concentration field shows a \emph{rough
diffusive interface} as confirmed by molecular dynamics simulations
shown in the top right panel of the figure. The deterministic equations
of diffusive mixing reduce to a one dimensional model due to the translational
symmetry along the $x$ axes, and would yield a \emph{flat} diffusive
interface as illustrated in the bottom left panel of the figure. However,
fluctuating hydrodynamics correctly reproduces the interface roughness,
as illustrated in the bottom right panel of the figure and demonstrated
quantitatively below.

We consider here a binary hard-disk mixture in two dimensions. We
use arbitrary (molecular) units of length, time and mass for convenience.
All hard disks had a diameter $\sigma=1$ in arbitrary units, and
we set the temperature at $k_{B}T=1$. The molecular mass for the
first fluid component was fixed at $m_{1}=1$, and for the second
component at $m_{2}=Rm_{1}$. For mass ratio $R=1$, the two types
of disks are mechanically-identical and therefore the species label
is simply a red-blue coloring of the particles. In this case $\bar{\rho}_{2}=\bar{\rho}_{1}$
and the low Mach number equations reduce to the incompressible equations
of fluctuating hydrodynamics with a passively-advected concentration
field. For the case of unequal particle masses, mechanical equilibrium
is obtained if the pressures in the two fluid components are the same.
It is well-known from statistical mechanics that for hard disks or
hard spheres the pressure is
\[
P=Y\left(\phi\right)\cdot n\cdot k_{B}T
\]
where $n=N/V$ is the number density, and $Y\left(\phi\right)$ is
a prefactor that only depends on the packing fraction $\phi=n\left(\pi\sigma^{2}/4\right)$
and not on the molecular mass. Therefore, for a mixture of disks or
spheres with equal diameters, at constant pressure, the number density
and the packing fraction $\phi$ are constant independent of the composition.
The equation of state at constant pressure and temperature is therefore
\[
1=\frac{n_{1}}{n}+\frac{n_{2}}{n}=\frac{\rho_{1}}{nm_{1}}+\frac{\rho_{2}}{nm_{2}},
\]
which is exactly of the form (\ref{eq:EOS_quasi_incomp}) with $\bar{\rho}_{1}=nm_{1}$
and $\bar{\rho}_{2}=nm_{2}$. The chemical potential of such a mixture
has the same concentration dependence as a low-density gas mixture
\cite{Landau:StatPhys1},
\[
\mu_{c}^{-1}k_{B}T=c\left(1-c\right)\left[cm_{2}+\left(1-c\right)m_{1}\right].
\]

\subsection{Hard Disk Molecular Dynamics}

In order to validate the predictions of our low Mach number model,
we performed Hard Disk Molecular Dynamics (HDMD) simulations of diffusive
mixing using a modification of the public-domain code developed by
the authors of Ref. \cite{MRJ_HS_4D}. We used a packing fraction
of $\phi=0.6$ for all simulations reported here. This packing fraction
is close to the freezing transition point but is known to be safely
in the (dense) gas phase (there is no liquid phase for a hard-disk
fluid). The initial particle positions were generated using a nonequilibrium
molecular dynamics simulation as in the hard-particle packing algorithm
described in Ref. \cite{Event_Driven_HE}. After the initial configuration
was generated the disks were assigned a species according to their
$y$ coordinate, and the mixing simulation performed using event-driven
molecular dynamics.

In order to convert the particle data to hydrodynamic data comparable
to that generated by the fluctuating hydrodynamics simulations, we
employed a grid of $N_{c}^{2}$ hydrodynamic cells that were each
a square of linear dimension $L_{c}=10\sigma$. At the chosen packing
fraction $\phi=0.6$ this corresponds to about $76$ disks per hydrodynamic
cell, which is deemed a reasonable level of coarse-graining for the
equations of fluctuating hydrodynamics to be a reasonably-accurate
model, while still keeping the computational demands of the simulations
manageable. We performed HDMD simulations for systems of size $N_{c}=64$
and $N_{c}=128$ cells, and simulated the mixing process to a final
simulation time of $t=5,800$ units. The largest system simulated
had about $1.25$ million disks (each simulation took about 5 days
of CPU time), which is well into the ``hydrodynamic'' rather than
``molecular'' scale.

Every $58$ units of time, particle data was converted to hydrodynamic
data for the purposes of analysis and comparison to hydrodynamic calculations.
There is not a unique way of coarse-graining particle data to hydrodynamic
data \cite{StagerredFluctHydro,DiscreteDiffusion_Espanol}; however,
we believe that the large-scale (giant) concentration fluctuations
studied here are \emph{not} affected by the particular choice. We
therefore used a simple method consistent with the philosophy of finite-volume
conservative discretizations. Specifically, we coarse-grained the
particle information by sorting the particles into hydrodynamic cells
based on the position of their centroid, as if they were point particles.
We then calculated $\rho_{1}$ and $\rho_{2}$ in each cell based
on the total mass of each species contained inside the given cell.
Since all particles have equal diameter other definitions that take
into account the particle shape and size give similar results.

\subsection{Hard Disk Hydrodynamics}

We now turn to hydrodynamic simulations of the diffusive mixing of
hard disks. Our hydrodynamic calculations use the same grid of cells
used to convert particle to hydrodynamic data. The only input required
for the hydrodynamic calculations, in addition to those provided by
equilibrium statistical mechanics, are the transport coefficients
of the fluid as a function of concentration, specifically, the shear
viscosity $\eta$ and the diffusion coefficient $\chi$.

The values for the transport coefficients used in the spatio-temporal
discretization, as explained in Refs. \cite{DiffusionRenormalization,StokesLaw}
and detailed in Appendix \ref{sec:TransportMD}, are not material
constants independent of the discretization. Rather, they are \emph{bare}
transport values $\eta_{0}$ and $\chi_{0}$ measured at the length
scales of the grid size. We assumed that the bare transport coefficients
obey the same scaling with the mass ratio $R$ as predicted by Enskog
kinetic theory (\ref{eq:eta_R_scaling},\ref{eq:chi_R_scaling}).
As explained in Appendix \ref{sec:TransportMD}, theoretical arguments
and molecular dynamics results suggest that renormalization effects
for viscosity are small and can be safely neglected. We have therefore
fixed the viscosity in the hydrodynamic calculations based on the
molecular dynamics estimate $\eta_{0}=2.5$ for the pure fluid with
molecular mass $m=1$ (see Section \ref{sub:Viscosity}). However,
the bare diffusion coefficient is strongly dependent on the size of
the hydrodynamic cells (held fixed in our calculations at $\D x=\D y=10$),
and on whether filtering (see Appendix \ref{AppendixFiltering}) is
used. Therefore, the value of $\chi_{0}$ needs to be adjusted based
on the spatial discretization, in such a way as to match the behavior
of the molecular dynamics simulations at length scales much larger
than the grid spacing. We describe the exact procedure we used to
accomplish this in Section \ref{sub:DiffusionCoeffient}.

The time step in our explicit algorithm is limited by the viscous
CFL number $\alpha_{\nu}=\nu\D t/\D x^{2}<1/4$. Since the hydrodynamic
calculations are much faster compared to the particle simulations,
we used the more expensive RK3 temporal integrator with a relatively
small time step $\D t=1.45$, corresponding to $\alpha_{\nu}\approx0.05$
for $c=1$. For $R=1$ and $N_{c}=64$ we employed a larger time step,
$\D t=3.625$ ($\alpha_{\nu}\approx0.125$), with no measurable temporal
discretization artifacts for the quantities studied here. We are therefore
confident that the discretization errors in this study are dominated
by spatial discretization artifacts. In future work we will explore
semi-implicit discretizations and study the effect of taking larger
time steps on temporal accuracy. Note that at these parameters for
$c=1$ the isothermal speed of sound is $c_{T}\approx5.1$ so that
a compressible scheme would require a time step on the order of $\D t\sim1$
(corresponding to advective CFL of about a half). By contrast, the
explicit low Mach number algorithm is stable for $\D t\lesssim7.5$.
This modest gain is due to the small hydrodynamic cell we use here
in order to compare to molecular dynamics. For mesoscopic hydrodynamic
cells the gain in time step size afforded by the low Mach formulation
will be several orders of magnitude larger.

For mass ratio $R=1$ and $R=2$, the hydrodynamic calculations were
initialized using statistically identical configurations as would
be obtained by coarse-graining the initial particle configuration.
This implies a sharp, step-like jump in concentration at $y=L/3$
and $y=2L/3$. Since our spatio-temporal discretization is not strictly
monotonicity-preserving, such sharp concentration gradients combined
with a small diffusion coefficient $\chi_{0}$ lead to a large cell
Peclet number. This may in turn lead to large deviations of concentration
outside of the allowed interval $0\leq c\leq1$ for larger mass ratios.
Therefore, for $R=4$ we smoothed the initial condition slightly so
that the sharp jump in concentration is spread over a few cells, and
also employed a $9$ point filter for the advection velocity ($w_{F}=4$,
see Appendix \ref{AppendixFiltering}). We verified that for $R=2$
using filtering only affects the large wavenumbers and does not appear
to affect the small wavenumbers we study here, provided the bare diffusion
coefficient $\chi_{0}$ is adjusted based on the specific filtering
width $w_{F}$.

\subsection{Comparison between Molecular Dynamics and Fluctuting Hydrodynamics
simulations}

In order to compare the molecular dynamics and the hydrodynamic simulations
we calculated several statistical quantities:
\begin{enumerate}
\item The averages of $\rho_{1}$ along the directions perpendicular to
the concentration gradient,
\begin{equation}
\rho_{1}^{\left(h\right)}\left(y\right)=L^{-1}\int_{x=0}^{L}\rho_{1}\left(x,y\right)dx,\label{eq:rho_1_h}
\end{equation}
where the integral is discretized as a direct sum over the hydrodynamic
cells. Note that it is statistically better to use conserved quantities
for such macroscopic averages than to use non-conserved variables
such as concentration \cite{UnbiasedEstimates_Garcia}.
\item The spectrum of the concentration averaged along the direction of
the gradient by computing the average
\[
c_{v}\left(x\right)=L^{-1}\int_{y=0}^{L}c\left(x,y\right)dy,
\]
and then taking the discrete Fourier transform. Intuitively, $c_{v}$
is a measure of the thickness of the red strip in Fig. \ref{fig:MixingIllustrationMass4},
and corresponds closely to what is measured in light scattering and
shadowgraphy experiments \cite{FluctHydroNonEq_Book,GiantFluctuations_Cannell}.
\item The discrete Fourier spectrum of the $y$-coordinate of the ``center-of-mass''
of concentration along the direction perpendicular to the gradient,
\[
h_{c}\left(x\right)=L^{-1}\int_{y=0}^{L}y\cdot c\left(x,y\right)dy.
\]
Intuitively, $h_{c}$ is a measure of the height of the centerline
of the red strip in Fig. \ref{fig:MixingIllustrationMass4}.
\end{enumerate}
All quantities were sampled at certain pre-specified time points in
a number of statistically-independent simulations $N_{s}$ and then
means and standard deviations calculated from the $N_{s}$ data points.
For system of size $N_{c}=64$ cells we used $N_{s}=64$ simulations,
and for systems of size $N_{c}=128$ we used $N_{s}=32$ simulations.
By far the majority of the computational cost was in performing the
HDMD simulations.

\subsubsection{Average Concentration Profiles}

Once $\chi_{0}$ and $\chi_{\text{eff}}$ were estimated based on
simulations of a constant-density ($R=1$) fluid (see Section \ref{sub:DiffusionCoeffient}),
kinetic theory (\ref{eq:eta_R_scaling},\ref{eq:chi_R_scaling}) can
be used to estimate them for different density ratios. In Fig. \ref{fig:spreading_64x64}
we show $\rho_{1}^{\left(h\right)}\left(y\right)$ for mass ratio
$R=2$, showing good agreement between HDMD and hydrodynamics, especially
when fluctuations are accounted for. For $R=4$ a direct comparison
is difficult because the initial condition was slightly different
in the hydrodynamic simulations due to the need to smooth the sharp
concentration gradient for numerical reasons, as explained earlier.
This difference strongly affects the shape of $\rho_{1}^{\left(h\right)}\left(y\right)$
at early times, however, it does not significantly modify the roughness
of the interface, which we study next.

\subsubsection{Interface Roughness}

The most interesting contribution of fluctuations to the diffusive
mixing process is the appearance of giant concentration fluctuations
in the presence of large concentration gradients, as evidenced in
the roughness of the interface between the two fluids during the early
stages of the mixing in Fig. \ref{fig:MixingIllustrationMass4}. In
order to quantify this interface roughness we used the one-dimensional
power spectra
\[
S_{c}\left(k_{x}\right)=\av{\hat{c}_{v}\hat{c}_{v}^{\star}}\mbox{ and }S_{h}\left(k_{x}\right)=\av{\hat{h}_{c}\hat{h}_{c}^{\star}}.
\]
Note that here we do not correct the discrete wavenumber for the spatial
discretization artifacts and continue to use $k_{x}$ instead of $k_{\perp}$. 

\begin{figure}
\begin{centering}
\includegraphics[width=0.49\textwidth]{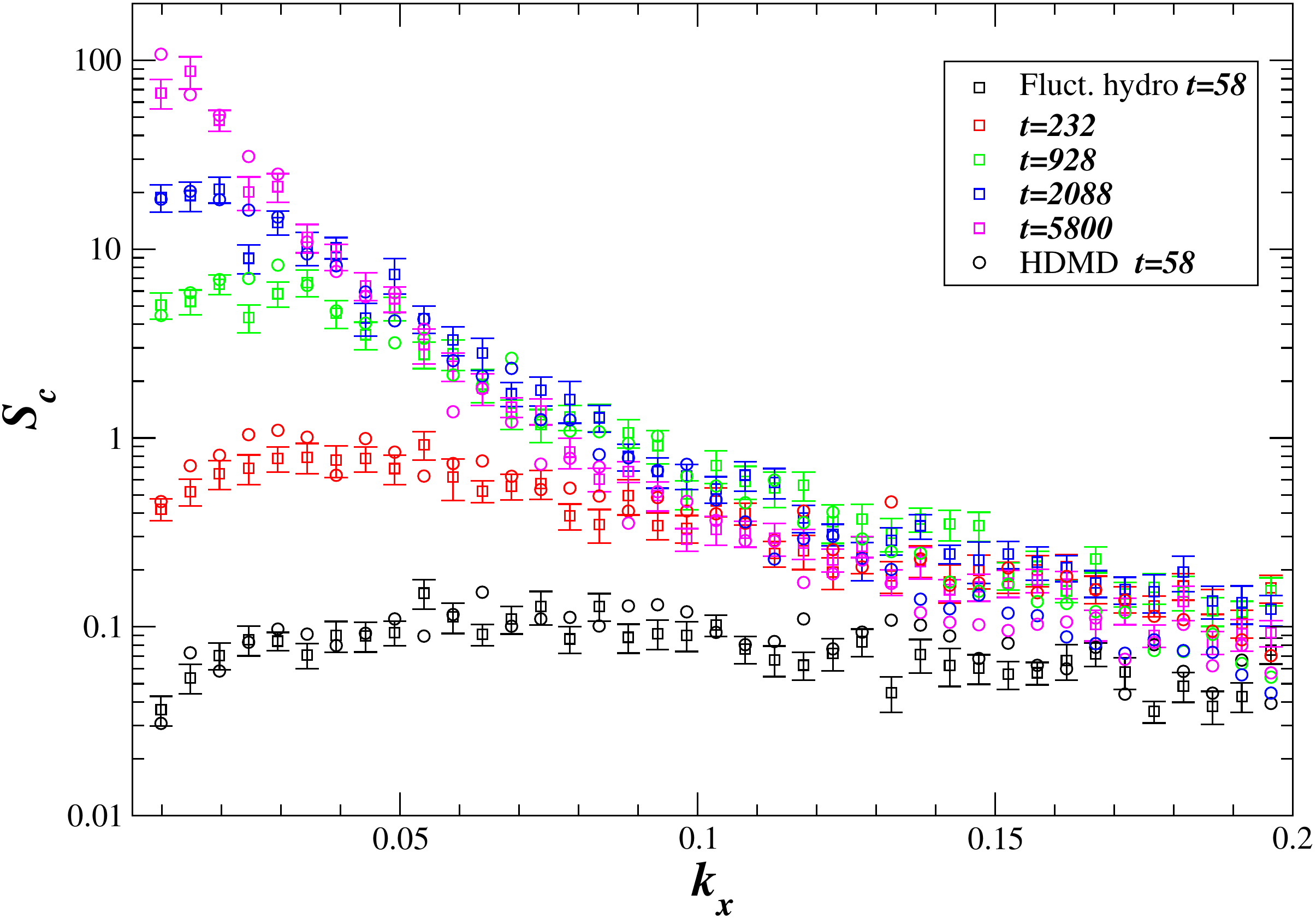}\includegraphics[width=0.49\textwidth]{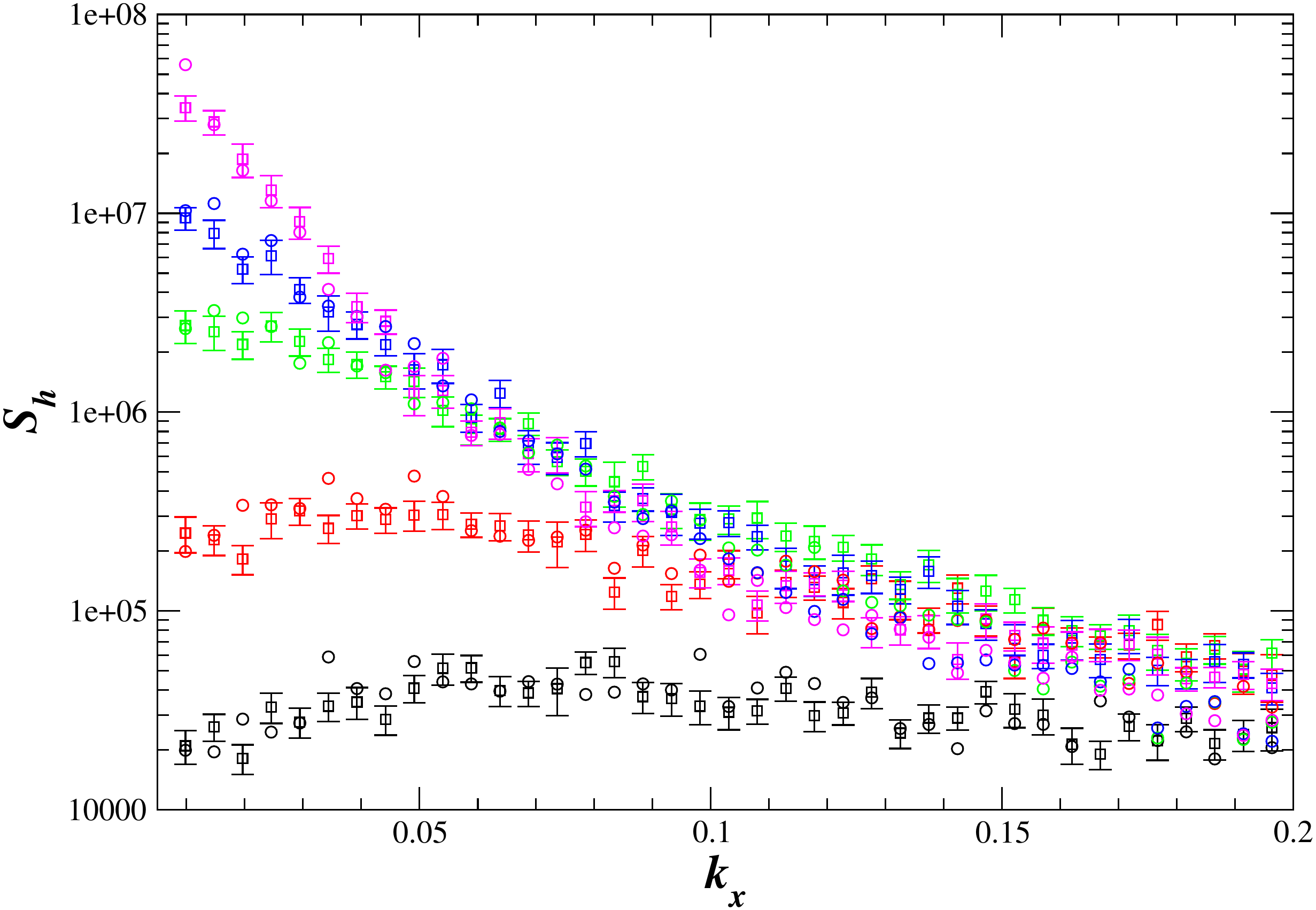}
\par\end{centering}

\caption{\label{fig:S_k_mass_1}Discrete spatial spectrum of the interface
fluctuations for $R=1$ and $N_{c}=128$ (averaged over 32 simulations)
at several points in time (drawn with different colors, as indicated
in the legend), for fluctuating hydrodynamics (FH, squares with error
bars) and HDMD (circles, error bars comparable to those for squares).
Note that the largest wavenumber supported by the grid is $k_{\max}=\pi/\D x\approx0.314$.
The larger wavenumbers are however dominated by spatial truncation
errors and the filter employed (if any) and we do not show them here.
(Left panel) Spectrum $S_{c}\left(k_{x}\right)$ of the vertically-averaged
concentration. (Right panel) Spectrum $S_{h}\left(k_{x}\right)$ of
the position of the vertical ``center-of-mass'' of concentration.}
\end{figure}

\begin{figure}
\begin{centering}
\includegraphics[width=0.49\textwidth]{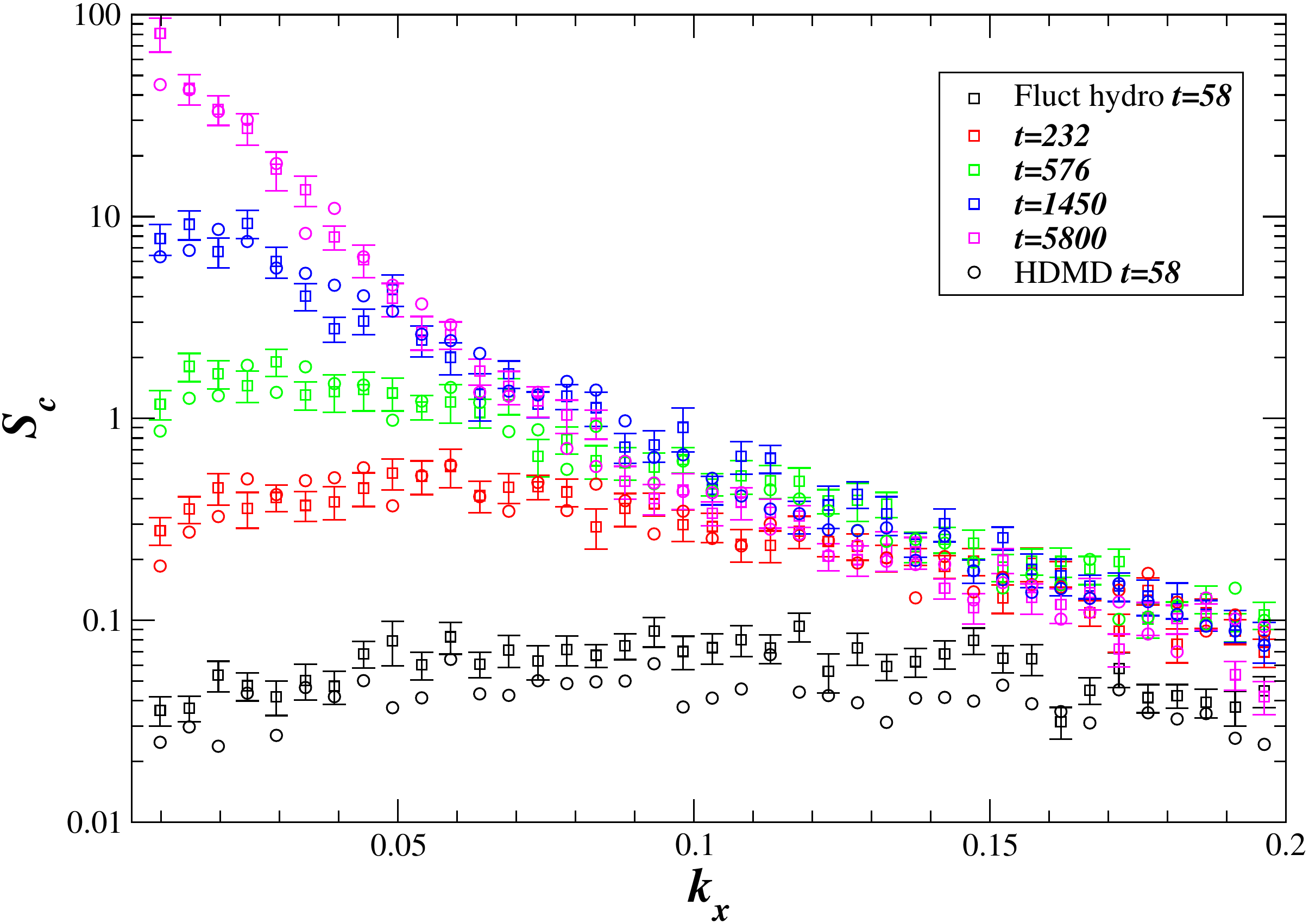}\includegraphics[width=0.49\textwidth]{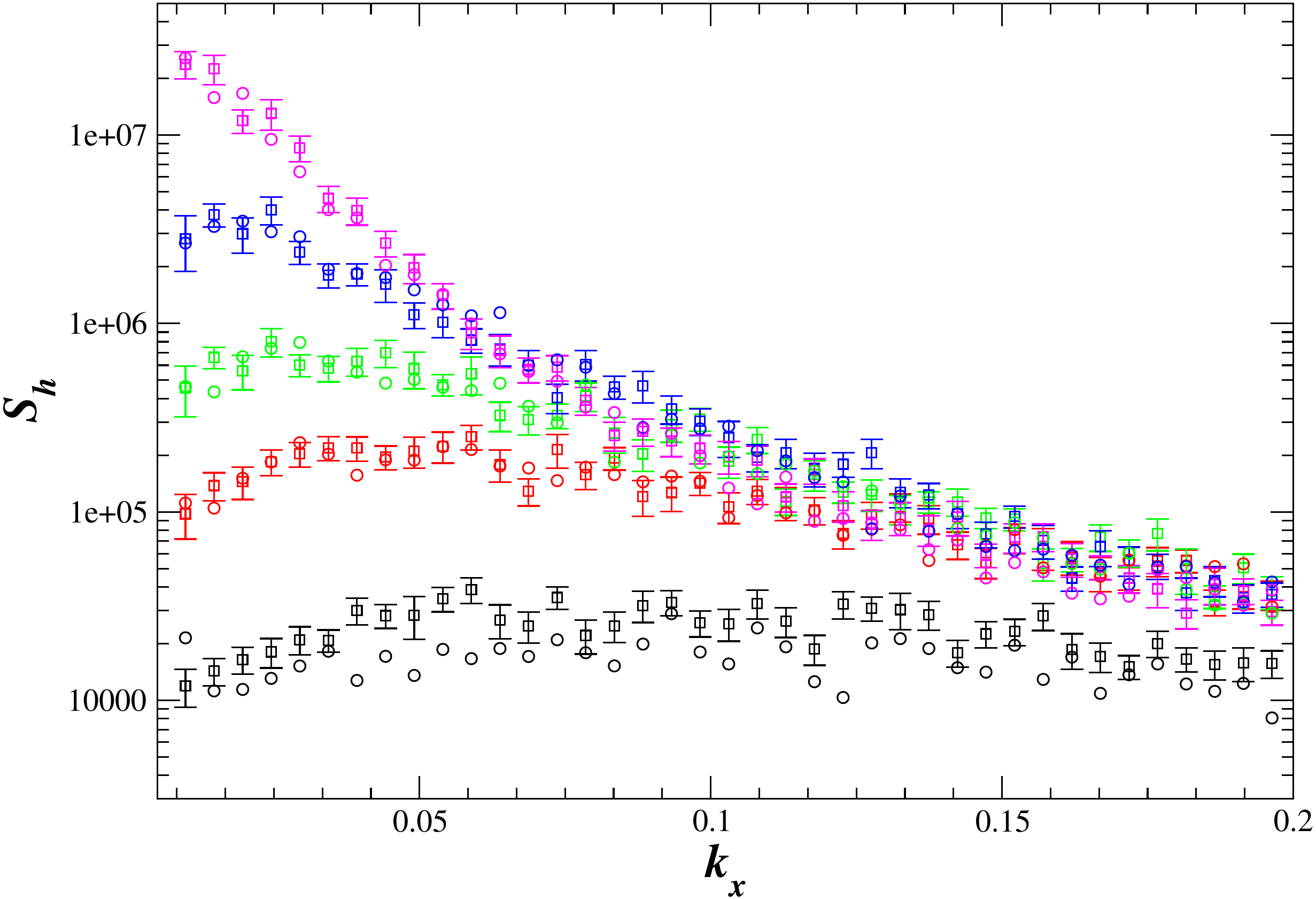}
\par\end{centering}

\caption{\label{fig:S_k_mass_4}Same as Fig. \ref{fig:S_k_mass_1} but for
density ratio $R=4$.}
\end{figure}

The temporal evolution of the spectra $S_{c}$ and $S_{h}$ is shown
in Fig. \ref{fig:S_k_mass_1} for mass ratio $R=1$, and in Fig. \ref{fig:S_k_mass_4}
for mass ratio $R=4$, for both HDMD and low Mach number fluctuating
hydrodynamics (note that deterministic hydrodynamics would give identically
zero for any spectral quantity). We observe an excellent agreement
between the two, including the correct initial evolution of the interface
fluctuations.

\begin{figure}
\begin{centering}
\includegraphics[width=0.49\textwidth]{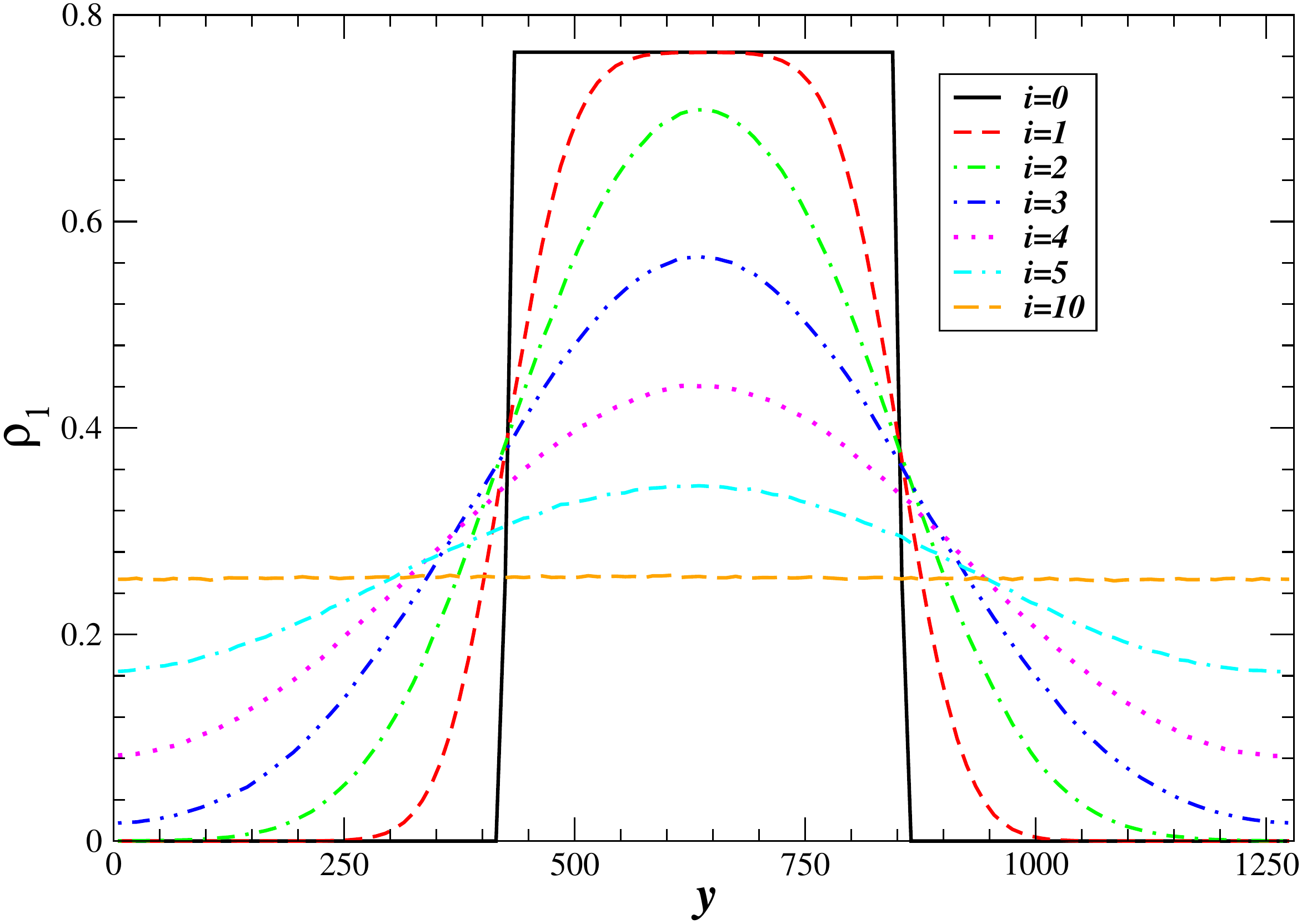}\includegraphics[width=0.49\textwidth]{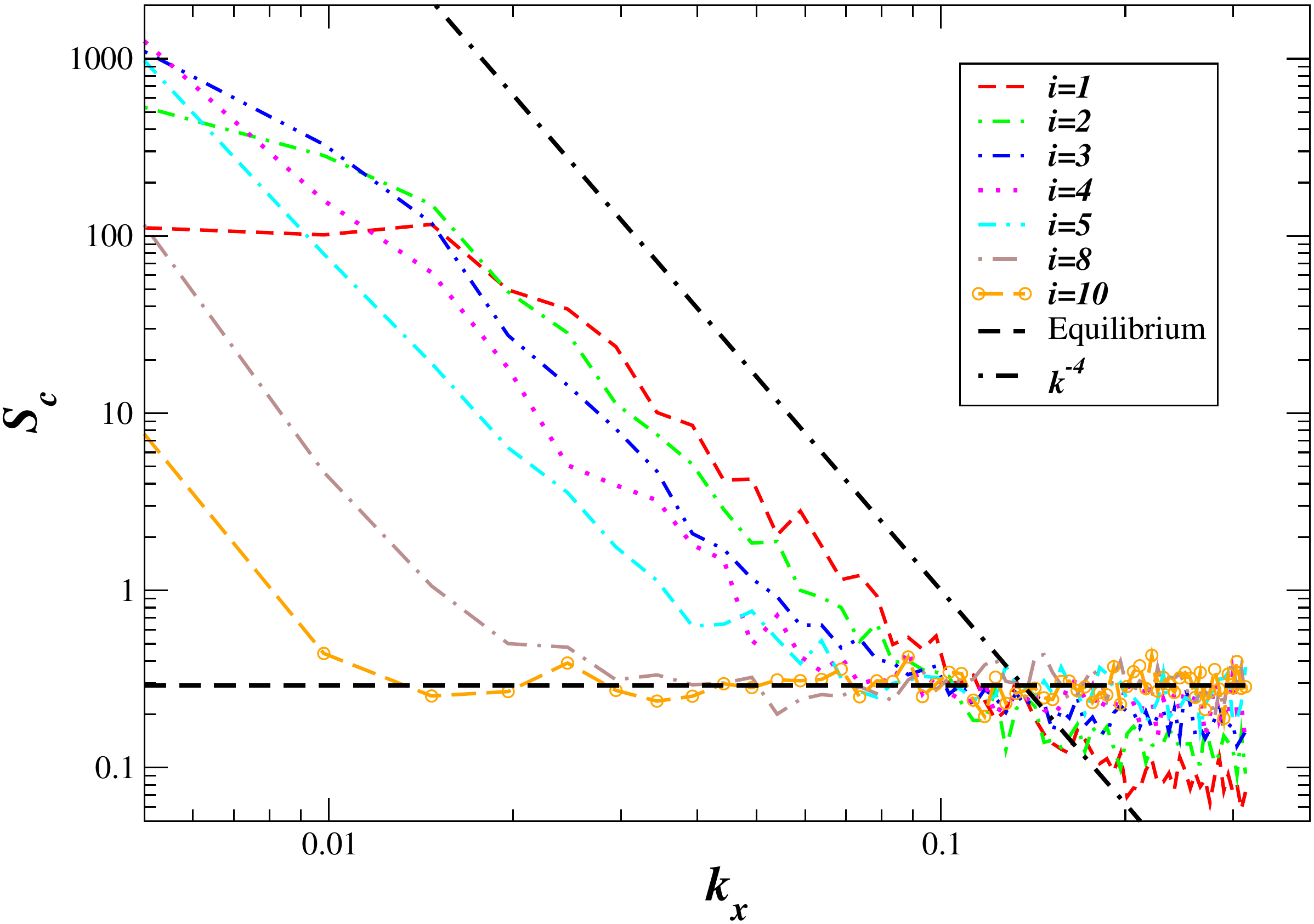}
\par\end{centering}

\caption{\label{fig:S_k_mass_1_long}Mixing to a time 128 times longer than
previous results, with results reported at time intervals $t=7424\, i^{2}$
for $i=1,\dots,10$. These long simulations are only feasible for
the fluctuating hydrodynamics code, and employ a somewhat larger time
step $\D t=3.625$. (\emph{Left}) Horizontally-averaged $\rho_{1}$,
as shown for the shorter runs in the left panel of Fig. \ref{fig:spreading_64x64}.
(\emph{Right}) The spectrum of interface fluctuations $S_{c}\left(k_{x}\right)$,
as shown in the left panels of Figs. \ref{fig:S_k_mass_1} and \ref{fig:S_k_mass_4}
for the shorter runs. The theoretical estimates for the spectrum of
equilibrium fluctuations, which is independent of wavenumber, is also
shown. We also indicate the theoretical prediction for the power-law
of the spectrum of steady-state nonequilibrium fluctuations under
an applied concentration gradient, $S_{c}\sim k^{-4}.$}
\end{figure}

Note that for a finite system, eventually complete mixing will take
place and the concentration fluctuations will have to revert to their
equilibrium spectrum, which is flat in Fourier space instead of the
power-law behavior seen out of equilibrium. In Fig. \ref{fig:S_k_mass_1_long}
we show results for mixing up to a time $t=7.42\cdot10^{5}$ (this
is 128 times longer than those described above). These long simulations
are only feasible for the fluctuating hydrodynamics code, and employ
a somewhat larger time step $\D t=3.625$. The results clearly show
that at late times the spectrum of the fluctuations reverts to the
equilibrium one; however, this takes some time even after the mixing
is essentially complete. Linearized incompressible fluctuating hydrodynamics
\cite{GiantFluctuations_Theory,LLNS_Staggered} predicts that at steady
state the spectrum of nonequilibrium concentration fluctuations is
a power law with exponent $-4$, $S_{c}\sim\left(\nabla c\right)^{2}k^{-4}.$
The dynamically-evolving spectra in the right panel of Fig. \ref{fig:S_k_mass_1_long}
show approximately such power-law behavior for intermediate times
and wavenumbers.

\subsection{Hard Sphere Fluctuating Hydrodynamics Simulations}

\begin{figure}
\begin{centering}
\includegraphics[width=0.75\textwidth]{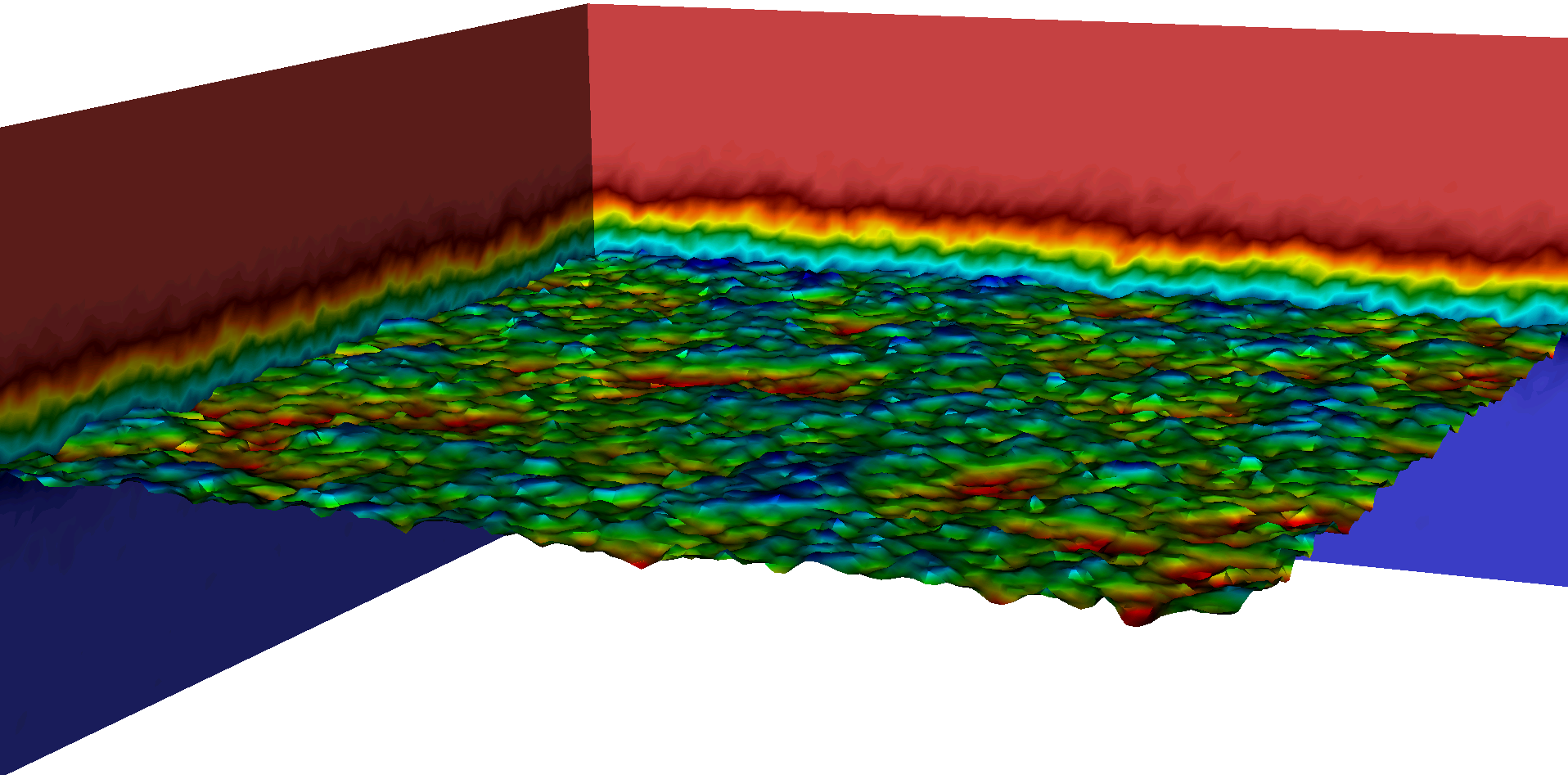}
\par\end{centering}

\begin{centering}
\includegraphics[width=0.49\textwidth]{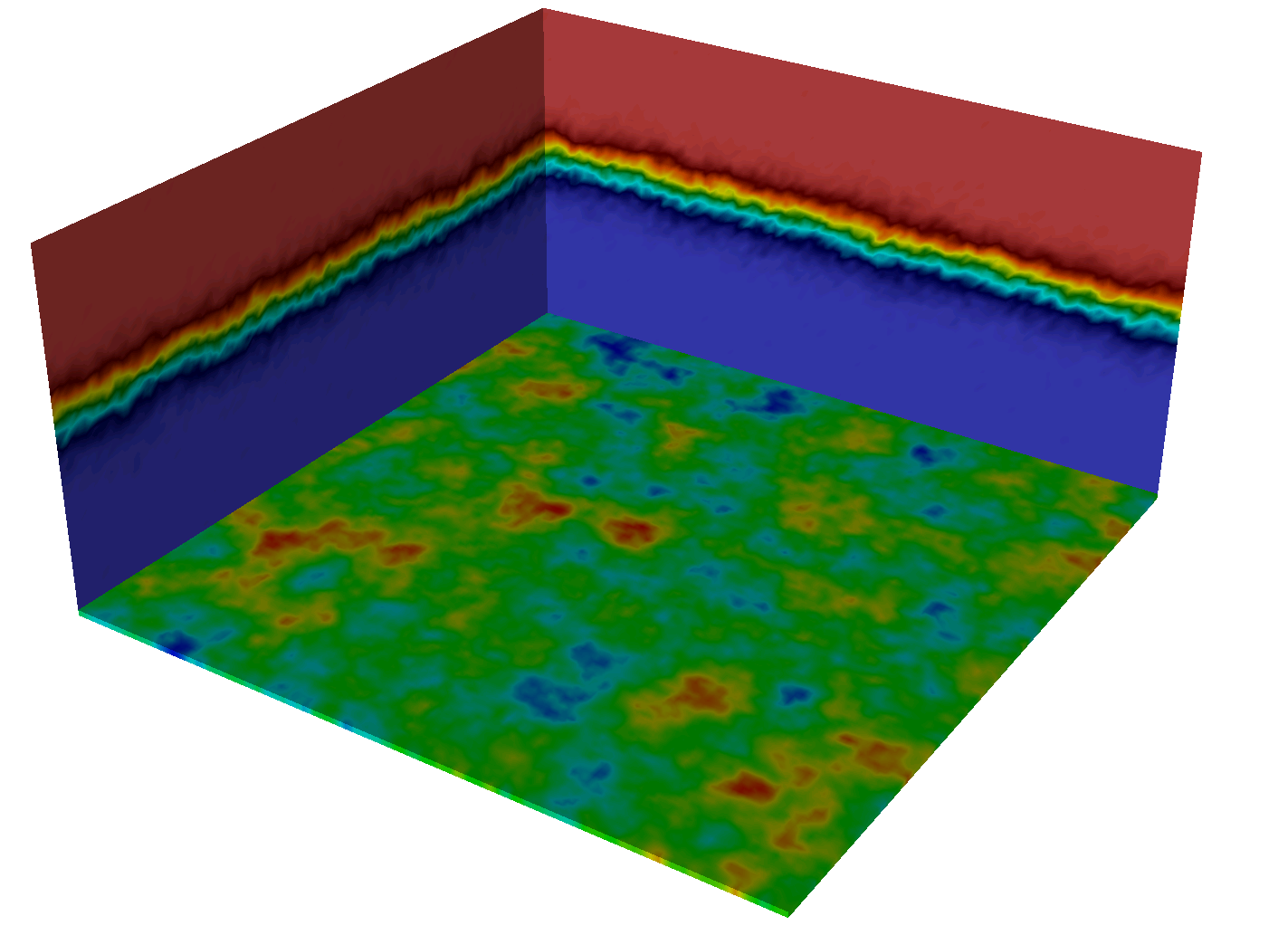}\includegraphics[width=0.49\textwidth]{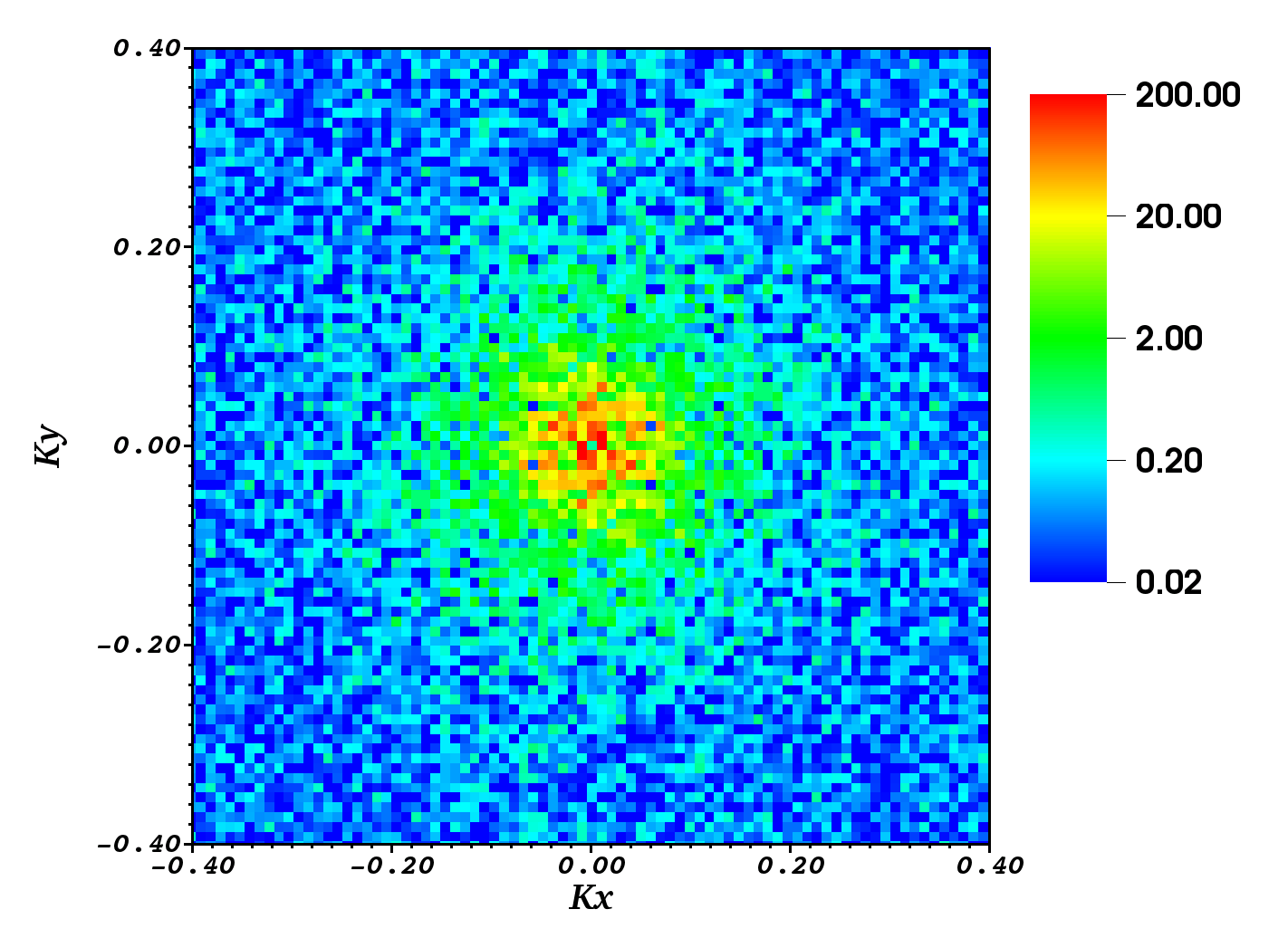}
\par\end{centering}

\caption{\label{fig:S_k_mass_4_3D}Diffusive mixing in three dimensions similar
to that illustrated in Fig. \ref{fig:MixingIllustrationMass4} for
two dimensions. Parameters are based on Enskog kinetic theory for
a hard-sphere fluid at packing fraction $\phi=0.45$, and there is
no gravity. The mixing starts with the top half being one species
and the bottom half another species, with density ratio $R=4$, and
concentration is kept fixed at the top and bottom boundaries while
the side boundaries are periodic. A snapshot taken at time $t=5,000$
is shown. (\emph{Top panel}) The side panes show two dimensional slices
for the concentration $c$. The approximated contour surface $c=0.2$
is shown with color based on surface height to illustrate the rough
diffusive interface. (Bottom left panel) Similar as top panel but
bottom pane shows vertically-averaged concentration $c_{v}\left(x,z\right)$,
illustrating the giant concentration fluctuations. (Bottom right panel)
The Fourier spectrum $S_{c}\left(k_{x},k_{y}\right)$ of $c_{v}$.
The color axes is logarithmic and clearly shows the appearance of
large scale (small wavenumber) fluctuations, as also seen in Fig.
\ref{fig:S_k_mass_4} in two dimensions.}
\end{figure}

In order to illustrate the appearance of giant fluctuations in three
dimensions we performed simulations of mixing in a mixture of hard
spheres with equal diameters, $\sigma=1$, and mass ratio $R=4$.
The packing density was chosen to be $\phi=0.45$, which corresponds
to a very dense gas, but is still well below the freezing point $\phi_{f}=0.49$.
For the hydrodynamic simulations we used cubic cells of dimension
$\D x=5$, which corresponds to about $107$ particles per hydrodynamic
cell on average. In Fig. \ref{fig:S_k_mass_4_3D} we show results
from a single simulation with a grid of size $128\times64\times128$
cells, which would correspond to about $10^{8}$ particles. This makes
molecular dynamics simulations infeasible, and makes hydrodynamic
calculations an invaluable tool in studying the mixing process at
these mesoscopic scales.

In the hydrodynamic simulations we used bare transport coefficient
values based on Enskog kinetic theory for the hard-sphere fluid \cite{HardSphereTransport_Review}.
For the single-component fluid with molecular mass $m=1$, this theory
gives $\eta_{0}\approx2.32$ and $\chi_{0}\approx0.053$, which corresponds
to a bare Schmidt number $S_{c}=\nu_{0}/\chi_{0}\approx51$. We employed
the same model dependence of bare transport coefficients on concentration
as for hard disks, see Eqs. (\ref{eq:eta_R_scaling},\ref{eq:chi_R_scaling}).
The time step was set at $\D t=1$ (corresponding to viscous CFL number
$\beta=\nu_{0}\D t/\D x^{2}\approx0.1$). In three dimensions, the
cell Peclet number is reduced with decreasing $\D x$ and we did not
find it necessary to employ any filtering. 

Instead of the fully periodic domain used in the two dimensional hard-disk
simulations, here we employ the fixed-concentration boundary conditions
(\ref{eq:v_n_BC}) and set $c(y=0;\, t)=0$ at the bottom and $c(y=L_{y};\, t)=1$
at the top boundary. This emulates the sort of ``open'' or ``reservoir''
boundaries \cite{OpenMD_Rafa} that mimic conditions in experimental
studies of diffusive mixing \cite{GiantFluctuations_Cannell}. The
initial condition is a fully phase-separated mixture with $c=1$ for
$y\geq L/2$, and $c=0$ otherwise. As the mixing process continues
the diffusive interface roughens and giant concentrations appear,
as illustrated in Fig. \ref{fig:S_k_mass_4_3D} and also observed
experimentally in water-glycerol mixtures in Ref. \cite{GiantFluctuations_Cannell}.
In three dimensions, however, the diffusive interface roughness is
much smaller than in two dimensions, being on the order of only 20
molecular diameters for the snapshot shown in the figure. This illustrates
the importance of dimensionality when including thermal fluctuations.
In particular, unlike in deterministic fluid dynamics, in fluctuating
hydrodynamics one cannot simply eliminate dimensions from consideration
even in simple geometries.

Approximate theory based on the Boussinesq approximation and linearization
of the equations of fluctuating hydrodynamics has been developed in
Ref. \cite{GiantFluctuations_Theory} and applied in the analysis
of experimental results on mixing in a water-glycerol mixture in the
presence of gravity \cite{GiantFluctuations_Cannell}. The simulations
reported here do not make the sort of approximations necessary in
analytical theories and can in principle be used to study the mixing
process quantitatively. However, it is important to emphasize that
in realistic liquids, such as a water-glycerol mixture, the Schmidt
number is on the order of a thousand. This makes explicit time stepping
schemes that fully resolve the dynamics of the velocity fluctuations
infeasible. In future work we will consider semi-implicit type stepping
methods that relax the severe time stepping restrictions present in
the explicit schemes considered here.

\section{\label{sec:Conclusions}Conclusions}

The behavior of fluids is strongly affected by thermal fluctuations
at scales from the microscopic to the macroscopic. Fluctuating hydrodynamics
is a powerful coarse-grained model for fluid dynamics at mesoscopic
and macroscopic scales, at both a theoretical and a computational
level. Theoretical calculations are rather complicated in the presence
of realistic spatial inhomogeneities and nontrivial boundary conditions.
In numerical simulations, those effects can readily be handled, however,
the large separation of time scales between different physical processes
poses a fundamental difficulty. Compressible fluctuating hydrodynamics
bridges the gap between molecular and hydrodynamic scales. At spatial
scales not much larger than molecular, sound and momentum and heat
diffusion occur at comparable time scales in both gases and liquids.
At mesoscopic and larger length scales, fast pressure fluctuations
due to thermally-actuated sound waves are much faster than diffusive
processes. It is therefore necessary to eliminate sound modes from
the compressible equations. In the deterministic context this is accomplished
using low Mach number asymptotic expansion.

For homogeneous simple fluids or mixtures of dynamically-identical
fluids the zeroth order low Mach equations are the well-known incompressible
Navier-Stokes equations, in which pressure is a Lagrange multiplier
enforcing a divergence-free velocity field. In mixtures of dissimilar
fluids, local changes in composition and temperature cause local expansion
and contraction of the fluid and thus a nonzero velocity divergence.
In this paper we proposed low Mach number fluctuating equations for
isothermal binary mixtures of incompressible fluids with different
density, or a mixture of low-density gases with different molecular
masses. These equations are a straightforward generalization of the
widely-used incompressible fluctuating Navier-Stokes equations. In
the low Mach number equations the incompressibility constraint $\nabla\cdot\V v=0$
is replaced by $\nabla\cdot\V v=-\beta\left(Dc/Dt\right)$, which
ensures that compositional changes are accompanied by density changes
in agreement with the fluid equation of state (EOS) at constant pressure
and temperature. This seemingly simple generalization poses many non-trivial
analytical and numerical challenges, some of which we addressed in
this paper.

At the analytical level the low Mach number fluctuating equations
are different from the incompressible equations because the velocity
divergence is directly coupled to the time derivative of the concentration
fluctuations. This means that at thermodynamic equilibrium the velocity
is not only white in space, a well-known difficulty with the standard
equations of fluctuating hydrodynamics, but is also white in time,
adding a novel type of difficulty that has not heretofore been recognized.
The unphysically fast fluctuations in velocity are caused by the unphysical
assumption of infinite separation of time scales between the sound
and the diffusive modes. This unphysical assumption also underlies
the incompressible fluctuating Navier-Stokes equations, however, in
the incompressible limit $\beta\rightarrow0$ the problem is not apparent
because the component of velocity that is white in time disappears.
Here we analyzed the low Mach equations at the linearized level, and
showed that they reproduce the slow diffusive fluctuations in the
full compressible equations, while eliminating the fast pressure fluctuations.
At the formal level, we suggest that a generalized Hodge decomposition
can be used to separate the vortical (solenoidal) modes of velocity
as the independently fluctuating variable, coupled with a gauge formulation
used to treat the divergence constraint. Such nonlinear analysis is
deferred for future research, and here we relied on the fact that
the temporal discretization regularizes the short-time dynamics at
time scales faster than the time step size $\D t$.

At the numerical level, the low Mach number equations pose several
distinct challenges. The first challenge is to construct conservative
spatial discretizations in which density is advected in a locally-conservative
manner while still maintaining the equation of state constraint relating
the local densities and composition. We accomplish this here by using
a specially-chosen model EOS that is linear yet still rather versatile
in practice, and by advecting densities using a velocity that obeys
a discrete divergence constraint. We note that for this simplified
case, the system can be modeled using only the concentration to describe
the thermodynamic state. However, for more general low Mach number
models maintaining a full thermodynamic representation of the state
independent of the constraint leads to more robust numerics. As in
incompressible hydrodynamics, enforcing this constraint requires a
Poisson pressure solver that dominates the computational cost of the
algorithm. A second challenge is to construct temporal integrators
that are at least second-order in time. We accomplish this here by
formally introducing an unconstrained gauge formulation of the equations,
while at the same time taking advantage of the gauge degree of freedom
to avoid ever explicitly dealing with the gauge variable. The present
temporal discretizations are purely explicit and are similar in spirit
to an explicit projection method. A third and remaining challenge
is to design efficient temporal integrators that handle momentum diffusion,
the second-fastest physical process, semi-implicitly. This poses well-known
challenges even in the incompressible setting. These challenges were
bypassed in recently-developed temporal integrators for the incompressible
fluctuating Navier-Stokes equations \cite{LLNS_Staggered} by avoiding
the splitting inherent in projection methods. Extending this type
of Stokes-system approach to the low Mach equations will be the subject
of future research.

One of the principal motivations for developing the low Mach number
equations and our numerical implementation was to model recent experiments
on the development of giant concentration fluctuations in the presence
of sharp concentration gradients. We first studied giant fluctuations
in a time-independent or static setting, as observed experimentally
by inducing a constant concentration gradient via a constant applied
temperature gradient. Our simulations show that under conditions employed
in experimental studies of the diffusive mixing of water and glycerol,
it is reasonable to employ the Boussinesq approximation. The results
also indicate that the constant-transport-coefficient approximation
that is commonly used in theoretical calculations is appropriate if
the diffusion coefficient follows a Stokes-Einstein relation, but
should be used with caution in general.

We continued our study of giant concentration fluctuations by simulating
the temporal evolution of a rough diffusive interface during the diffusive
mixing of hard disk fluids. Comparison between computationally-intensive
event-driven molecular dynamics simulations and our hydrodynamic calculations
demonstrated that the low Mach number equations of fluctuating hydrodynamics
provide an accurate coarse-grained model of fluid mixing. Special
care must be exercised, however, in choosing the bare transport coefficients,
especially the concentration diffusion coefficient, as these are renormalized
by the fluctuations and can be strongly grid-dependent \cite{DiffusionRenormalization_PRL,DiffusionRenormalization,StokesLaw}.
Some questions remain about how to define and measure the bare transport
coefficients from microscopic simulations, but we show that simply
comparing particle and hydrodynamic calculations at large scales is
a robust technique.

The strong coupling between velocity fluctuations and diffusive transport
means that deterministic models have limited utility at mesoscopic
scales, and even macroscopic scales in two-dimensions. This implies
that standard fluorescent techniques for measuring diffusion coefficients,
such as fluorescence correlation spectroscopy (FCS) and fluorescence
recovery after photobleaching (FRAP) \cite{FREP_3D}, may not in fact
be measuring material constants but rather geometry-dependent values
\cite{StokesLaw}. Fluctuating hydrodynamic simulations of typical
experimental simulations, however, are still out of reach due to the
very large separation of time scales between mass and momentum diffusion.
Surpassing this limitation requires the development of a semi-implicit
temporal discretization that is stable for large time steps. Furthermore,
it is also necessary to develop novel mathematical models and algorithms
that are not only stable but also accurate in the presence of such
large separation of scales. This is a nontrivial challenge if thermal
fluctuations are to be included consistently, and will be the subject
of future research.
\begin{acknowledgments}
We would like to thank Boyce Griffith and Mingchao Cai for helpful
comments. J. Bell, A. Nonaka and A. Garcia were supported by the DOE
Applied Mathematics Program of the DOE Office of Advanced Scientific
Computing Research under the U.S. Department of Energy under contract
No. DE-AC02-05CH11231. A. Donev was supported in part by the National
Science Foundation under grant DMS-1115341 and the Office of Science
of the U.S. Department of Energy through Early Career award number
DE-SC0008271. T. Fai wishes to acknowledge the support of the DOE
Computational Science Graduate Fellowship, under grant number DE-FG02-97ER25308.
Y. Sun was supported by the National Science Foundation under award
OCI 1047734.
\end{acknowledgments}
\appendix

\section*{Appendix}

\section{\label{sec:LinearizedAnalysis}Linearized Analysis}

As discussed in more depth in Ref. \cite{LLNS_Staggered}, there are
fundamental mathematical difficulties with the interpretation of the
nonlinear equations of fluctuating hydrodynamics due to the roughness
of the fluctuating fields. It should be remembered, however, that
these equations are coarse-grained models with the coarse-graining
length scale set by the size of the hydrodynamic cells used in discretizing
the equations \cite{DiscreteLLNS_Espanol}. The spatial discretization
removes the small length scales from the stochastic forcing and regularizes
the equations. It is important to point out, however, that imposing
such a small-scale regularization\emph{ }(smoothing) of the stochastic
forcing also requires a suitable renormalization of the transport
coefficients \cite{DiffusionRenormalization_I,DiffusionRenormalization_PRL,StokesLaw},
as we discuss in more detail in Section \ref{sec:MixingMD}.

As long as there are sufficiently many molecules per hydrodynamic
cell the fluctuations in the spatially-discrete hydrodynamic variables
will be small and the behavior of the nonlinear equations will closely
follow that of the \emph{linearized} equations of fluctuating hydrodynamics
\cite{LLNS_Staggered}, which can be given a precise meaning \cite{DaPratoBook}.
It is therefore crucial to understand the linearized equations from
a theoretical perspective, and to analyze the behavior of the numerical
schemes in the linearized setting \cite{LLNS_S_k}.

\subsection{\label{sub:CompressibleSpectra}Compressible Equations}

Some of the most important quantities predicted by the fluctuating
hydrodynamics equations are the equilibrium structure factors (static
covariances) of the fluctuating fields. These can be obtained by linearizing
the compressible equations (\ref{LLNS_primitive}) around a uniform
reference state, $\rho=\rho_{0}+\d{\rho}$, $c=c_{0}+\d c$, $\V v=\d{\V v}$,
$P=P_{0}+\d P$ where 
\[
\d P=c_{T}^{2}\left[\left(\d{\rho}\right)-\beta\rho\left(\d c\right)\right],
\]
and then applying a spatial Fourier transform \cite{FluctHydroNonEq_Book,LLNS_S_k}.
Owing to fluctuation-dissipation balance the static structure factors
are independent of the wavevector $\V k$ at thermodynamic equilibrium,
\begin{eqnarray}
S_{\rho,\rho}\left(\V k\right)= & \av{\left(\widehat{\d{\rho}}\right)\left(\widehat{\d{\rho}}\right)^{\star}} & =\frac{\rho_{0}k_{B}T_{0}}{c_{T}^{2}}+\beta^{2}\frac{\rho_{0}k_{B}T_{0}}{\mu_{c}}\nonumber \\
\M S_{\V v,\V v}\left(\V k\right)= & \av{(\widehat{\delta\V v})(\widehat{\d{\V v}})^{\star}} & =\rho_{0}^{-1}k_{B}T_{0}\,\M I\nonumber \\
S_{c,c}\left(\V k\right)= & \av{\left(\widehat{\d c}\right)\left(\widehat{\d c}\right)^{\star}} & =\frac{k_{B}T_{0}}{\rho_{0}\mu_{c}}.\label{eq:S_equilibrium}
\end{eqnarray}
Note that density fluctuations do not vanish even in the incompressible
limit $c_{T}\rightarrow\infty$ unless $\beta=0$. While fluctuations
in $\rho_{1}$ and $\rho_{2}$ are uncorrelated, the fluctuations
in concentration and density are \emph{correlated} even at equilibrium,
\[
S_{c,\rho}=\av{\left(\widehat{\d{\rho}}\right)\left(\widehat{\d c}\right)^{\star}}=\beta\frac{k_{B}T_{0}}{\mu_{c}}=\rho_{0}\beta S_{c,c}.
\]
We will see below that the low Mach equations correctly reproduce
the static covariances of density and concentration in the limit $c_{T}\rightarrow\infty$.

The dynamics of the equilibrium fluctuations can also be studied by
applying a Fourier-Laplace transform in time in order to obtain the
dynamic structure factors (equilibrium correlation functions) as a
function of wavenumber $\V k$ and wavefrequency $\omega$ \cite{FluctHydroNonEq_Book,LLNS_S_k}.
It is well-known that the dynamic spectrum of density fluctuations
$S_{\rho,\rho}\left(\V k,\omega\right)$ exhibits three peaks for
a given $\V k$, one central Rayleigh peak at small frequencies (slow
concentration fluctuations), and two symmetric Brillouin peaks centered
around $\omega\approx\pm c_{T}k$. As the fluid becomes less compressible
(i.e., the speed of sound increases), there is an increasing separation
of time-scales between the side and central spectral peaks. As we
will see below, the low Mach equations reproduce the central peaks
in the dynamic structure factors only, eliminating the side peaks
and the associated stiff dynamics.

\subsection{\label{sub:GiantTheory}Low Mach Equations}

We now examine the spatio-temporal correlations of the steady-state
fluctuations in the low Mach number equations (\ref{eq:momentum_eq},\ref{eq:rho1_eq},\ref{eq:div_v_constraint},\ref{eq:rho_eq}).
In order to model the nonequilibrium setting in which giant concentration
fluctuations are observed, we include a constant background concentration
gradient in the equations. Note that a density gradient will accompany
a concentration gradient, and this can introduce some additional terms
in $\V F$ depending on how $\rho\chi$ depends on concentration.
For simplicity, we assume $\rho\chi$ is a constant so that the diffusive
term $\grad\cdot\V F$ in (\ref{eq:rho1_eq}) is simply $\rho\chi\grad^{2}c$.
We also assume the viscosity $\eta$ is spatially constant, to get
the simplified coupled velocity-concentration equations, 
\begin{align}
D_{t}\V v= & -\rho^{-1}\grad\pi+\nu\grad^{2}\V v+\rho^{-1}\left(\grad\cdot\M{\Sigma}\right)+\V g\nonumber \\
D_{t}c= & \chi\grad^{2}c+\rho^{-1}\left(\grad\cdot\M{\Psi}\right)\nonumber \\
\grad\cdot\V v= & -\beta D_{t}c.\label{eq:simpl_primitive_eqs}
\end{align}
where $\nu=\eta/\rho$ and $\rho=\rho(c)$ is given by (\ref{eq:EOS_quasi_incomp}).

We linearize the equations (\ref{eq:simpl_primitive_eqs}) around
a steady state, $c=\bar{c}+\d c$, $\V v=\bar{\V v}+\d{\V v}=\d{\V v}$,
and $\pi=\bar{\pi}+\d{\pi}$, where the reference state is in mechanical
equilibrium, $\bar{\rho}^{-1}\grad\bar{\pi}=\V g$. We denote the
background concentration gradient with $\V h=\grad\bar{c}$. We additionally
assume that the reference state varies very weakly on length scales
of order of the wavelength, an in particular, that $\bar{\rho}$ and
$\bar{c}$ are essentially constant. This allows us to drop the bars
from the notation and employ a \emph{quasi-periodic} or weak-gradient
approximation \cite{GiantFluctuations_Theory,DiffusionRenormalization}.
In the linear approximation, the EOS constraint relates density and
concentration fluctuations, $\d{\rho}=\rho\beta\left(\d c\right)$.
The term $\V v\cdot\grad\V v$ is second order in the fluctuations
and drops out, but the advective term $\V v\cdot\grad c$ leads to
a term $\left(\d{\V v}\right)\cdot\V h$ in the concentration equation.
The forcing term due to gravity becomes $\rho^{-1}\left(\d{\rho}\right)\V g=\beta\left(\d c\right)\V g$.
After a spatial Fourier transform, the linearized form of (\ref{eq:simpl_primitive_eqs})
becomes a collection of stochastic differential equations, one system
of linear additive-noise equations per wavenumber,
\begin{eqnarray}
\partial_{t}\left(\widehat{\delta\V v}\right) & = & -i\rho^{-1}\V k\left(\widehat{\delta\pi}\right)-\nu\, k^{2}\left(\widehat{\delta\V v}\right)+i\rho^{-1}\V k\cdot\widehat{\M{\Sigma}}+\beta\V g\left(\widehat{\delta c}\right)\label{eq:v_t_linearized}\\
\partial_{t}\left(\widehat{\delta c}\right) & = & -\V h\cdot\left(\widehat{\delta\V v}\right)-\chi k^{2}\left(\widehat{\delta c}\right)+i\rho^{-1}\left(\V k\cdot\hat{\M{\Psi}}\right)\label{eq:c_t_linearized}\\
\hat{\V k}\cdot\left(\widehat{\delta\V v}\right) & = & -\beta\left[i\chi k\left(\widehat{\delta c}\right)+\rho^{-1}\left(\hat{\V k}\cdot\hat{\M{\Psi}}\right)\right].\label{eq:div_v_linearized}
\end{eqnarray}
Replacing the right hand side of (\ref{eq:div_v_linearized}) with
zero leads to the incompressible approximation used in Ref. \cite{GiantFluctuations_Theory},
corresponding to the Boussinesq approximation of taking the limit
$\beta\rightarrow0$ while keeping the product $\beta g$ constant.

\subsubsection{Equilibrium Fluctuations}

Let us first compare the dynamics of the equilibrium fluctuations
($\V h=\V 0$) in the low Mach equations with those in the complete
compressible equations. For simplicity of notation we will continue
to use the hat symbol to denote the space-time Fourier transform.

In the wavenumber-frequency $(\V k,\omega)$ Fourier domain, the concentration
fluctuations in the absence of a gradient are obtained from (\ref{eq:c_t_linearized}),
\[
\widehat{\delta c}\left(\V k,\omega\right)=\frac{i\rho^{-1}k}{i\omega+\chi k^{2}}\left(\hat{\V k}\cdot\hat{\M{\Psi}}\right),
\]
which is the same as the compressible equations. The density fluctuations
follow the concentration fluctuations, $\widehat{\delta\rho}=\rho\beta\,\widehat{\delta c}$,
and the dynamic structure factor for density shows the same central
Rayleigh peak as obtained from the isothermal compressible equations
\cite{FluctHydroNonEq_Book},
\[
S_{\rho,\rho}\left(\V k,\omega\right)=\frac{\beta^{2}k^{2}}{\omega^{2}+\chi^{2}k^{4}}\av{\hat{\M{\Psi}}\hat{\M{\Psi}}^{\star}}=\beta^{2}\left(\rho\mu_{c}^{-1}k_{B}T\right)\frac{2\chi k^{2}}{\omega^{2}+\chi^{2}k^{4}},
\]
where we used Eq. (\ref{stoch_flux_covariance}) for the covariance
of $\hat{\M{\Psi}}$. This shows that the low Mach number equations
correctly reproduce the slow fluctuations (small $\omega$) in density
and concentration, while eliminating the side Brillouin peaks associated
with the fast isentropic pressure fluctuations.

The fluctuations in velocity, however, are different between the compressible
and low Mach number equations. Let us first examine the transverse
(solenoidal) component of velocity $\widehat{\delta\V v}_{s}=\widehat{\M{\mathcal{P}}}\widehat{\delta\V v}$,
where $\M{\mathcal{P}}$ is the constant-density orthogonal projection
onto the space of divergence-free velocity fields( $\widehat{\M{\mathcal{P}}}=\M I-k^{-2}(\V k\V k^{\star})$
in Fourier space). Applying the projection operator to the velocity
equation (\ref{eq:v_t_linearized}) shows that the fluctuations of
the solenoidal modes are the same as in the incompressible approximation,
\[
\partial_{t}\left(\widehat{\delta\V v}_{s}\right)=-\nu\, k^{2}\left(\widehat{\delta\V v}_{s}\right)+i\rho^{-1}\V k\cdot\widehat{\M{\mathcal{P}}}\widehat{\M{\Sigma}}+\beta\widehat{\M{\mathcal{P}}}\V g\left(\widehat{\delta c}\right).
\]
The fluctuations of the compressive velocity component $\widehat{\delta v}_{l}=\hat{\V k}\cdot\left(\widehat{\delta\V v}\right)$,
on the other hand, are driven by the stochastic mass flux $\hat{\M{\Psi}}$,
as seen from eq. (\ref{eq:div_v_linearized}) at thermodynamic equilibrium,
\[
\widehat{\delta v}_{l}=\frac{i\omega\beta\rho^{-1}}{i\omega+\chi k^{2}}\left(\hat{\V k}\cdot\hat{\M{\Psi}}\right).
\]
The dynamic structure factor (space-time Fourier spectrum) of the
longitudinal component 
\[
S_{v,v}^{(l)}=\av{\left(\widehat{\delta v}_{l}\right)\left(\widehat{\delta v}_{l}\right)^{\star}}\sim\frac{\beta^{2}\omega^{2}}{\left(\omega^{2}+\chi^{2}k^{4}\right)}
\]
does not decay to zero as $\omega\rightarrow\infty$. This indicates
that the fluctuations of velocity are not only white in space but
also white in time. In the incompressible approximation $\beta\rightarrow0$
so that the longitudinal velocity fluctuations vanish and the static
spectrum of the velocity fluctuations is equal to the projection operator,
$\M S_{\V v,\V v}=\widehat{\M{\mathcal{P}}}$ \cite{LLNS_Staggered}.
In the compressible equations, the dynamic structure factor for the
longitudinal component of velocity decays to zero as $\omega\rightarrow\infty$
because it has two sound (Brillouin) peaks centered around $\omega\approx c_{T}k$,
in addition to the central diffusive (Rayleigh) peak. The low Mach
number equations reproduce the central peak (slow fluctuations) correctly,
replacing the side peaks with a flat spectrum for large $\omega$.
The origin of this unphysical behavior is the unjustified assumption
of infinite separation of time scales between the propagation of sound
and the diffusion of mass, momentum and energy. In reality, the same
molecular motion underlies all of these processes and the incompressible
or the low Mach number equations cannot be expected to reproduce the
correct physical behavior at very short time scales ($\omega\gtrsim c_{T}k$).

\subsubsection{\label{sub:GiantFourier}Nonequilibrium Fluctuations}

If we neglect the term involving $\hat{\M{\Psi}}$ in (\ref{eq:div_v_linearized})
and eliminate the Lagrange multiplier (non-thermodynamic pressure)
$\pi$ using (\ref{eq:div_v_linearized}), we obtain the linearized
velocity equation in Fourier space 
\begin{eqnarray}
\partial_{t}\left(\widehat{\delta\V v}\right) & = & -\nu\, k^{2}\left(\widehat{\delta\V v}\right)+i\rho^{-1}\V k\cdot\widehat{\M{\mathcal{P}}}\widehat{\M{\Sigma}}+\beta\left(\widehat{\delta c}\right)\widehat{\M{\mathcal{P}}}\V g\nonumber \\
 &  & -i\beta\chi\left[\V h\cdot\left(\widehat{\delta\V v}\right)\right]\V k+i\beta\chi\left(\nu-\chi\right)k^{2}\left(\widehat{\delta c}\right)\V k.\label{eq:vhat_t_proj}
\end{eqnarray}
It is straightforward to obtain the steady-state covariances (static
structure factors) in the presence of a concentration gradient from
the linearized system of velocity-concentration equations (\ref{eq:c_t_linearized},\ref{eq:vhat_t_proj})
\cite{LLNS_S_k}. The procedure amounts to solving a linear system
for three covariances (velocity-velocity, concentration-concentration,
and velocity-concentration). These types of calculations are particularly
well-suited for modern computer algebra systems like Maple and can
be carried out for arbitrary wavenumber and background concentration
gradient. We omit the full solution for brevity.

Experiments measure the steady-state spectrum of concentration fluctuations
averaged along the gradient \cite{GiantFluctuations_Cannell,FractalDiffusion_Microgravity},
and we will therefore focus on wavenumbers perpendicular to the gradient,
$\V k\cdot\V h=0$. A straightforward calculation shows that the concentration
fluctuations are enhanced as the square of the applied gradient,
\begin{equation}
S_{c,c}\left(\V k\right)=\av{(\widehat{\delta c})(\widehat{\d c})^{\star}}=\frac{k_{B}T_{0}}{\rho_{0}\mu_{c}}+\frac{\nu k_{B}T}{\rho(\nu+\chi)\left[\left(\nu\chi k_{\perp}^{4}+h_{\parallel}g\beta\right)+\beta^{2}\frac{\chi^{3}\nu}{(\nu+\chi)^{2}}k_{\perp}^{2}h_{\perp}^{2}\right]}\, h_{\parallel}^{2},\label{eq:S_c_c}
\end{equation}
where \textbf{$\perp$} and $\parallel$ denote the perpendicular
and parallel component relative to gravity, respectively. The term
in the denominator involving $h_{\perp}$ comes from the low Mach
number constraint (\ref{eq:div_v_constraint}) and is usually negligible
since the concentration gradient is parallel to gravity or $\chi/\nu\ll1$.
Without this term the result (\ref{eq:S_c_c}) is the same result
as obtained in \cite{GiantFluctuations_Theory}, and shows that fluctuations
at wavenumbers below $k_{\perp}^{4}=h_{\parallel}g\beta/\left(\nu\chi\right)$
are suppressed by gravity, as we study numerically in Section \ref{sec:GiantFluct}.

\section{\label{AppendixFiltering}Spatial Filtering}

In our spatial discretization, we use centered differencing for the
advective terms because this leads to a skew-adjoint discretization
of advection \cite{ConservativeDifferences_Incompressible} that maintains
discrete fluctuation-dissipation balance in the spatially-discretized
stochastic equations \cite{LLNS_S_k,DFDB}. It is well-known that
centered discretizations of advection do not preserve monotonicity
properties of the underlying PDEs in the deterministic setting, unlike
one-sided (upwind) discretizations. Therefore, our spatio-temporal
discretization can lead to unphysical oscillations of the concentration
and density in cases where the cell Peclet number $\mbox{Pe}=\D x\norm{\V v}/\chi$
is large.

In the deterministic setting, $\mbox{Pe}$ can always be decreased
by reducing $\D x$ and resolving the fine scale dissipative features
of the flow. However, in the stochastic setting, the magnitude of
the fluctuating velocities at equilibrium is
\[
\av{\left(\d v\right)^{2}}\sim\frac{k_{B}T}{\rho\D V},
\]
where $\D V$ is the volume of the hydrodynamic cell. Therefore, in
two dimensions the characteristic advection velocity magnitude is
$\norm{\V v}\sim\D x^{-1}$. This means that in two dimensions $\mbox{Pe}$
is independent of the grid size and reducing $\D x$ cannot fix problems
that may arise due to a large cell Peclet number. For some of the
simulations reported in Section \ref{sec:MixingMD}, we have found
it necessary to implement a spatial filtering procedure to reduce
the magnitude of the fluctuating velocities while preserving their
spectrum as well as possible at small wavenumbers.

The filtering procedure consists of applying a local averaging operation
to the spatially-discretized random fields $\M{\mathcal{W}}$ and
$\widetilde{\M{\mathcal{W}}}$ independently along each Cartesian
direction. This local averaging smooths the random forcing and thus
reduces the spectrum of the random forcing at larger wavenumbers.
The specific filters we use are taken from Ref. \cite{CompressibleFiltering}.
For stencil width $w_{F}=2$, filtering a discrete field $W$ in one
dimension takes the form
\[
W_{i}\leftarrow\frac{5}{8}W_{i}+\frac{1}{4}\left(W_{i-1}+W_{i+1}\right)-\frac{1}{16}\left(W_{i-2}+W_{i+2}\right).
\]
In Fourier space, for discrete wavenumber $\D k=k\D x$ this local
averaging multiplies the spectrum of $W$ by $\mathcal{F}\left(\D k\right)=1+O\left(\D k^{4}\right)$
and therefore maintains the second-order accuracy of the spatial discretization.
At the same time, the filtering reduces the variance of the fluctuating
fields by about a factor of two in one dimension (a larger factor
in two dimensions). The spectrum of the fluctuations can be preserved
even more accurately if a stencil of width $w_{F}=4$ is used for
the local averaging,
\[
W_{i}\leftarrow\frac{93}{128}W_{i}+\frac{7}{32}\left(W_{i-1}+W_{i+1}\right)-\frac{7}{64}\left(W_{i-2}+W_{i+2}\right)+\frac{1}{32}\left(W_{i-3}+W_{i+3}\right)-\frac{1}{256}\left(W_{i-4}+W_{i+4}\right),
\]
giving a sixth-order accurate filter $\mathcal{F}\left(\D k\right)=1+O\left(\D k^{8}\right)$
and a reduction of the variance by about a third in one dimension.
In two and three dimensions the filtering operators are simple tensor
products of one-dimensional filtering operators. Note that we only
use these filters with periodic boundary conditions. One can, of course,
also use Fourier transform techniques to filter out high frequency
components from the stochastic mass and momentum fluxes.

\section{\label{sec:TransportMD}Extracting Transport Properties from Molecular
Dynamics}

The hydrodynamic simulations described in Section \ref{sec:MixingMD}
require as input transport coefficients, notably, the shear viscosity
$\eta$ and diffusion coefficient $\chi$, which need to be extracted
from the underlying microscopic (molecular) dynamics. This is a very
delicate and important step that has not, to our knowledge, been carefully
performed in previous studies. In this Appendix we give details about
the procedure we developed for this purpose.

\subsection{\label{sub:Viscosity}Viscosity $\nu$}

As discussed in more detail in Refs. \cite{DiffusionRenormalization,StokesLaw},
the transport coefficients in fluctuating hydrodynamics are not universal
material constants but rather depend on the spatial scale (degree
of coarse-graining) under question. We emphasize that this scale-dependent
renormalization is not a molecular scale effect but rather an effect
arising out of hydrodynamic fluctuations, and persists even at the
hydrodynamic scales we are examining here. The best way to define
and measure transport coefficients is by examining the dynamics of
\emph{equilibrium} fluctuations, specifically, by examining the \emph{dynamic
structure factors} of the hydrodynamic fields \cite{FluctHydroNonEq_Book},
i.e., the equilibrium averages of the spatio-temporal Fourier spectra
of the fluctuating hydrodynamic fields. For a hydrodynamic variable
$\xi$ that is transported by a purely diffusive process, the spectrum
of the fluctuations at a given wavenumber $k$ and wavefrequency $\omega$
is expected to be a Lorentzian peak of the form
\[
S_{x}\left(k,\omega\right)=\av{\hat{x}\left(k,\omega\right)\hat{x}^{\star}\left(k,\omega\right)}\sim\left[\omega^{2}+\zeta^{2}k^{4}\right]^{-1},
\]
where in general the diffusion constant $\zeta\left(k\right)$ depends
on the the wavenumber $k$ (wavelength $\lambda=2\pi/k$). We can
therefore estimate the diffusion coefficient $\chi$ by fitting a
Lorentzian peak to $S_{c}\left(k,\omega\right)$ for different $k$'s
(i.e., $\xi\equiv c$). Similarly, we can estimate the kinematic viscosity
$\nu=\eta/\rho$ by fitting a Lorentzian curve to dynamic structure
factors for the scaled vorticity, $\xi\equiv k^{-1}\left(\grad\times\V v\right)_{z}$.

\begin{figure}
\begin{centering}
\includegraphics[width=0.75\textwidth]{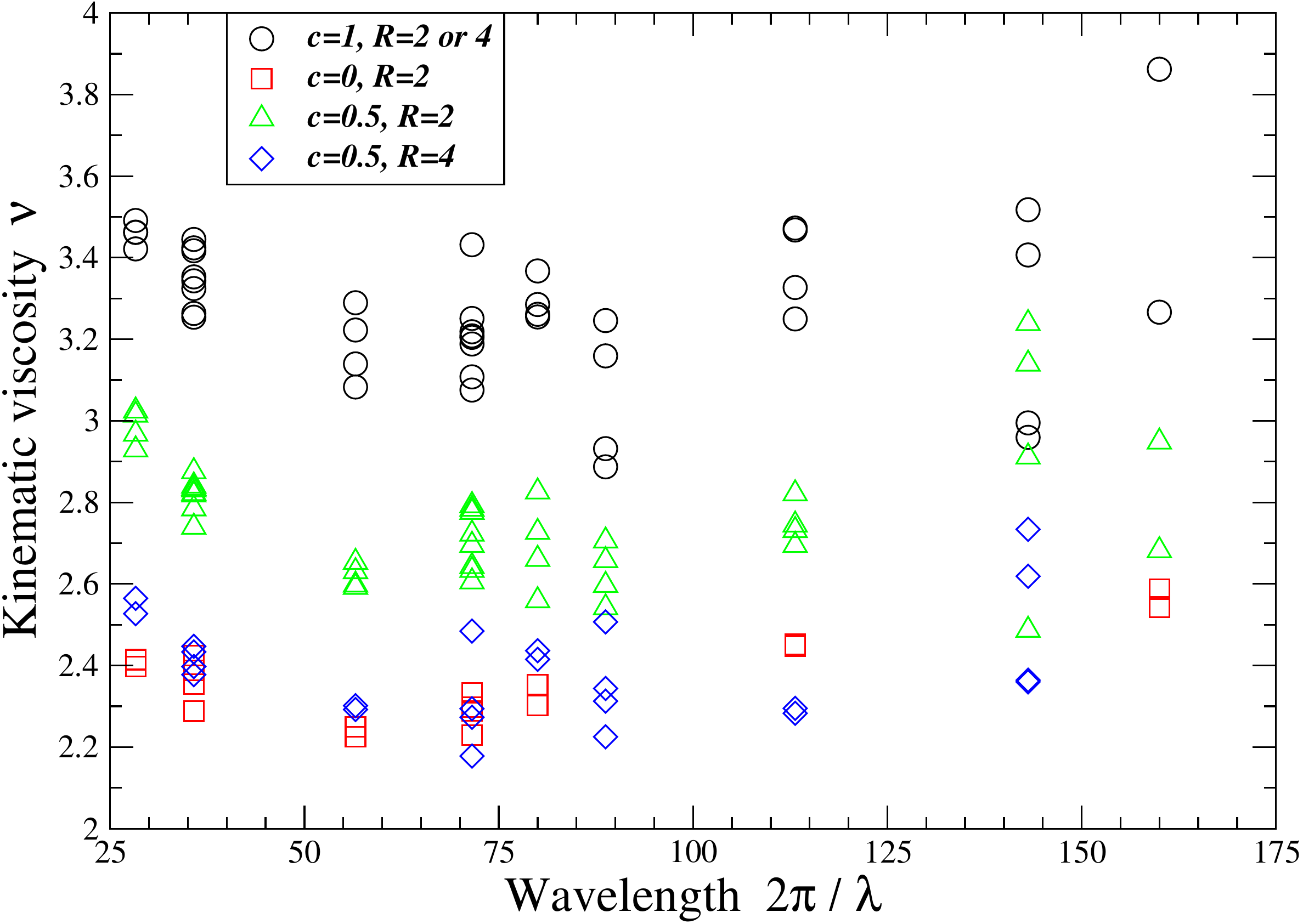}
\par\end{centering}

\caption{\label{fig:S_kw_HDMD}Estimates of the momentum diffusion coefficient
(viscosity) $\nu=\eta/\rho$ obtained from the width of the central
peak in the dynamic structure factor of vorticity. A collection of
$24$ distinct discrete wavenumbers $\V k$ were used and the width
of the peaks estimated using a nonlinear least squares Lorentzian
fit.}
\end{figure}

We performed long equilibrium molecular dynamics simulations of systems
corresponding to a grid of $N_{c}=32$ hydrodynamic cells, and then
calculated the discrete spatio-temporal Fourier spectrum of the hydrodynamic
fields at a collection of discrete wavenumbers $\V k$. Since these
simulations are at equilibrium, the systems are well-mixed, specifically,
the initial configurations were generated by randomly assigning a
species label to each particle. We then performed a nonlinear least
squares Lorentzian fit in $\omega$ for each $\V k$ and estimated
the width of the Lorentzian peak. The results for the dynamics of
the equilibrium vorticity fluctuations are shown in Fig. \ref{fig:S_kw_HDMD}.
We see that kinematic viscosity is relatively constant for a broad
range of wavelengths, consistent with fluctuating hydrodynamics calculations
\cite{Renormalization_Viscosity} and previous molecular dynamics
simulations \cite{Hard_Disks_Transport}. For the pure component one
fluid, $c=1$, with density $\rho\approx0.764$ the figure shows $\nu\approx3.3$.
We therefore used $\eta_{1}\approx0.764\cdot3.3\approx2.5$ in all
of the hydrodynamic runs reported in Section \ref{sec:MixingMD}.
This is about $20\%$ higher than the prediction of the simple Enskog
kinetic theory \cite{Enskog_2DHS}, $\eta\approx2.06$, and is consistent
with the estimates reported in Ref. \cite{Hard_Disks_Transport}.
Because of the diffusion coefficient is small at the densities we
study, more specifically, because the Schmidt number $S_{c}=\nu/\chi$
is larger than 10, we were unable to obtain reliable estimates for
$\chi\left(k\right)$ from the dynamic structure factor for concentration.

Simple dimensional analysis or kinetic theory shows that $\eta\sim\sqrt{m}$.
Since the disks of the two species have equal diameters the viscosity
of the pure second fluid component is
\begin{equation}
\eta_{2}=\eta_{1}\sqrt{\frac{m_{2}}{m_{1}}}=\eta_{1}\sqrt{R}.\label{eq:eta_R_scaling}
\end{equation}
There is no simple theory that accurately predicts the concentration
dependence of the viscosity of a hard disk mixture at higher densities
\cite{HardSphereTransport_Review}. To our knowledge there is no published
Enskog kinetic theory calculations for hard-disk mixtures in two dimensions,
even for the simpler case of equal diameters. As an approximation
to the true dependence, we employed a simple linear interpolation
of the \emph{kinematic} viscosity $\nu(c)=\eta(c)/\rho$ as a function
of the mass concentration $c$ between the two known values $\nu_{1}=\nu\left(c=1\right)\approx3.3$
and $\nu_{2}=\nu\left(c=0\right)=\nu_{1}/\sqrt{R}$. The numerical
results for mixtures with mass ratios $R=2$ and $R=4$ in Fig. \ref{fig:S_kw_HDMD}
are consistent with this approximation to within the large error bars.
For example, for $c=1/2$ and $R=4$ the interpolation gives $\nu=3\cdot3.3/4\approx2.5$
which is in reasonable agreement with the numerical estimate.

\subsection{\label{sub:DiffusionCoeffient}Diffusion Coefficient $\chi$}

For the inter-species diffusion coefficient $\chi$, which we emphasize
is distinct from the self-diffusion coefficients for particles of
either species, Enskog kinetic theory predicts no concentration dependence
and a simple scaling with the mass ratio \cite{HardSphereTransport_Review},
\begin{equation}
\chi\left(R\right)=\chi\left(R=1\right)\sqrt{\frac{1+R}{2R}}.\label{eq:chi_R_scaling}
\end{equation}
This particular dependence on mass ratio $R$ comes from the fact
that the average relative speed between particles of different species
is $\sim\sqrt{k_{B}T/m_{R}}$, where $m_{R}=2m_{1}m_{2}/\left(m_{1}+m_{2}\right)$
is the reduced molecular mass. We have assumed in our hydrodynamic
calculations that the diffusion coefficient is independent of concentration
and follows (\ref{eq:chi_R_scaling}). The only input to the hydrodynamic
calculation is the bare self-diffusion coefficient for the pure component
fluid, $\chi_{0}\left(R=1\right)$. Diffusion is strongly renormalized
by thermal fluctuations, and fluctuating hydrodynamics theory and
simulations predict a strong dependence of the diffusion coefficient
$\chi$ on the wavelength \cite{DiffusionRenormalization}, consistent
with molecular dynamics results \cite{Hard_Disks_Transport}.

\begin{figure}
\begin{centering}
\includegraphics[width=0.49\textwidth]{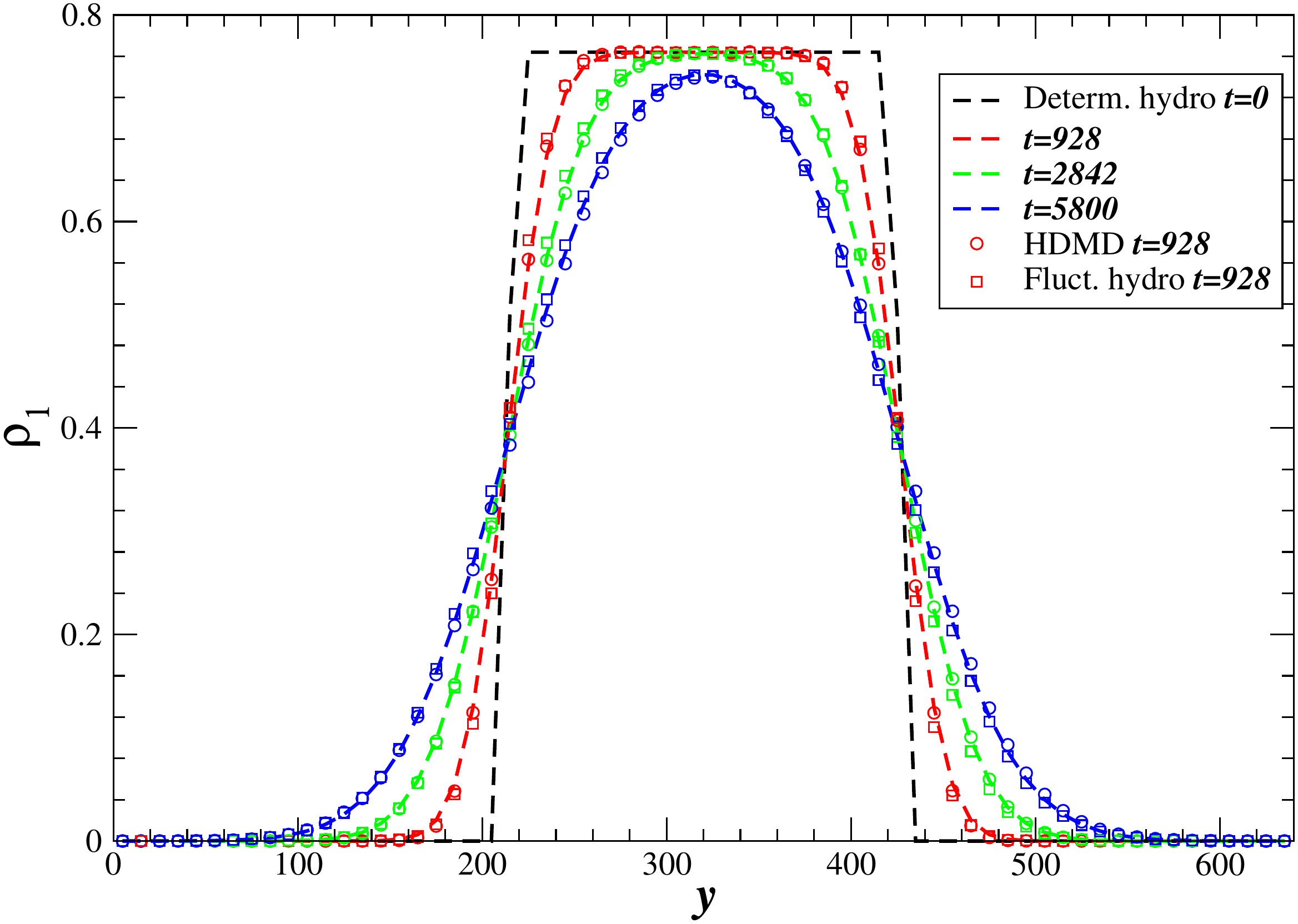}\includegraphics[width=0.49\textwidth]{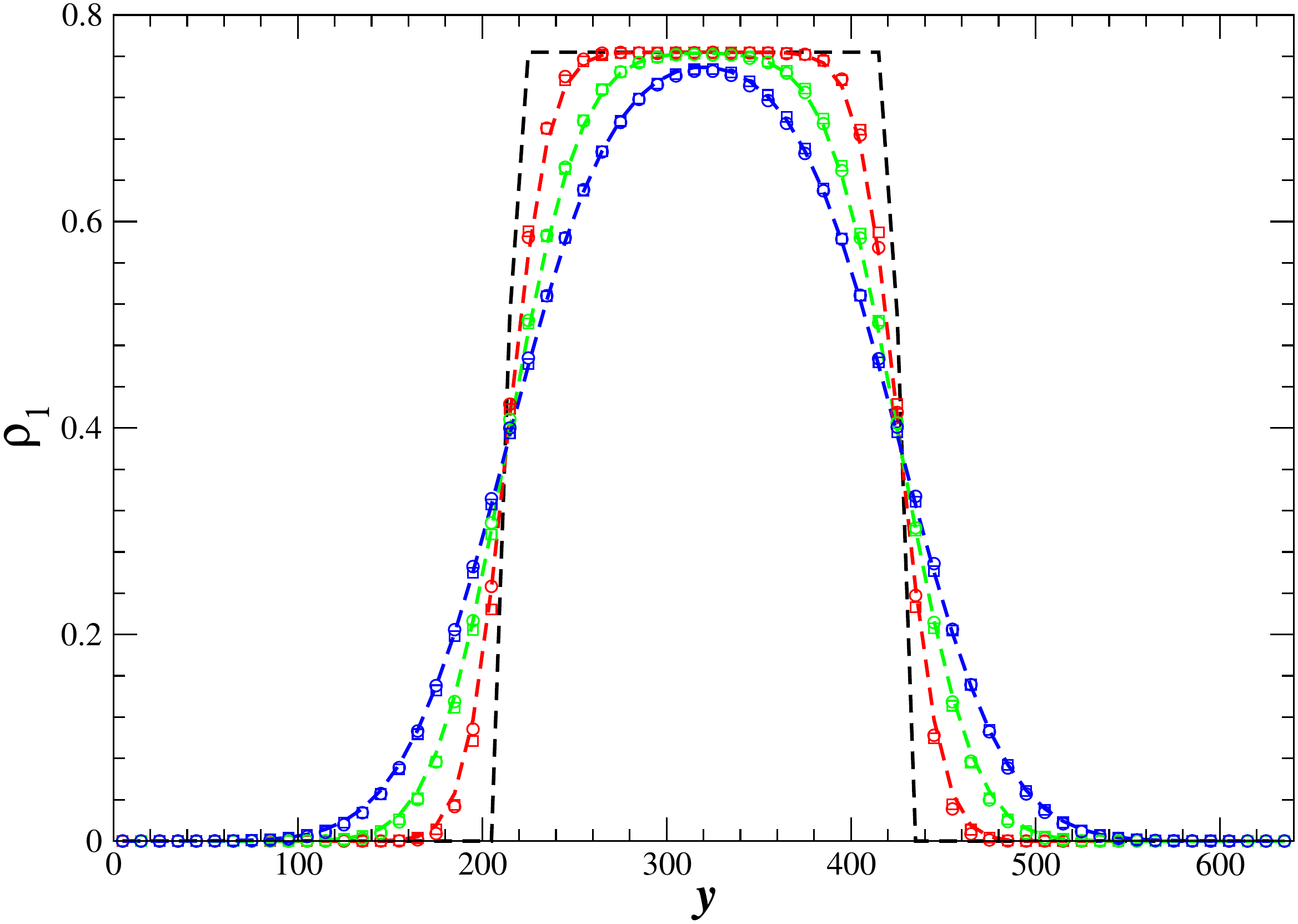}
\par\end{centering}

\caption{\label{fig:spreading_64x64}(Left panel) Diffusive evolution of the
horizontally-averaged density $\rho_{1}^{\left(h\right)}\left(y\right)$
for a system of size $N_{c}=64$ hydrodynamic cells and density ratio
$R=1$, as obtained from HDMD simulations (circles, averaged over
64 runs), deterministic hydrodynamics with $\chi_{\text{eff}}=0.2$
(dashed lines), and fluctuating hydrodynamics with $\chi_{0}=0.09$
(squares, averaged over 64 runs). Error bars are comparable to the
symbol size and not shown for clarity. (Right panel) Same as the left
panel except the density ratio is $R=2$ and the transport coefficients
are adjusted according to (\ref{eq:eta_R_scaling},\ref{eq:chi_R_scaling}).}
\end{figure}

In order to estimate the appropriate value of the bare diffusion coefficient
$\chi_{0}$ we numerically solved an inverse problem. Using simple
bisection, we looked for the value of $\chi_{0}$ that leads to best
agreement for the average or ``macroscopic'' diffusion (mixing)
between the particle and continuum simulations. Specifically, we calculated
the density of the first species $\rho_{1}^{\left(h\right)}\left(y\right)$
along the $y$-direction by averaging $\rho_{1}$ in each horizontal
row of hydrodynamic cells, see Eq. (\ref{eq:rho_1_h}). The results
for $\rho_{1}^{\left(h\right)}$ for mass ratios $R=1$ and $R=4$
are shown in Fig. \ref{fig:spreading_64x64} at different points in
time for systems of size $N_{c}=64$ cells. The figures show the expected
sort of diffusive mixing profile, and is exactly what would be used
in experiments to measure diffusion coefficients using fluorescent
techniques such as Fluorescence Recovery After Photo-bleaching (FRAP)
\cite{FREP_3D}. This macroscopic measurement smooths over the fluctuations
(roughness) of the diffusive interface and only measures an effective
diffusion coefficient at the scale of the domain length $L$. If deterministic
hydrodynamics is employed, $\rho_{1}^{\left(h\right)}\left(y\right)$
is the solution of a one-dimensional system of equations obtained
by simply deleting the stochastic forcing and the $x$-dependence
in the low Mach equations. Instead of solving this system analytically,
we employed our spatio-temporal discretization with fluctuations turned
off, and with an effective diffusion coefficient $\chi=\chi_{\text{eff}}$
that accounts for the renormalization of the diffusion coefficient
by the thermal fluctuations.

By matching the profile $\rho_{1}^{\left(h\right)}\left(y\right)$
between the HDMD and the fluctuating and deterministic hydrodynamic
simulations at mass ratio $R=1$ and system size $N_{c}=64$ cells,
we obtained estimates for the bare $\chi_{0}$ and the renormalized
$\chi_{\text{eff}}$ coefficient (see Fig. \ref{fig:spreading_64x64}).
The best estimate for the bare diffusion coefficient based on this
matching in the absence of filtering is $\chi_{0}=0.09\pm0.01$. This
compares reasonably-well to the prediction of Enskog theory \cite{Enskog_2DHS}
of $\chi\approx0.08$, as well as to the measurement of the self-diffusion
coefficient for a periodic system with 169 disks reported in Ref.
\cite{Hard_Disks_Transport}, $\chi\approx0.14$ (recall that a single
hydrodynamic cell in our case contains about 76 particles). When a
5-point filter is employed the estimate is $\chi_{0}\left(w_{F}=2\right)\approx0.12$
and when a 9-point filter is employed $\chi_{0}\left(w_{F}=4\right)\approx0.11$.
The estimated renormalized diffusion coefficient is much larger, $\chi_{\text{eff}}\approx0.20\pm0.01$,
consistent with a rough estimate based on the simple theory presented
in Ref. \cite{DiffusionRenormalization},
\[
\chi_{\text{eff}}\approx\chi_{0}+\frac{k_{B}T}{4\pi\rho\left(\nu+\chi_{0}\right)}\ln\left(\frac{N_{c}}{3}\right)\approx\begin{cases}
0.18 & \mbox{ for }N_{c}=64\\
0.20 & \mbox{ for }N_{c}=128
\end{cases}.
\]
To within statistical accuracy we were not able to detect the increase
in the estimated diffusion coefficients when using the larger systems
of size $N_{c}=128$ cells, however, for $N_{c}=32$ it was clear
that $\chi_{\text{eff}}$ is reduced. 

It is important to emphasize that $\chi_{\text{eff}}$ is not a material
constant but rather depends on the details of the problem in question,
in particular, the system geometry and size and boundary conditions
\cite{StokesLaw}. By contrast, $\chi_{0}$ is a constant for a given
spatial discretization, and one can use the same number for different
scenarios so long as the hydrodynamic cell size and the filter are
kept fixed. Unlike deterministic hydrodynamics, which presents an
incomplete picture of diffusion, fluctuating hydrodynamics correctly
accounts for the important contribution of the thermal velocity fluctuations
and the roughness of the diffusive interface seen in Fig. \ref{fig:MixingIllustrationMass4}.


\begin{thebibliography}{10}

\bibitem{DiffusionRenormalization_PRL}
A.~Donev, A.~L. Garcia, Anton de~la Fuente, and J.~B. Bell.
\newblock {Diffusive Transport by Thermal Velocity Fluctuations}.
\newblock {\em Phys. Rev. Lett.}, 106(20):204501, 2011.

\bibitem{FractalDiffusion_Microgravity}
A.~Vailati, R.~Cerbino, S.~Mazzoni, C.~J. Takacs, D.~S. Cannell, and M.~Giglio.
\newblock {Fractal fronts of diffusion in microgravity}.
\newblock {\em Nature Communications}, 2:290, 2011.

\bibitem{Nanofluidics_Review}
L.~Bocquet and E.~Charlaix.
\newblock {Nanofluidics, from bulk to interfaces}.
\newblock {\em Chemical Society Reviews}, 39(3):1073--1095, 2010.

\bibitem{Nanofluids_Review}
L.~Wang and M.~Quintard.
\newblock Nanofluids of the future.
\newblock {\em Advances in Transport Phenomena}, pages 179--243, 2009.

\bibitem{SIBM_Biomembrane}
A.~Naji, P.~J. Atzberger, and Frank L.~H. Brown.
\newblock Hybrid elastic and discrete-particle approach to biomembrane dynamics
  with application to the mobility of curved integral membrane proteins.
\newblock {\em Phys. Rev. Lett.}, 102(13):138102, 2009.

\bibitem{BrownainMotor_Peskin}
C.S. Peskin, G.M. Odell, and G.F. Oster.
\newblock {Cellular motions and thermal fluctuations: the Brownian ratchet}.
\newblock {\em Biophysical Journal}, 65(1):316--324, 1993.

\bibitem{Nanopore_Fluctuations}
F.~Detcheverry and L.~Bocquet.
\newblock Thermal fluctuations in nanofluidic transport.
\newblock {\em Phys. Rev. Lett.}, 109:024501, 2012.

\bibitem{CapillaryNanowaves}
R.~Delgado-Buscalioni, E.~Chacon, and P.~Tarazona.
\newblock Hydrodynamics of nanoscopic capillary waves.
\newblock {\em Phys. Rev. Lett.}, 101(10):106102, 2008.

\bibitem{StagerredFluct_Inhomogeneous}
B.~Z. Shang, N.~K. Voulgarakis, and J.-W. Chu.
\newblock {Fluctuating hydrodynamics for multiscale simulation of inhomogeneous
  fluids: Mapping all-atom molecular dynamics to capillary waves}.
\newblock {\em J. Chem. Phys.}, 135:044111, 2011.

\bibitem{BreakupNanojets}
M.~Moseler and U.~Landman.
\newblock Formation, stability, and breakup of nanojets.
\newblock {\em Science}, 289(5482):1165--1169, 2000.

\bibitem{DropletSpreading}
B.~Davidovitch, E.~Moro, and H.A. Stone.
\newblock {Spreading of viscous fluid drops on a solid substrate assisted by
  thermal fluctuations}.
\newblock {\em Phys. Rev. letters}, 95(24):244505, 2005.

\bibitem{DropFormationFluctuations}
Y.~Hennequin, D.~G. A.~L. Aarts, J.~H. van~der Wiel, G.~Wegdam, J.~Eggers,
  H.~N.~W. Lekkerkerker, and D.~Bonn.
\newblock Drop formation by thermal fluctuations at an ultralow surface
  tension.
\newblock {\em Phys. Rev. Lett.}, 97(24):244502, 2006.

\bibitem{ParticleMesoscaleHydrodynamics}
H.~Noguchi, N.~Kikuchi, and G.~Gompper.
\newblock {Particle-based mesoscale hydrodynamic techniques}.
\newblock {\em Europhysics Letters}, 78:10005, 2007.

\bibitem{SHSD_PRL}
A.~Donev, A.~L. Garcia, and B.~J. Alder.
\newblock {Stochastic Hard-Sphere Dynamics for Hydrodynamics of Non-Ideal
  Fluids}.
\newblock {\em Phys. Rev. Lett}, 101:075902, 2008.

\bibitem{Landau:Fluid}
L.D. Landau and E.M. Lifshitz.
\newblock {\em Fluid Mechanics}, volume~6 of {\em Course of Theoretical
  Physics}.
\newblock Pergamon Press, Oxford, England, 1959.

\bibitem{FluctHydroNonEq_Book}
J.~M.~O.~De Zarate and J.~V. Sengers.
\newblock {\em {Hydrodynamic fluctuations in fluids and fluid mixtures}}.
\newblock Elsevier Science Ltd, 2006.

\bibitem{OttingerBook}
H.~C. {\"O}ttinger.
\newblock {\em Beyond equilibrium thermodynamics}.
\newblock Wiley Online Library, 2005.

\bibitem{LB_SoftMatter_Review}
B.~D{\"u}nweg and A.J.C. Ladd.
\newblock {Lattice Boltzmann simulations of soft matter systems}.
\newblock {\em Adv. Comp. Sim. for Soft Matter Sciences III}, pages 89--166,
  2009.

\bibitem{LB_OrderParameters}
Sumesh~P. T., Pagonabarraga I., and R.~Adhikari.
\newblock Lattice-boltzmann-langevin simulations of binary mixtures.
\newblock {\em Phys. Rev. E}, 84:046709, 2011.

\bibitem{SELM}
P.~J. Atzberger.
\newblock {Stochastic Eulerian-Lagrangian Methods for Fluid-Structure
  Interactions with Thermal Fluctuations}.
\newblock {\em J. Comp. Phys.}, 230:2821--2837, 2011.

\bibitem{Bell:09}
J.B. Bell, A.~Garcia, and S.~Williams.
\newblock Computational fluctuating fluid dynamics.
\newblock {\em ESAIM: M2AN}, 44(5):1085--1105, 2010.

\bibitem{LLNS_Staggered}
F.~Balboa Usabiaga, J.~B. Bell, R.~Delgado-Buscalioni, A.~Donev, T.~G. Fai,
  B.~E. Griffith, and C.~S. Peskin.
\newblock {Staggered Schemes for Fluctuating Hydrodynamics}.
\newblock {\em SIAM J. Multiscale Modeling and Simulation}, 10(4):1369--1408,
  2012.

\bibitem{IncompressibleLimit_Majda}
S.~Klainerman and A.~Majda.
\newblock Compressible and incompressible fluids.
\newblock {\em Communications on Pure and Applied Mathematics}, 35(5):629--651,
  1982.

\bibitem{ZeroMachCombustion}
A.~Majda and J.~Sethian.
\newblock {The derivation and numerical solution of the equations for zero Mach
  number combustion}.
\newblock {\em Combustion science and technology}, 42(3):185--205, 1985.

\bibitem{Cahn-Hilliard_QuasiIncomp}
J.~Lowengrub and L.~Truskinovsky.
\newblock {Quasi-incompressible Cahn-Hilliard fluids and topological
  transitions}.
\newblock {\em Proceedings of the Royal Society of London A}, 454(1978):2617,
  1998.

\bibitem{ZeroMach_Buoyancy}
R.~G. Rehm and H.~R. Baum.
\newblock The equations of motion for thermally driven buoyant flows.
\newblock {\em N. B. S. J. Res.}, 83:297--308, 1978.

\bibitem{LowMachAdaptive}
R.B. Pember, L.H. Howell, J.B. Bell, P.~Colella, W.Y. Crutchfield, W.A.
  Fiveland, and J.P. Jessee.
\newblock {An adaptive projection method for unsteady, low-Mach number
  combustion}.
\newblock {\em Combustion Science and Technology}, 140(1-6):123--168, 1998.

\bibitem{ZeroMach_Klein}
T.~Schneider, N.~Botta, KJ~Geratz, and R.~Klein.
\newblock {Extension of finite volume compressible flow solvers to
  multi-dimensional, variable density zero Mach number flows}.
\newblock {\em J. Comp. Phys.}, 155(2):248--286, 1999.

\bibitem{LaminarFlowChemistry}
M.~S. Day and J.~B. Bell.
\newblock {Numerical simulation of laminar reacting flows with complex
  chemistry}.
\newblock {\em Combustion Theory and Modelling}, 4(4):535--556, 2000.

\bibitem{LowMach_FiniteDifference}
F.~Nicoud.
\newblock {Conservative high-order finite-difference schemes for low-Mach
  number flows}.
\newblock {\em J. Comp. Phys.}, 158(1):71--97, 2000.

\bibitem{LowMachAcoustics}
B.~Muller.
\newblock {Low-Mach-number asymptotics of the Navier-Stokes equations}.
\newblock {\em Journal of Engineering Mathematics}, 34(1):97--109, 1998.

\bibitem{StokesKrylov}
M.~Cai, A.~J. Nonaka, J.~B. Bell, B.~E. Griffith, and A.~Donev.
\newblock {Efficient Variable-Coefficient Finite-Volume Stokes Solvers}.
\newblock Arxiv preprint 1308.4605, 2013.

\bibitem{DSMC_Fluctuations_Shear}
A.~L. Garcia, M.~Malek Mansour, G.~C. Lie, M.~Mareschal, and E.~Clementi.
\newblock {Hydrodynamic fluctuations in a dilute gas under shear}.
\newblock {\em Phys. Rev. A}, 36(9):4348--4355, 1987.

\bibitem{Mareschal:92}
M.~{Mareschal}, M.M. {Mansour}, G.~{Sonnino}, and E.~{Kestemont}.
\newblock {Dynamic structure factor in a nonequilibrium fluid: A
  molecular-dynamics approach}.
\newblock {\em Phys. Rev. A}, 45:7180--7183, May 1992.

\bibitem{LongRangeCorrelations_MD}
J.~R. Dorfman, T.~R. Kirkpatrick, and J.~V. Sengers.
\newblock {Generic long-range correlations in molecular fluids}.
\newblock {\em Annual Review of Physical Chemistry}, 45(1):213--239, 1994.

\bibitem{Zarate:04}
J.M.~Ortiz de~Zarate and J.V. Sengers.
\newblock On the physical origin of long-ranged fluctuations in fluids in
  thermal nonequilibrium states.
\newblock {\em J. Stat. Phys.}, 115:1341--59, 2004.

\bibitem{GiantFluctuations_Theory}
A.~Vailati and M.~Giglio.
\newblock {Nonequilibrium fluctuations in time-dependent diffusion processes}.
\newblock {\em Phys. Rev. E}, 58(4):4361--4371, 1998.

\bibitem{TemperatureGradient_Cannell}
C.~J. Takacs, G.~Nikolaenko, and D.~S. Cannell.
\newblock Dynamics of long-wavelength fluctuations in a fluid layer heated from
  above.
\newblock {\em Phys. Rev. Lett.}, 100(23):234502, 2008.

\bibitem{GiantFluctuations_ThinFilms}
D.~{Brogioli}.
\newblock {Giant fluctuations in diffusion in freely-suspended liquid films}.
\newblock {\em ArXiv e-prints}, 2011.

\bibitem{GiantFluctuations_Nature}
A.~Vailati and M.~Giglio.
\newblock {Giant fluctuations in a free diffusion process}.
\newblock {\em Nature}, 390(6657):262--265, 1997.

\bibitem{GiantFluctuations_Cannell}
F.~Croccolo, D.~Brogioli, A.~Vailati, M.~Giglio, and D.~S. Cannell.
\newblock Nondiffusive decay of gradient-driven fluctuations in a
  free-diffusion process.
\newblock {\em Phys. Rev. E}, 76(4):041112, 2007.

\bibitem{LLNS_S_k}
A.~Donev, E.~Vanden-Eijnden, A.~L. Garcia, and J.~B. Bell.
\newblock {On the Accuracy of Explicit Finite-Volume Schemes for Fluctuating
  Hydrodynamics}.
\newblock {\em CAMCOS}, 5(2):149--197, 2010.

\bibitem{DiscreteLLNS_Espanol}
P.~Espa{\~n}ol, J.G. Anero, and I.~Z{\'u}{\~n}iga.
\newblock {Microscopic derivation of discrete hydrodynamics}.
\newblock {\em J. Chem. Phys.}, 131:244117, 2009.

\bibitem{DiffusionJSTAT}
A.~Donev, T.~G. Fai, and E.~Vanden-Eijnden.
\newblock A reversible mesoscopic model of diffusion in liquids: from giant
  fluctuations to fick's law.
\newblock {\em Journal of Statistical Mechanics: Theory and Experiment},
  2014(4):P04004, 2014.

\bibitem{DFDB}
S.~Delong, B.~E. Griffith, E.~Vanden-Eijnden, and A.~Donev.
\newblock {Temporal Integrators for Fluctuating Hydrodynamics}.
\newblock {\em Phys. Rev. E}, 87(3):033302, 2013.

\bibitem{AdiabaticElimination_1}
C.~W. Gardiner and M.~L. Steyn-Ross.
\newblock {Adiabatic elimination in stochastic systems. I-III}.
\newblock {\em Phys. Rev. A}, 29:2814--2844, 1984.

\bibitem{MoriZwanzig_ConstrainedMD}
C.~Hij{\'o}n, P.~Espa{\~n}ol, E.~Vanden-Eijnden, and R.~Delgado-Buscalioni.
\newblock Mori-zwanzig formalism as a practical computational tool.
\newblock {\em Faraday Discuss.}, 144:301--322, 2009.

\bibitem{Landau:StatPhys1}
L.D. Landau and E.M. Lifshitz.
\newblock {\em Statistical Physics}, volume~5 of {\em Course of Theoretical
  Physics}.
\newblock Pergamon, third ed., part 1 edition, 1980.

\bibitem{AnelasticApproximation}
D.R. Durran.
\newblock Improving the anelastic approximation.
\newblock {\em J. Atmos. Sci}, 46(11):1453--1461, 1989.

\bibitem{GaugeIncompressible_E}
W.~E and J.G. Liu.
\newblock Gauge method for viscous incompressible flows.
\newblock {\em Commun. Math. Sci.}, 1(2):317--332, 2003.

\bibitem{almgren-iamr}
A.~S. Almgren, J.~B. Bell, P.~Colella, L.~H. Howell, and M.~L. Welcome.
\newblock A conservative adaptive projection method for the variable density
  incompressible {N}avier-{S}tokes equations.
\newblock {\em J. Comput. Phys.}, 142:1--46, May 1998.

\bibitem{HarWel65}
F.H. Harlow and J.E. Welch.
\newblock Numerical calculation of time-dependent viscous incompressible flow
  of fluids with free surfaces.
\newblock {\em Physics of Fluids}, 8:2182--2189, 1965.

\bibitem{Projection4thOrder_SDC}
Ann~S Almgren, AJ~Aspden, John~B Bell, and ML~Minion.
\newblock On the use of higher-order projection methods for incompressible
  turbulent flow.
\newblock {\em SIAM Journal on Scientific Computing}, 35(1):B25--B42, 2013.

\bibitem{LowMach_DiscreteCompatibility}
P.~Rauwoens, J.~Vierendeels, E.~Dick, and B.~Merci.
\newblock {A conservative discrete compatibility-constraint low-Mach
  pressure-correction algorithm for time-accurate simulations of variable
  density flows}.
\newblock {\em J. Comp. Phys.}, 228(13):4714--4744, 2009.

\bibitem{ConservativeDifferences_Incompressible}
Y.~Morinishi, T.S. Lund, O.V. Vasilyev, and P.~Moin.
\newblock {Fully conservative higher order finite difference schemes for
  incompressible flow}.
\newblock {\em J. Comp. Phys.}, 143(1):90--124, 1998.

\bibitem{bellColellaGlaz:1989}
J.~B. Bell, P.~Colella, and H.~M. Glaz.
\newblock A second order projection method for the incompressible
  {N}avier-{S}tokes equations.
\newblock {\em J. Comp. Phys.}, 85(2):257--283, 1989.

\bibitem{GiantFluctuations_Summary}
A~Vailati, R~Cerbino, S~Mazzoni, M~Giglio, C~J Takacs, and D~S Cannell.
\newblock {Gradient-driven fluctuations in microgravity.}
\newblock {\em Journal of physics. Condensed matter}, 24(28):284134, 2012.

\bibitem{WaterGlycerolDiffusion}
Gerardino D'Errico, Ornella Ortona, Fabio Capuano, and Vincenzo Vitagliano.
\newblock Diffusion coefficients for the binary system glycerol+ water at 25 c.
  a velocity correlation study.
\newblock {\em Journal of Chemical \& Engineering Data}, 49(6):1665--1670,
  2004.

\bibitem{DiffusionRenormalization}
A.~Donev, A.~L. Garcia, Anton de~la Fuente, and J.~B. Bell.
\newblock {Enhancement of Diffusive Transport by Nonequilibrium Thermal
  Fluctuations}.
\newblock {\em J. of Statistical Mechanics: Theory and Experiment},
  2011:P06014, 2011.

\bibitem{Nanopore_FluctuationsPRE}
Fran{\c{c}}ois Detcheverry and Lyd{\'e}ric Bocquet.
\newblock Thermal fluctuations of hydrodynamic flows in nanochannels.
\newblock {\em Physical Review E}, 88(1):012106, 2013.

\bibitem{StokesLaw}
A.~Donev, T.~G. Fai, and E.~Vanden-Eijnden.
\newblock {Reversible Diffusion by Thermal Fluctuations}.
\newblock Arxiv preprint 1306.3158, 2013.

\bibitem{SoretDiffusion_Croccolo}
F.~Croccolo, H.~Bataller, and F.~Scheffold.
\newblock {A light scattering study of non equilibrium fluctuations in liquid
  mixtures to measure the Soret and mass diffusion coefficient}.
\newblock {\em J. Chem. Phys.}, 137:234202, 2012.

\bibitem{MRJ_HS_4D}
M.~Skoge, A.~Donev, F.~H. Stillinger, and S.~Torquato.
\newblock {Packing Hyperspheres in High-Dimensional Euclidean Spaces}.
\newblock {\em Phys. Rev. E}, 74:041127, 2006.

\bibitem{Event_Driven_HE}
A.~Donev, S.~Torquato, and F.~H. Stillinger.
\newblock {Neighbor List Collision-Driven Molecular Dynamics Simulation for
  Nonspherical Particles: {I.} Algorithmic Details {II.} Applications to
  Ellipses and Ellipsoids}.
\newblock {\em J. Comp. Phys.}, 202(2):737--764, 765--793, 2005.

\bibitem{StagerredFluctHydro}
N.~K. Voulgarakis and J.-W. Chu.
\newblock {Bridging fluctuating hydrodynamics and molecular dynamics
  simulations of fluids}.
\newblock {\em J. Chem. Phys.}, 130(13):134111, 2009.

\bibitem{DiscreteDiffusion_Espanol}
P.~Espa{\~n}ol and I.~Z{\'u}{\~n}iga.
\newblock On the definition of discrete hydrodynamic variables.
\newblock {\em J. Chem. Phys}, 131:164106, 2009.

\bibitem{UnbiasedEstimates_Garcia}
A.~L. Garcia.
\newblock {Estimating hydrodynamic quantities in the presence of microscopic
  fluctuations}.
\newblock {\em Communications in Applied Mathematics and Computational
  Science}, 1:53--78, 2006.

\bibitem{HardSphereTransport_Review}
C.~M. Silva and H.~Liu.
\newblock {\em Theory and Simulation of Hard-Sphere Fluids and Related
  Systems}, chapter Modelling of Transport Properties of Hard Sphere Fluids and
  Related Systems, and its Applications, pages 383--492.
\newblock Springer, 2008.

\bibitem{OpenMD_Rafa}
R.~Delgado-Buscalioni.
\newblock Tools for multiscale simulation of liquids using open molecular
  dynamics.
\newblock {\em Numerical Analysis of Multiscale Computations}, pages 145--166,
  2012.

\bibitem{FREP_3D}
K.~Braeckmans, L.~Peeters, N.N. Sanders, S.C. De~Smedt, and J.~Demeester.
\newblock Three-dimensional fluorescence recovery after photobleaching with the
  confocal scanning laser microscope.
\newblock {\em Biophysical journal}, 85(4):2240, 2003.

\bibitem{DiffusionRenormalization_I}
D.~Bedeaux and P.~Mazur.
\newblock {Renormalization of the diffusion coefficient in a fluctuating fluid
  I}.
\newblock {\em Physica}, 73:431--458, 1974.

\bibitem{DaPratoBook}
G.~Da Prato.
\newblock {\em {Kolmogorov equations for stochastic PDEs}}.
\newblock Birkhauser, 2004.

\bibitem{CompressibleFiltering}
C.~A. Kennedy and M.~H. Carpenter.
\newblock {Several new numerical methods for compressible shear-layer
  simulations}.
\newblock {\em Applied Numerical Mathematics}, 14(4):397--433, 1994.

\bibitem{Renormalization_Viscosity}
J.~Lutsko and J.~W. Dufty.
\newblock Mode-coupling contributions to the nonlinear shear viscosity.
\newblock {\em Phys. Rev. A}, 32:1229--1231, 1985.

\bibitem{Hard_Disks_Transport}
R.~Garc\'\i{}a-Rojo, S.~Luding, and J.~J. Brey.
\newblock Transport coefficients for dense hard-disk systems.
\newblock {\em Phys. Rev. E}, 74(6):061305, 2006.

\bibitem{Enskog_2DHS}
David~M. Gass.
\newblock Enskog theory for a rigid disk fluid.
\newblock {\em J. Chem. Phys.}, 54(5):1898--1902, 1971.

\end{thebibliography}

\end{document}